\documentclass[11pt]{article}
\usepackage{amsmath,amssymb}
\usepackage{graphicx}

\newcommand{\qed}{{\unskip\nobreak\hfil\penalty50\hskip2em\vadjust{}
       \nobreak\hfil$\Box$\parfillskip=0pt\finalhyphendemerits=0\par}}

\newtheorem{theorem}{Theorem}[section]
\newtheorem{lemma}{Lemma}[section]
\newtheorem{definition}{Definition}[section]
\newtheorem{cor}{Corollary}
\newtheorem{prop}{Proposition}

\newcommand{{\Z}}{\mathbb Z}
\newcommand{\R}{\mathbb R}
\newcommand{\N}{\mathbb N}

\newcommand{\dist}{\vert \vert}
\renewcommand{\P} {{\mathcal P}}

\newcommand{\cir}{\Gamma_0}
\newcommand{\reg}{{\rm RG}\big( \Gamma_0 \big)}

\newcommand{\vcir}{V\big( \Gamma_0 \big)}

\newcommand{\acon}{{\rm AREA}_{\bo{0},n^2}}

\newcommand{\intg}{{\rm INT}\big( \Gamma_0 \big)}
\newcommand{\exc}{{\rm EXC}\big( \Gamma_0 \big)}

\newcommand{\globdis}{{\rm GD}\big( \Gamma_0 \big)}
\newcommand{\mc}{\mathcal}
\newcommand{\bo}{\mathbf}
\newcommand{\argu}{{\rm arg}}

\newcommand{\axy}{A_{\bo{x},\bo{y}}}
\newcommand{\axye}{E\big(A_{\bo{x},\bo{y}}\big)}
\newcommand{\axyo}{A_{\bo{x_0},\bo{y_0}}}
\newcommand{\axyoe}{E\big(A_{\bo{x_0},\bo{y_0}}\big)}
\newcommand{\xo}{\bo{x_0}}
\newcommand{\yo}{\bo{y_0}}

\newcommand{\zoz}{\{0,1 \}^{E(\Z^2)}}
\newcommand{\aarg}[2]{A_{\bo{#1},\bo{#2}}}

\newcommand{\ang}{\angle}
\newcommand{\qzero}{q_0}
\newcommand{\clu}[1]{C_{\pi/2 - \qzero}^F \big( #1 \big)}
\newcommand{\clum}[1]{C_{\pi/2 - \qzero}^B \big( #1 \big)} 
\newcommand{\clus}[2]{C_{\pi/2 - \qzero}^F \big( #1 \big) \cup C_{\pi/2 - \qzero}^B \big( #2 \big)} 
 
\newcommand{\club}[1]{C_{\pi/2 - 2\qzero}^F \big( #1 \big)}
\newcommand{\clumb}[1]{C_{\pi/2 - 2\qzero}^B \big( #1 \big)}

\newcommand{\cdkp}{{\rm CRG}_{\delta,K,\phi}^{\bo{x},\bo{y}}}
\newcommand{\maxozrg}{{\rm MAXREG}}
\newcommand{\mar}{\theta_{\rm RG}^{\rm MAX}\big(\cir\big)}

\newcommand{\csec}{W_{\bo{u},\epsilon/2}}
\newcommand{\csece}{E\big( W_{\bo{u},\epsilon/2}  \big)}
\newcommand{\csech}{W_{\bo{u},\epsilon/4}}

\newcommand{\sigmac}{\sigma_{(u,\epsilon/2)}}
\newcommand{\ccup}{W_{\bo{e_2},\pi/3}}
\newcommand{\ccdown}{W_{-\bo{e_2},\pi/3}}

\newcommand{\sroapp}{\xi_{Q,H,(h,0)}}

\newcommand{\sweepm}{{\rm SWEEP}^-}
\newcommand{\sweepp}{{\rm SWEEP}^+}
\newcommand{\search}{{\rm SEARCH}}
\newcommand{\ccone}{c_1}
\newcommand{\cctwo}{C_1}

\newcommand{\cthr}{C_4}
\newcommand{\ccona}{C_6}

\newcommand{\cgac}{2}
\newcommand{\cfpo}{\hat{C}}
\newcommand{\cfp}{\overline{C}}

\newcommand{\conka}{C_{\rm rwm}}

\newcommand{\clemlu}{C'}
\newcommand{\ccthr}{C_3}
\newcommand{\ccfour}{C_5}
\newcommand{\conepsstar}{C_4}
\newcommand{\cwinf}{w_{{\rm min}}}
\newcommand{\cwsup}{w_{{\rm max}}}
\newcommand{\clemgac}{C_{\rm gac}}
\newcommand{\clemkac}{C_2}
\newcommand{\clemkam}{m_0}
\newcommand{\crwm}{C}
\newcommand{\cposen}{c_{\rm be}}
\newcommand{\cpi}{20\pi}
\newcommand{\cgenbig}{C_*}

\newcommand{\perpu}[1]{#1^{\perp}}

\newcommand{\circl}{\overline\Gamma_0}
\newcommand{\marcl}{\theta_{\rm RG}^{\rm MAX}\big(\circl \big)}
\newcommand{\regcl}{{\rm RG}\big( \circl \big)}

\newcommand{\centre}{{\rm cen}}
\newcommand{\tcir}{\tilde\Gamma_0}
\newcommand{\area}[1]{{\rm AREA}_{\bo{0},#1}}
\newcommand{\areacon}{\big\vert \intg \big\vert \geq n^2}

\newcommand{\wulff}{\mathcal{W}_\beta}

\newcommand{\sentier}{\tau}

\def\build#1_#2^#3{\mathrel{ \mathop{\kern 0pt#1}\limits_{#2}^{#3}}}

\begin{document}
\title{Phase separation in random cluster models III: \\ circuit regularity}
\author{Alan Hammond\thanks{Department of Statistics, University of Oxford. This work was undertaken during a visit to the Theory Group at Microscoft Research in Redmond, WA, and at Ecole Normale Superieure in Paris.}} 
 \maketitle

\begin{abstract}
We study the droplet that results from conditioning the subcritical Fortuin-Kasteleyn random cluster model on the presence of an open circuit $\Gamma$ encircling the origin and enclosing an area of at least (or exactly) $n^2$. 
In this paper, we prove that the resulting circuit is highly regular: we define a notion of 
a regeneration site in such a way that, for any such site $\bo{v} \in \Gamma$, the circuit $\Gamma$ cuts through the radial line segment 
through $\bo{v}$ only at $\bo{v}$.
We show that, provided that the conditioned circuit is centred at the origin in a natural sense, the set of regeneration sites reaches into all parts of the circuit, with maximal distance from one such site to the next being at most logarithmic in $n$ with high probability.
The result provides a flexible control on the conditioned circuit that permits the use of surgical techniques to bound its fluctuations, and, as such, it plays a crucial role in the derivation of bounds on the local fluctuation of the circuit carried out in \cite{hammondone} and \cite{hammondtwo}. 
\end{abstract}

\setlength{\baselineskip}{16pt}
%\begin{center}
%{\Large\bf The fluctuation of the droplet boundary in random cluster models}\bs\bs\\
%{\Large\bf Alan Hammond}\bs\bs\\
%\end{center}
\begin{section}{Introduction}
Phase separation refers to the study of the geometry of the random curve that forms the boundary between 
two distinct populations of spins in a statistical mechanical model such as percolation, the Potts model or the random cluster model. For example, if the two-dimensional Ising model at supercritical inverse temperature $\beta > \beta_c$
in a large box with negative boundary conditions is conditioned by the presence of a significant excess of plus signs, then those excess signs typically gather together in a single droplet having the opposite magnetisation to its exterior. The object of study of phase separation is then the droplet boundary. As explained in \cite{alexcube} and \cite{hammondone}, a close relative of this problem is the behaviour of the circuit that arises by conditioning a subcritical random cluster model on the presence of a circuit encircling the origin and trapping a high area.

The aim of this paper is to prove that the circuit that results from this conditioning has a high degree of regularity. The original Ornstein-Zernike picture proposed that such a circuit would macroscopically resemble a dilate of an isoperimetrically extremal curve (the Wulff curve), and, that, microscopically, the circuit would be divided into a large number of small irreducible fragments, the circuit within each fragment following the direction of the macroscopic Wulff curve up to a fixed angular deviation. 

In essence, we establish that the irreducible fragments are of at most a size logarithmic in the area that the circuit is conditioned to trap. Such an understanding of the conditioned circuit is highly convenient for the purpose of analysing the circuit by means of, for example, surgical techniques. Beyond its intrinsic interest, the present work thus provides the foundation for the techniques in \cite{hammondone} 
and \cite{hammondtwo}
that yield strong conclusions regarding the local deviation of the conditioned circuit.

We recall the definition of the random cluster model.
\begin{definition}
For $\Lambda \subseteq \Z^2$, let $E(\Lambda)$ denote the set of nearest-neighbour edges whose endpoints lie in $\Lambda$
and write $\partial_{\rm int} \big( \Lambda \big)$ for the interior vertex boundary of $\Lambda$, namely, the subset of $\Lambda$ each of whose elements is an endpoint of some element of 
 $E(\Lambda)^c$. Fix a choice of $\Lambda \subseteq \Z^2$ that is finite. 
The free random cluster model on $\Lambda$ with parameters $p \in [0,1]$ and $q > 0$ on $\Lambda$ 
is the probability space over $\eta: E(\Lambda) \to \{0,1\}$ with measure
$$
\phi_{p,q}^f(\eta) = \frac{1}{Z_{p,q}}   p^{\sum_e \eta(e)} 
\big( 1 - p \big)^{\sum_e (1 - \eta(e))} q^{k(\eta)},
$$
where $k(\eta)$ denotes the number of connected components in the subgraph of 
$\big(\Lambda,E(\Lambda)\big)$ containing all vertices and all edges $e$ such that $\eta(e) = 1$. (The constant $Z_{p,q}$ is a normalization.) The wired random cluster model 
$\phi_{p,q}^w$ is defined similarly, with $k(\eta)$ now denoting the number of such connected components none of whose edges touch $\partial_{\rm int} \big( \Lambda \big)$.

For parameter choices $p \in [0,1]$ and $q \geq 1$, 
either type of random cluster measure $\P$ satisfies the FKG inequality: suppose that 
$f,g:  \{0,1\}^{E(\Lambda)} \to \R$
are increasing functions with respect to the natural partial order on $\{0,1\}^{E(\Lambda)}$. 
Then $\mathbb{E}_\P \big( fg \big) \geq \mathbb{E}_\P \big( f \big) \mathbb{E}_\P \big( g \big)$,
where $\mathbb{E}_\P$ denotes expectation with respect to $\P$.

Consequently, we define the infinite-volume free and wired random cluster measures $\P^f$ and $\P^w$ as limits of the finite-volume counterparts taken along any increasing sequence of finite sets $\Lambda \uparrow \Z^2$. The measures $\P^f$ and $\P^w$ are defined on the space of functions $\eta: E(\Z^2) \to \{0,1\}$ with the product $\sigma$-algebra. In a realization $\eta$, the edges $e \in E(\Z^2)$
 such that $\eta(e) = 1$ are called open; the remainder are called closed.  A subset of $E(\Z^2)$ will be called open (or closed) if all of its elements are open (or closed).
We will record a realization in the form $\omega \in \zoz$, where 
the set of coordinates that are equal to $1$ under $\omega$ is the set of open edges under $\eta$.
Any $\omega \in \zoz$ will be called a configuration. For $\bo{x},\bo{y} \in \Z^d$, we write $\bo{x} \leftrightarrow \bo{y}$ to indicate that $\bo{x}$ and $\bo{y}$ lie in a common connected component of open edges.  
 
Set $\beta \in (0,\infty)$ according to $p = 1 - \exp\{ - 2 \beta \}$. 
In this way, the infinite volume measures are parameterized by $\P_{\beta,q}^w$ and $\P_{\beta,q}^f$ with $\beta > 0$ and $q \geq 1$.
For any $q \geq 1$, $\P^w_{\beta,q} = \P^f_{\beta,q}$ for all but at most countably many values of $\beta$ \cite{grimmett}. %Writing $\big\{ 0 \leftrightarrow \infty \big\}$ for the event that the connected component of open edges in which the origin lies is infinite, 
We may thus define
$$
\beta_c^1 = \inf \big\{  \beta > 0: \P^*_{\beta,q} \big(  0 \leftrightarrow \infty \big) > 0 \big\}
$$
obtaining the same value whether we choose $* = w$ or $* = f$.
\end{definition}
There is a unique random cluster model for each subcritical $\beta < \beta_c^1$ \cite{grimmett}, that we will denote by $\P_{\beta,q}$.
\begin{definition}
Let $\hat{\beta}_c$ denote the supremum over $\beta > 0$ such that the following holds: letting $\Lambda = \big\{ -N,\ldots,N \big\}^d$, there exist constants $C > c > 0$ such that, for any $N$,
$$
\P^w_{\beta,q} \Big( \bo{0} \leftrightarrow \Z^d \setminus  \Lambda_N \Big) \leq C \exp \big\{ - c N  \big\}.
$$
\end{definition}
In the two-dimensional case that is the subject of this article, 
it has been established 
that $\beta_c^1 = \hat{\beta}_c$ for $q=1$, $q=2$  
and for $q$ sufficiently high, by \cite{alexmix}, and respectively \cite{ab}, \cite{abf} and \cite{lmmrs}.
A recent advance \cite{beffaraduminilcopin} showed that, on the square lattice, in fact, $\beta_c^1 = \hat{\beta}_c$ holds for all $q \geq 1$. The common value, which is $2^{-1} \log \big( 1 + \sqrt{q} \big)$, we will denote by $\beta_c$.

The droplet boundary is now defined:
\begin{definition}\label{defcirsur}
A circuit %\hfff{gencir} 
$\Gamma$ is a nearest-neighbour path in $\Z^2$ whose endpoint coincides with its start point, but for which no other vertex is visited twice.
We set $E(\Gamma)$ equal to the set of nearest-neighbour edges between successive elements of $\Gamma$. 
For notational convenience, when we write $\Gamma$, we refer to the closed subset of $\R^2$ given by the union of the topologically closed intervals corresponding to the elements of $E(\Gamma)$. 
We set $V(\Gamma) = \Gamma \cap \Z^2$.

Let $\omega \in \zoz$. A circuit $\Gamma$ is called open if $E(\Gamma)$ is open. 
For any circuit $\Gamma$, we write %\hfff{intcir} 
${\rm INT} \big( \Gamma \big)$ for the bounded component of $\R^2 \setminus \Gamma$, that is, for the set of points enclosed by $\Gamma$. 

An open circuit $\Gamma$ is called outermost if any open circuit $\Gamma'$ satisfying 
${\rm INT} \big( \Gamma \big) \subseteq {\rm INT} \big( \Gamma' \big)$ is equal to $\Gamma$. Note that, if
a point $\bo{z} \in \R^2$ is enclosed by a positive but finite number of open circuits in a configuration 
$\omega \in \zoz$ , it
is enclosed by a unique outermost open circuit.

We write %\hfff{outcir} 
$\cir$ for the outermost open circuit $\Gamma$ for which $\bo{0} \in {\rm INT}(\Gamma)$, taking $\cir = \emptyset$ if no such circuit exists.
\end{definition}
\noindent{\bf Remark.} 
Under any subcritical random cluster measure $P = \P_{\beta,q}$, with $\beta < \hat{\beta}_c$, there is an exponential decay in distance for the probability that two points lie in the same open cluster. (See Theorem $A$ of \cite{civ}.) As such, $P$-a.s., no point in $\R^2$ is surrounded by infinitely many open circuits, so that $\cir$ exists (and is non-empty) whenever $\bo{0}$ is surrounded by an open circuit.

We will formulate a notion of regeneration site for the conditioned circuit such that the circuit cuts through the semi-infinite line segment from the origin through any regeneration site only at that site. 
We will prove circuit regularity in the form that every subpath of the circuit of logarithmic diameter contains such sites.
Clearly, to formulate such a statement, it is necessary to centre the circuit appropriately, so that the origin is not close to some part of the circuit. We now explain the convention that we adopt. 
\begin{subsubsection}{The Wulff shape and circuit centering}\label{seccircen}
The macroscopic profile of the conditioned circuit is given by the boundary of the Wulff shape.
\begin{definition}
 We define the {\it inverse correlation length}: for $\bo{x} \in \R^2$,
$$
\xi(\bo{x}) = - \lim_{k \to \infty} k^{-1} \log P \big( \bo{0} \leftrightarrow \lfloor k\bo{x} \rfloor  \big),
$$
where $\lfloor \bo{y} \rfloor \in \Z^2$ is the component-wise integer part of $\bo{y} \in \R^2$. 
\end{definition}
\begin{definition}
The unit-area Wulff shape $\wulff$ is the compact set given by
$$
\wulff 
= \lambda \bigcap_{\bo{u} \in S^1} \Big\{ \bo{t} \in \R^2: \big( \bo{t},\bo{u} \big) \leq \xi\big(\bo{u}\big)  \Big\},
$$
with the dilation factor $\lambda > 0$ chosen to ensure that $\big\vert \wulff \big\vert = 1$.
\end{definition}
Global deviations of the conditioned circuit from the Wulff shape may be measured in the following way.
\begin{definition}\label{defglobdis}
Let $\Gamma \subseteq \R^2$ denote a circuit.  Define its
global distortion ${\rm GD} \big( \Gamma \big)$ (from an factor $n$ dilate of the Wulff shape boundary) by means of 
\begin{equation}\label{eqdefgd}
{\rm GD} \big( \Gamma \big) = \inf_{\bo{z} \in \Z^2} d_H \Big(  n \partial \wulff + \bo{z}, \cir  \Big),
\end{equation}
where $d_H$ denotes the Hausdorff distance on sets in $\R^2$. 
\end{definition}
%(In a general context, this would be a peculiar definition. However, we will work with this quantity only in the case of circuits that are conditioned to trap an area of at least, or exactly, $n^2$.)
\begin{definition}\label{defncentr}
Let $\Gamma \subseteq \R^2$ denote a circuit. 
The lattice point $\bo{z}$ attaining the minimum in (\ref{eqdefgd}) will be called the centre $\centre(\Gamma)$ of $\Gamma$, with a fixed rule used to break ties. 
\end{definition}
We work with circuits centred at the origin:
\begin{definition}
We write $\acon$ for the event $\big\{ \areacon \big\} \cap \big\{ \centre(\cir) = \bo{0} \big\}$.
When $\acon$ occurs, we write $\tcir = n \partial \wulff$ for the Wulff shape dilate attaining the infimum that defines $\globdis$.
\end{definition}
\end{subsubsection}
\begin{subsubsection}{The radial definition of regeneration structure}\label{secrrg}
We now define, for a circuit $\Gamma$, the notion of a $\Gamma$-regeneration site: 
\begin{definition}
Let $\qzero > 0$ and $c_0 \in \big(0,\qzero/2 \big)$ denote fixed constants. (The precise conditions that fix these values will be stated in Section \ref{secfp}.) For $\bo{x},\bo{y} \in \R^2$, write $\ang \big( \bo{x},\bo{y} \big) \in [0,\pi]$
for the angle between $\bo{x}$ and $\bo{y}$.  
The forward cone $C^F_{\pi/2 - \qzero} \big( \bo{v} \big)$ denotes the set of vectors $\bo{w} \in \R^2$ for which $\ang\big( \bo{w} - \bo{v}, \bo{v}^{\perp} \big) \leq \pi/2 - \qzero$, where $\bo{v}^{\perp}$ denotes the vector obtained from $\bo{v}$ by a counterclockwise right-angle rotation. The backward cone $C^B_{\pi/2 - \qzero} \big( \bo{v} \big)$
 denotes the set of vectors $\bo{w} \in \R^2$ for which $\ang\big( \bo{w} - \bo{v}, - \bo{v}^{\perp} \big) \leq \pi/2 - \qzero$.

A site $\bo{v} \in V(\Gamma)$ in a circuit $\Gamma$ for which $\centre(\Gamma)$ is called a $\Gamma$-regeneration site if
\begin{equation}\label{eqdefreg}
\Gamma \cap 
 W_{\bo{v},c_0} \subseteq C^F_{\pi/2 - \qzero} \big(  \bo{v} \big) \cup C^B_{\pi/2 - \qzero} \big( \bo{v} \big).
\end{equation}
where, for $c \in [0,\pi)$,
\begin{equation}\label{wnot}
W_{\bo{v},c} = \Big\{ \bo{z} \in \R^2: \argu(\bo{v}) - c \leq  \argu(\bo{z}) \leq  \argu(\bo{v}) + c \Big\}
\cup \big\{ \bo{0} \big\}
\end{equation}
denotes the cone of points whose angular displacement from $\bo{v}$ is at most $c$. 
See Figure 1.
%\ref{figregdef}. 
We write ${\rm RG}\big( \Gamma \big)$ for the set of $\Gamma$-regeneration sites.
\end{definition} 
\begin{figure}\label{figregdef}
\begin{center}
\includegraphics[width=0.3\textwidth]{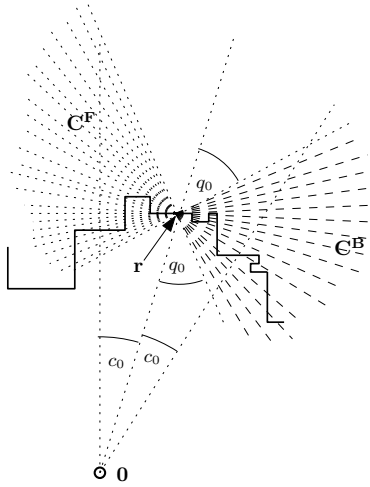} \\
\end{center}
\caption{A $\cir$-regeneration site and the nearby circuit.}
\end{figure}
\end{subsubsection}
\begin{subsection}{Statement of results}
We monitor angular segments in which $\cir$-regeneration sites are absent:
\begin{definition}\label{defmar}
We write $\mar \in [0,2\pi)$ for the angle of the largest angular sector rooted at the origin that contains no $\cir$-regeneration sites. That is,
\begin{equation}\label{eqdefmar}
\mar = \sup \Big\{ r \in [0,2\pi): \exists \bo{a} \in S^1, 
 W_{\bo{a},r/2} \big( \bo{0} \big) \cap \reg = \emptyset  \Big\}.
\end{equation}
\end{definition}
The principal result of this paper 
is the following statement on the regularity of the conditioned circuit:
\begin{theorem}\label{thmmaxrg}
Let $P = \P_{\beta,q}$ with $\beta < \hat{\beta}_c$ and $q \geq 1$. 
There exist $c > 0$ and $C > 0$ such that
$$
P \Big(  \mar > u/n   \Big\vert \acon \Big)
 \leq  \exp \Big\{ - c u \Big\} 
$$
for $C \log n   \leq u \leq c n$.
\end{theorem}
In considering boundary regularity under the conditioning $\areacon$, an alternative is to study the open cluster to which $\cir$ belongs. We denote this set by $\circl$. We define the set $\regcl$ of {\it cluster regeneration sites} 
according to (\ref{eqdefreg}), in which $\cir$ is replaced by $\circl$, and $\marcl$ by (\ref{eqdefmar}) with the same change.
\begin{theorem}\label{thmmaxrgcl}
There exist $c_0,\qzero > 0$, $c > 0$ and $C > 0$ such that the following holds.
Let $P = \P_{\beta,q}$ with $\beta < \hat{\beta}_c$ and $q \geq 1$. 
There exist  $c > 0$ and $C > 0$ such that,
for $C \log n   \leq u \leq c n$,
$$
P \Big(  \marcl > u/n   \Big\vert \areacon \Big)
 \leq  \exp \Big\{ - c u \Big\}.  
$$ 
\end{theorem}
This formulation of the theorem is valuable in \cite{hammondtwo} for proving lower bounds on the local fluctuation of the conditioned circuit. 
It is easy to see that $\regcl \subseteq \reg$, so that Theorem \ref{thmmaxrgcl}
is a strengthening of Theorem \ref{thmmaxrg}. Its proof, however, is the same
as that of Theorem \ref{thmmaxrg}, except for appropriate substitutions of $\circl$ for $\cir$, and with some objects becoming clusters, instead of paths. 
. 
%It is a slight distraction to keep track of the two sets $\cir$ and $\circl$ throughout the proof, so we prefer to present the proof of Theorem \ref{thmmaxrg} and comment when needed on the changes required to obtain Theorem \ref{thmmaxrgcl}. 

Finally, we record the extension of the theorems for fixed area conditioning. 
\begin{cor}
Theorems \ref{thmmaxrg} and \ref{thmmaxrgcl} hold with verbatim statements for the conditional measure
$P \big( \cdot \big\vert \vert {\rm INT}(\cir)\vert = n^2 \big)$. 
\end{cor}
\noindent{\bf Proof}. This is implied by the proof of the analogous statement, Theorem 1.3 of \cite{hammondone}. \qed
\end{subsection}
\begin{subsubsection}{Related definitions and approaches}
Alexander \cite{alexcube} proved that a range of subcritical measures that include the random cluster measures lack super-logarithmically sized bottlenecks, where a bottleneck refers to a section of the conditioned circuit that doubles back to a distance much shorter than its length, in this way, excluding one of the behaviours that may cause a long gap between regeneration sites as we have defined them. 

The theory of local limit analysis and regeneration structure of the connected component arising when a subcritical random cluster model is conditioned to connect two distant points has been developed by \cite{civ}. Indeed, we will make use of some aspects of this theory in our approach. (See Section \ref{secoz}.) Related deductions such as an invariance principle for the phase separation boundary are presented in \cite{greenbergioffe}. The adaptation of Ruelle's theory of spectral operators required for the local limit analysis is given in \cite{civoz}.

The paper \cite{hrynivioffe} studies a model of self-avoiding polygons in the first quadrant in the plane, in which a polygon is penalized exponentially according to its total length, and is then conditioned on enclosing a 
high area in the region between the polygon and the coordinate axes. Deductions are made that are similar to Theorem \ref{thmmaxrg}
for a definition of regeneration site that stipulates that the polygon cut through the vertical (or horizontal) line in question at a unique point. 
The apparatus developed by \cite{hrynivioffe} is a powerful one, achieving a sharp asymptotic formula for the partition function of the model. 
We do not derive such formulae in this paper, but rather obtain a result on circuit regularity 
in which all the complexities of the planar nature of the problem are present, 
for a radially specified notion of regeneration site,
that is valid in the class of subcritical random cluster models. Of course, we are also powerfully motivated to do so by the need for such a result provided by \cite{hammondone} and \cite{hammondtwo}, in which are  derived  logarithmically sharp uniform radial and longitudinal local deviation results for the conditioned circuit. It would be very natural to seek to further understand fluctuation in the conditioned circuit by combining local limit analysis ideas from \cite{hrynivioffe} and the surgical techniques presented in this paper and in \cite{hammondone} and \cite{hammondtwo}.
\end{subsubsection}
\begin{subsection}{Extensions}
As discussed in \cite{hammondone}, the natural scaling for local deviation of the conditioned circuit from its convex envelope is typically $n^{1/3}$ radially and $n^{2/3}$ longitudinally.
The absence of superlogarithmic backtracking known by virtue of Theorem \ref{thmmaxrg} means that, in scaling the curve by these factors in orthogonal directions about a given point, the limiting law of the curve, if it exists, may be parametrized as a function of a single variable. A natural extension of the present paper is to prove the existence of a limiting law and to find it explicitly, for example, as a solution of a stochastic differential equation, by exploiting the established regularity of the circuit.
The analogous question for a Brownian bridge $B:[-T,T] \to [0,\infty)$ conditioned to remain above the semi-circle of radius $T$ centred at $(0,0)$ has been found to satisfy a stochastic differential equation with a drift term expressed in terms of the Airy function \cite{ferrarispohn}. 

Another natural way of measuring deviation in the circuit is from an appropriate dilation of the Wulff shape. 
In this case, deviation has a typical order of $n^{1/2}$.
Central limit theorems for this deviation have been proved, for a one-dimensional random walk under whose trajectory a high area is captured in \cite{dobhryone}, 
and in the phase boundary of a two-dimensional Ising model at low temperature, in \cite{dobhrytwo}. 
\end{subsection}
%\end{subsection}
\begin{subsection}{The structure of the paper}
Section \ref{secnott} introduces notations and some required tools. 

The argument is divided into three parts.
A first step establishes that 
regeneration sites occur in any sector of angle that is uniformly positive in $n$. This result appears as Proposition \ref{proprgmac} in Section \ref{secmac}. 
In the second step, that appears as Proposition \ref{propag} in Section \ref{secagm}, a lower bound is provided on the area-excess captured by $\cir$ over the mandated $n^2$ for the conditional measure $P \big( \cdot \big\vert  \acon \big)$.
These first two steps are preliminaries for the proof of Theorem \ref{thmmaxrg} itself, which appears in the final Section \ref{secproof}.  Each of these three sections begins with an overview of the proof in question.

The role of surgical techniques is critical in all of these arguments. 
Lemmas from \cite{hammondone} are restated here to reduce the need to continually cross-reference with that paper.  We emphasise, however, that \cite{hammondone} functions as an introduction to the techniques used in this paper, and strongly recommend that the reader examine the overview Section 1.1 of \cite{hammondone} before embarking on the proofs in the present article. We also mention that,
although Proposition \ref{proprgmac} is simply a tool whose statement is quite unsurprising,
 its proof is notationally a little cumbersome. 
The heart of the argument in this paper appears in Section \ref{secproof}.
A suggestion for realizing the aim of understanding the principal ideas of this paper is to begin 
by reading the introduction and Section 1.1 of \cite{hammondone}, then briefly peruse the proof of Proposition 2 in Section 3 of \cite{hammondone}, (since this proof uses the relevant surgical technique and forms a template for other arguments), and then the statements and perhaps the sketch proofs 
of Propositions \ref{proprgmac}  and \ref{propag} in the present article, (that are to be found at the beginning of Sections \ref{secmac} and \ref{secagm}), before turning to the actual proof of Theorem \ref{thmmaxrg} in  Section \ref{secproof}. \\
\noindent{\bf Acknowledgments.} I would like to thank Kenneth Alexander, Dmitry Ioffe, Yuval Peres and Senya Shlosman for helpful discussions.
\end{subsection}
\end{section}
\begin{section}{Notation and tools}\label{secnott}
\begin{subsection}{Notation}\label{secnot}
\begin{definition}\label{defpathedge}
Elements of $\R^2$ will be denoted by boldface symbols. 
By a discrete path, we mean a list of elements of $\Z^2$, each being a nearest-neighbour of the preceding one, and without repetitions. 
In referring to a path, we mean a subset of $\R^2$ given by the union of the topologically closed edges formed from the set of consecutive pairs of vertices of some discrete path. 
(As such, a path is defined to be self-avoiding, including at its vertices.)
In a similar vein, any subset of $\R^2$ that is introduced as a connected set is understood to be a union of closed intervals $\big[ \bo{u},\bo{v} \big]$ corresponding to nearest-neighbour edges $(\bo{u},\bo{v})$.
For such a set $A$, we write $V(A) = A \cap \Z^2$ and $E(A)$ for the set of edges of which $A$
is comprised. 

For a general subset $A \subseteq \R^2$, we write $E(A)$ for the set of nearest-neighbour edges 
$(\bo{u},\bo{v}) \in E(\Z^2)$ such that $\big[ \bo{u},\bo{v} \big] \subseteq A$. (This is of course consistent with the preceding definition.) We write 
$E^*(A)$ for the set of nearest-neighbour edges $(\bo{u},\bo{v}) \in E(\Z^2)$ such that $\big[ \bo{u},\bo{v} \big] \cap A \not= \emptyset$. 
\end{definition} 
\begin{definition}\label{defclosedtri}
For  $\bo{x},\bo{y} \in \Z^2$, $\bo{y} \not= \bo{x}$,
we write $\ell_{\bo{x},\bo{y}}$ for the planar line containing $\bo{x}$ and $\bo{y}$, and $\ell^+_{\bo{x},\bo{y}}$ for the semi-infinite line segment that contains $\bo{y}$ and has endpoint $\bo{x}$. We write $\big[ \bo{x},\bo{y} \big]$ for the line segment whose endpoints are $\bo{x}$ and $\bo{y}$. 
We write $T_{\bo{0},\bo{x},\bo{y}}$ for the
closed triangle with vertices $\bo{0}$, $\bo{x}$ and $\bo{y}$.
For $\bo{x},\bo{y} \in \R^2$, 
we write $\ang\big(\bo{x},\bo{y} \big) \in [0,\pi]$ for the angle between these two vectors.
Borrowing complex notation, we set $\argu\big(\bo{x}\big)$ for the argument of $\bo{x}$.
In many derivations, the cones, line segments and points in question all lie in a cone, rooted at the origin, whose aperture has angle strictly less than $2\pi$. As such, it is understood that $\argu$
denotes a continuous branch of the argument that is defined throughout the region under consideration.
\end{definition}
Sometimes we wish to specify a cone by a pair of boundary points, and sometimes by the argument-values of its boundary lines:
\begin{definition}
For $\bo{x},\bo{y} \in \Z^2$, $\argu(\bo{x}) < \argu(\bo{y})$, write
$$
 A_{\bo{x},\bo{y}} = \Big\{ \bo{z} \in \R^2: \argu\big( \bo{x} \big) \leq  \argu\big( \bo{z} \big) \leq  \argu\big( \bo{y} \big)  \Big\} \cup \big\{ \bo{0} \big\}.
$$ 
Recall that, in (\ref{wnot}), we  specified a cone by the argument-values of its boundary lines.
Extending this notation,
for any $\bo{x} \in \Z^2$ and $c \in [0,\pi)$, we write 
$W_{\bo{v},c}\big( \bo{x} \big) = \bo{x} + W_{\bo{v},c}$. 
We also write, for $\bo{x} \in \R^2$ and $c \in (0,2\pi)$,
$$
W_{\bo{v},c}^+  = \Big\{ \bo{z} \in \R^2: \argu(\bo{v})  \leq  \argu(\bo{z}) \leq  \argu(\bo{v}) + c \Big\}
\cup \big\{ \bo{0} \big\}
$$
and 
$$
W_{\bo{v},c}^-  = \Big\{ \bo{z} \in \R^2: \argu(\bo{v}) - c \leq  \argu(\bo{z}) \leq  \argu(\bo{v})  \Big\}
\cup \big\{ \bo{0} \big\}.
$$
\end{definition}
\begin{definition}
For $\bo{v} \in \R^2$, let $\bo{v}^{\perp} \in S^1 $ denote the direction vector 
obtained by a counterclockwise turn of $\pi/2$ from the direction of $\bo{v}$. 
\end{definition}
\begin{definition}\label{defmarg}
For $P$ a probability measure on $\zoz$ and for $\omega' \in \{ 0,1 \}^A$ for some $A \subseteq E(\Z^2)$, we write $P_{\omega'}$ for the conditional law of $P$ given $\omega\big\vert_A = \omega'$. We will also write $P \big( \cdot \big\vert \omega' \big)$ for $P_{\omega'}$.
\end{definition}
\begin{definition}
Given a subset $A \subseteq \R^2$, two elements $\bo{x},\bo{y} \in \Z^2 \cap A$, and a 
configuration $\omega \in \zoz$, 
we write $\bo{x} \build\leftrightarrow_{}^A \bo{y}$
for the event that there exists an $\omega$-open path from $\bo{x}$ to $\bo{y}$ all of whose edges lie in $E(A)$.  By the ($\omega$-)open component of $\bo{x}$ in $A$, we mean the connected subset of $A$ whose members lie in an edge belonging to an ($\omega$-)open path in $E(A)$ that begins at $\bo{x}$.  
\end{definition}
Throughout, the notation $\dist \cdot \dist$ and $d \big( \cdot , \cdot \big)$ refers to the Euclidean metric on $\R^2$.  For $\gamma \subseteq \R^2$, we set ${\rm diam}(\gamma) = \sup \big\{ d\big(\bo{x},\bo{y}\big): \bo{x},\bo{y} \in \gamma \big\}$. For $K > 0$, we set $B_K = \big\{ \bo{x} \in \R^2: \dist \bo{x} \dist \leq K \big\}$.
\end{subsection}
\begin{subsection}{Decorrelation of well-separated regions: ratio-weak-mixing}\label{secrwm}
The following spatial decorrelation property is well-suited to analysing the conditioned circuit.
\begin{definition}\label{defrwm}
A probability measure $P$ on $\{ 0,1  \}^{E(\Z^2)}$ is said to satisfy the
ratio-weak-mixing property if, for some $\crwm,\lambda > 0$, and for all sets 
$\mc{D}, \mc{F} \subseteq E \big( \Z^2 \big)$,
$$
\sup \Big\{ \Big\vert \frac{P \big( D \cap F \big)}{P\big( D \big) P \big(
  F \big)} - 1 \Big\vert: D \in \sigma_{\mc{D}}, F \in \sigma_{\mc{F}}, 
  P \big( D \big) P \big( F \big) > 0   \Big\} 
\leq \crwm \sum_{x \in V(\mc{D}), y \in V(\mc{F})} e^{- \lambda \vert x - y \vert},
$$
whenever the right-hand-side of this expression is less than one. Here, for $A \subseteq E(\Z^2)$, $\sigma_A$ denotes the set of configuration events measurable with respect to the variables $\big\{ \omega(e): e \in A \big\}$.
\end{definition}
The ratio-weak-mixing property is satisfied by any  $P = \P_{\beta,q}$, with $\beta < \hat{\beta}_c$, that is, by any random cluster model with exponential decay of connectivity (Theorem 3.4, \cite{alexon}).
Any such measure trivially satisfies the following condition:
\begin{definition}\label{defbden}
A probability measure $P$ on $\{ 0,1  \}^{E(\Z^2)}$ satisfies
the bounded energy property if there exists a constant $\cposen > 0$ such that,  for any $\omega' \in \zoz$ and an edge $e \in E(\Z^2)$, the conditional probability that $\omega(e) = 1$ given the marginal $\omega' \big\vert_{E(\Z^2) \setminus \{ e \}}$ is bounded between $\cposen$ and $1 - \cposen$.
\end{definition}
The following tool (Lemma 2.1 of \cite{hammondone})
will establish near-independence of separated regions:
\begin{lemma}\label{lemkapab}
Given $A,B \subseteq E(\Z^2)$ and $m > 0$, let 
$$
\kappa_m \big( A,B \big) = \sum_{\bo{x} \in V(A),\bo{y} \in V(B), \vert\vert
  \bo{x} - \bo{y} \vert\vert \geq m} \exp \big\{ - \lambda \vert\vert
\bo{x} - \bo{y} \vert\vert \big\},
$$
where 
we write $V(A)$ for the set of vertices incident to one of the edges comprising $A$ (and similarly, of course, for $V(B)$).
Set $\phi_m \big( A,B \big) = \big\vert \big\{ \bo{a} \in A, \bo{b} \in B:
\vert\vert \bo{a} - \bo{b}
\vert\vert \leq m \big\} \big\vert$.
We say that $A$ and $B$ are $(m,C_0)$-well separated if $A \cap B = \emptyset$,
$\kappa_m \big( A,B \big) \leq 1/(2\crwm)$ and $\phi_m \big( A,B \big) \leq
C_0$. Here, $\crwm$ denotes the constant appearing in the definition (\ref{defrwm}) of ratio weak mixing. 
Let $P$ be a probability measure on $\{ 0,1  \}^{E(\Z^2)}$  satisfying the ratio-weak-mixing and bounded energy properties.
For $\clemkam \in \N$, $\clemkac \in \N$, there exists $\conka =
C\big(\clemkam,\clemkac\big)$ such that, if $A,B \subseteq E(\Z^2)$ are $(\clemkam,\clemkac)$-well
separated, then, for any $G \in \sigma_B$,
\begin{equation}\label{omubd}
\sup_{\omega \in \{ 0,1 \}^A} P_\omega \big( G \big) \leq \conka P(G) \qquad \quad \empty
\end{equation}
and
\begin{equation}\label{omlbd}
\inf_{\omega \in \{ 0,1 \}^A} P_\omega \big( G \big) \geq \conka^{-1} P(G). 
\end{equation}
\end{lemma}
\begin{subsubsection}{Large deviations of global distortion}
Recall the global distortion quantity from Definition \ref{defglobdis}.
A large deviations' estimate (Proposition 1 of \cite{hammondone})
on the macroscopic profile of the conditioned circuit will be valuable.
\begin{prop}\label{propglobdis}
Let $P = \P_{\beta,q}$ with $\beta < \hat{\beta}_c$  and $q \geq 1$. 
There exists $c > 0$
such that, for any  $\epsilon \in \big( 0, c \big)$, and for all $n \in \N$,
\begin{equation}\label{eqgd}
P \Big( \globdis > \epsilon n \Big\vert \big\vert {\rm INT} \big( \cir \big) \big\vert \geq n^2  \Big)  \leq \exp \big\{ - c \epsilon n \big\}.
\end{equation}
Under this measure, $\bo{0} \in {\rm INT} \big( \cir \big)$ except with exponentially decaying probability in $n$. Moreover, (\ref{eqgd}) holds under the conditional measure
$P \big( \cdot \big\vert \acon \big)$.
\end{prop}
The following (Lemma 2.3 of \cite{hammondone}) is an immediate consequence. 
\begin{lemma}\label{lemmac}
There exists $\epsilon > 0$, $\ccone > 0$ and $\cctwo > 0$ such that
$$
P \Big(   \cir \subseteq B_{\cctwo n} \setminus B_{\ccone n}    \Big\vert \acon \Big)
   \geq 1 - \exp \big\{ - \epsilon n \big\}. 
$$
\end{lemma}
\end{subsubsection}
\end{subsection}
\begin{subsection}{Ornstein-Zernike results for point-to-point connections}\label{secoz}
We recount the statements that we require from the theory \cite{civ} of point-to-point conditioned connections in a subcritical random cluster model. 

We record Theorem $A$ of \cite{civ} in the two-dimensional case:
 \begin{lemma}\label{lemozciv}
Let $P = \P_{\beta,q}$ with $\beta < \hat{\beta}_c$  and $q \geq 1$. Then
$$
P \Big(  \bo{0} \leftrightarrow \bo{x}  \Big) = 
\big\vert\big\vert  \bo{x}  \big\vert\big\vert^{- \frac{1}{2}}
 \Psi \big( \bo{n_x} \big) 
 \exp \Big\{   - \xi \big( \bo{n_x} \big)   \vert\vert \bo{n_x} \vert\vert \Big\}
   \Big( 1 + o(1) \Big),
$$
uniformly as $\bo{x} \to \infty$. The functions $\Psi$ and $\xi$ are positive, locally analytic functions on $S^1$, and $\bo{n_x} = \frac{\bo{x}}{\vert\vert \bo{x} \vert\vert}$. 
\end{lemma}
The following appears in Theorem B of \cite{civ}:
\begin{lemma}\label{lemozstr}
Let $P = \P_{\beta,q}$ with $\beta < \hat{\beta}_c$  and $q \geq 1$. 
Then $\wulff$ has a locally analytic, strictly convex boundary.
\end{lemma}
In \cite{civ}, by a refinement of the techniques of \cite{civearly},
local limit results 
such as Lemma \ref{lemozciv} are proved by using the following tool. 
The common open component $\overline\gamma_{\bo{x},\bo{y}}$ of $\bo{x}$ and $\bo{y}$ under $P \big( \cdot \big\vert  \xo \leftrightarrow \yo  \big)$ is split into pieces as follows. For $\delta > 0$, a point $\bo{v} \in \gamma_{\bo{x},\bo{y}}$ is called a $\delta$-cone point of $\overline\gamma_{\bo{x},\bo{y}}$ if the whole of $\overline\gamma_{\bo{x},\bo{y}}$ is contained in the union  $W_{-(\bo{y}-\bo{x}),\delta} \big( \bo{v} \big) \cup 
 W_{\bo{y}-\bo{x},\delta} \big( \bo{v} \big)$ of two aperture-$\delta$ cones emanating from $\bo{v}$ whose axes are aligned with $\bo{y} - \bo{x}$. It is argued that, typically, the set of such regeneration sites populates the cluster $\overline\gamma_{\bo{x},\bo{y}}$  under the conditioned measure in a very dense way. As such, the regeneration sites split the cluster into a large number of intervening ``sausages'', which, due to the exponential spatial decorrelation of subcritical random cluster measures, are in essence independent except at very close range. (Gaussian fluctuations then emerge due to the renewal structure apparent in this concatenation of small sausages.) 
One minor inconvenience with this definition of regeneration site is that, if $\bo{y} - \bo{x}$ is non-axial and $\delta > 0$ is small, then the above union of cones will not contain any of the edges of $\Z^2$ that neighbour $\bo{v}$, making the definition useless.  In \cite{civ}, $\delta$ is chosen large enough to avoid this problem. (See Section 2.9 of \cite{civ}.) It is crucial for our purpose that $\delta > 0$ may be chosen to be fixed but arbitrarily small. Hence, we must slightly alter the definition of cone site, adding to the union of cones a bounded ball about $\bo{v}$ so that the local lattice effect is handled. We now give the formal definition of this notion of regeneration site. Note that of course we have already defined a notion of regeneration site for the conditioned circuit $\cir$ that is our object of study in this article. As such, we use the term connection regeneration site to distinguish the two notions. 
%\begin{definition}
%Let $\bo{x},\bo{y} \in \Z^2$, and let $C$ be a connected set containing 
%$\bo{x}$ and $\bo{y}$. For $\delta > 0$, $\bo{v} \in V(C)$ will be called a
%$(\delta,K)$-cone point if
%$$
%C \setminus B_K(\bo{v}) \subseteq C_{-(\bo{y} - \bo{x}),\delta} \big(
%\bo{v} \big)
% \cup C_{\bo{y} - \bo{x},\delta} \big( \bo{v} \big).
%$$
%We write $\cdk$ for the set of $(\delta,K)$-cone points of $C$.
%The connected components of $C \setminus \cdk$ will be called
%$(\delta,K)$-regeneration clusters of $C$. We write $\maxozrg_C$ for the size
%of the largest such component.
%\end{definition}
%\begin{lemma}\label{lemozu}
%For each $\delta > 0$, there exists $K \in \N$ and $c > 0$ such that
%$$
%P \Big(  \cdk \big( C_{\bo{x},\bo{y}} \big) = \emptyset  \Big\vert \bo{x}
%\leftrightarrow \bo{y} \Big) \leq \exp \Big\{ - c \vert\vert \bo{x} - \bo{y} \vert\vert \Big\}.
%$$
%\end{lemma}
%\begin{lemma}
%For each $\delta > 0$, there exists $K \in \N$ and $c > 0$ such that
%$$
%P \Big( \maxozrg_{C_{\bo{x},\bo{y}}} \geq t \log \vert\vert \bo{x} - \bo{y}
%\vert\vert \Big\vert \bo{x}
%\leftrightarrow \bo{y} \Big) \leq \exp \big\{ - c t \big\}.
%$$
%\end{lemma}
\begin{definition}\label{defpathreg}
Let $\bo{x},\bo{y} \in \Z^2$.
For $\delta > 0$, fix $K \in \N$ large enough that there exists a path with edge-set contained in $E^*(B_K)$  that intersects both  
$\partial B_K \cap 
W_{-(\bo{y} - \bo{x}), \delta}\big( \bo{0} \big)$ and 
$\partial B_K \cap W_{\bo{y} - \bo{x}, \delta}\big( \bo{0} \big)$, and contains $\bo{0}$.
Let $\phi$ denote such a path.

Let $\gamma$ denote a connected set such that $\bo{x}, \bo{y} \in \gamma$.
A vertex $\bo{v} \in V(\gamma)$ will be called a $(\delta,K,\phi)$-connection regeneration site of $\gamma$ (or, in full, a 
$(\bo{x},\bo{y},\delta,K,\phi)$-connection regeneration site of $\gamma$)
if
$$
\gamma \setminus B_K (\bo{v}) \subseteq W_{-(\bo{y}-\bo{x}),\delta} \big( \bo{v} \big) \cup 
 W_{\bo{y}-\bo{x},\delta} \big( \bo{v} \big), \quad \textrm{and} \, \, \, E(\gamma) \cap E^* \big( B_K(\bo{v}) \big) = \bo{v} + \phi.
$$
We write ${\rm CRG}_{\delta,K,\phi}^{\bo{x},\bo{y}}\big(\gamma\big)$ 
for the set of such sites. We also write ${\rm CRG}_{\delta,K}^{\bo{x},\bo{y}}$ 
for the union of ${\rm CRG}_{\delta,K,\phi}^{\bo{x},\bo{y}}\big(\gamma\big)$  over all paths
$\phi \subseteq E^*(B_K)$ as above.

We call the connected components of $\gamma \setminus \cdkp\big(\gamma\big)$ the connection regeneration clusters of $\gamma$. We call such a cluster internal unless it contains either $\bo{x}$ or $\bo{y}$. An internal connection regeneration cluster $C$ of $\gamma$ has precisely two elements of $\cdkp\big(\gamma\big)$ in its boundary. These may be labelled $\bo{f}(C)$ and $\bo{b}(C)$ in such a way that $C \setminus B_K  \big(  \bo{b}(C) \cup \bo{f}(C) \big)
 \subseteq W_{-(\bo{y}-\bo{x}),\delta} \big( \bo{b}(C) \big) \cap  
 W_{\bo{y}-\bo{x},\delta} \big(  \bo{f}(C)  \big)$.  We define the displacement $\bo{d}\big( C \big)$ of an internal cluster $C$ by  $\bo{d}\big( C \big) = \bo{b} - \bo{f}$. In the case of the cluster $C$ containing $\bo{x}$ (or $\bo{y}$), we set $\bo{d}(C)$ equal to $\bo{r} - \bo{x}$ (or $\bo{y} - \bo{r}$), where 
$\bo{r}$ is the unique element of $\cdkp\big(\gamma\big)$ in the boundary of $C$. 

We write $\maxozrg_{\delta,K,\phi}^{\bo{x},\bo{y}}(\gamma)$ to be the maximum of $\vert\vert \bo{d}\big( C \big) \vert\vert$
over all connection regeneration clusters $C$ of $\gamma$.
% for the minimal value of $k \in \N$ such that there exists a sequence
%$\bo{x} = \bo{p_0},\bo{p_1},\ldots,\bo{p_r} = \bo{y}$ satisfying
%$\bo{p_i} \in \cdkp(C)$ for $1 \leq i \leq r-1$ and
%$\vert\vert \bo{p_i} - \bo{p_{i-1}} \vert\vert \leq k$ for $1 \leq i \leq r$.
\end{definition}
\noindent{\bf Remark.}
The use of balls $B_K$ to permit small-aperture cones to be used in the definition of connection regeneration sites necessitates some minor changes in the formulation of the proofs of \cite{civ}. In defining the percolation event $\Lambda = \mathcal{W} \mathcal{N}$ in Section 3.1 of \cite{civ}, the precise definition of the event $\mathcal{E}_N(\gamma)$ in (3.1) must be modified to specify a partition of closed edges across the $B_K$ balls about $\bo{f}$ and $\bo{b}$. To record the details of these changes would be unenlightening, and we leave them to a patient reader. \\
\begin{figure}\label{figpathreg}
\begin{center}
\includegraphics[width=0.8\textwidth]{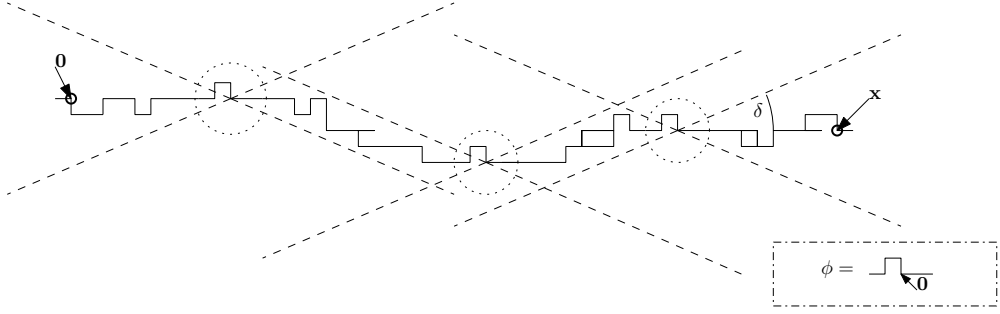} \\
\end{center}
\caption{Illustrating Definition \ref{defpathreg}. The cluster $\gamma_{\bo{0},\bo{x}}$ is depicted under a sample of $P \big( \cdot \big\vert  \bo{0} \leftrightarrow \bo{x}  \big)$. The constant $K$ equals $2$, and $\phi$ equals the six-edged connected set indicated in the separate sketch in the box. The three elements of ${\rm CRG}_{\delta,K,\phi}^{\bo{0},\bo{x}}\big( \gamma_{\bo{0},\bo{x}} \big)$ lie at the centres of the radius-$K$ circles.}
\end{figure}
In the next three lemmas, $\overline\gamma_{\bo{x},\bo{y}}$ denotes the common open component of $\bo{x}$ and $\bo{y}$ (which exists in the contexts in question).
The first two lemmas follow directly from (1.6) and (1.8) of \cite{civ}.
\begin{lemma}\label{lemmaxreg}
Let  $\phi \subseteq E^*(B_K)$ for which $\bo{0} \in \phi$
contain a path from 
$\partial B_K \cap W_{-(\bo{y}-\bo{x}),\delta}$ 
to  
$\partial B_K \cap W_{\bo{y}-\bo{x},\delta}$.
Then there exist $C > c > 0$ such that, for $t > C \log \vert\vert \bo{x} - \bo{y} \vert\vert$,
$$
P \Big( \maxozrg_{\delta,K,\phi}^{\bo{x},\bo{y}} \big( \overline\gamma_{\bo{x},\bo{y}} \big) > t \Big\vert \bo{x} \leftrightarrow \bo{y} \Big)
\leq \exp \big\{ - c t \big\}. 
$$
\end{lemma}
\begin{lemma}\label{lemoznor}
Let $P = \P_{\beta,q}$ with $\beta < \hat{\beta}_c$. 
For all $\delta > 0$, there exists $K = K(\delta) \in \N$ and $c = c(\delta) > 0$ such that,
for all $\bo{x},\bo{y} \in \Z^2$,
$$
P \bigg(  
   \overline\gamma_{\bo{x},\bo{y}} \subseteq
   \Big( W_{\bo{y}-\bo{x},\delta}\big( \bo{x} \big) \cap 
  W_{-(\bo{y}-\bo{x}),\delta}\big( \bo{x} \big)  \Big) \cup B_K\big(\bo{x}\big)
 \cup B_K\big(\bo{y}\big)
\Big\vert
 \bo{x} \leftrightarrow \bo{y}  \bigg) \geq c. 
$$
\end{lemma}
\begin{definition}\label{deffluc}
Let $\gamma$ denote a connected set containing $\bo{x},\bo{y} \in \Z^2$.
We write 
$$
 {\rm fluc}_{\bo{x},\bo{y}}\big( \gamma \big) = \sup \Big\{ d\big( \bo{z} , \big[ \bo{x},\bo{y} \big] \big): \bo{z} \in \gamma \Big\}.
$$
 \end{definition}
We require Lemma 2.6 of \cite{hammondone}, a bound on moderate fluctuations of conditioned connections (and a consequence of \cite{civ}):
\begin{lemma}\label{lemmdf}
Let $P = \P_{\beta,q}$ with $\beta < \hat{\beta}_c$. 
There exists a constant $c > 0$ such that, for all $\bo{x},\bo{x} \in \Z^2$ and $0 < t < c \vert\vert \bo{x} - \bo{y}
\vert\vert^{1/2}$,
$$
P \Big(   {\rm fluc}_{\bo{x},\bo{y}}\big( \overline\gamma_{\bo{x},\bo{y}} \big) 
 \geq  \vert\vert \bo{x} - \bo{y}
\vert\vert^{1/2} t   \Big\vert \bo{x} \leftrightarrow \bo{y} \Big)
 \leq \exp \big\{ - c t^2 \big\}. 
$$
\end{lemma}
\end{subsection}
\begin{subsection}{Some comments on the required hypotheses}
Most of the arguments in this paper and its companions \cite{hammondone} and \cite{hammondtwo} use hypotheses that are a little weaker than insisting that $P = \P_{\beta,q}$, (with $\beta < \hat{\beta}_c$ and $q \geq 1$), be a subcritical random cluster measure. The basic hypotheses to which we will allude in the proofs are:
\begin{itemize}\label{assump}
\item $P$ satisfies the ratio weak mixing property (\ref{defrwm}),
\item $P$ has exponential decay of connectivity; that is, 
namely, that there exists $c > 0$ such that $P_{\omega} \big( \bo{0} \to \partial B_n \big) \leq \exp \big\{ -cn \big\}$ for all $n \in \N$ and $\omega \in \{0,1\}^{E(Z^2) \setminus E(B_n)}$,
\item $P$ satisfies the bounded energy property (\ref{defbden}),
\item $P$ is translation-invariant,
\item $P$ is invariant under axial symmetry.
\end{itemize}
As we have noted, the ratio weak mixing and bounded energy properties are satisfied by any $\P_{\beta,q}$, $\beta < \hat{\beta}_c$. Exponential decay of connectivity follows from $\beta < \hat{\beta}_c$ and the ratio weak mixing property, and the other listed hypotheses are trivially satisfied. See \cite{hammondone} for a discussion of the limited use of additional hypotheses in the sequence \cite{hammondone}, \cite{hammondtwo} and the present article.
\begin{subsubsection}{Final preliminaries}\label{secfp}
The definition of  $\cir$-regeneration site was specified in terms of two constants $\qzero > 0$ and
$c_0 \in \big( 0 , \qzero/2 \big)$. We now record the precise conditions that these two constants are required to satisfy.
\begin{definition}\label{defrg}
For $\bo{u} \in S^1$,
let $w_{\bo{u}}$ denote the counterclockwise-oriented unit tangent vector
to $\partial \wulff$ at $\partial \wulff \cap \big\{ t\bo{u}: t \geq 0  \big\}$.
Let $\qzero > 0$ satisfy
\begin{equation}\label{supang}
\sup_{\bo{z} \in S^1} {\rm ang} \big( w_{\bo{z}} , \bo{z}^{\perp} \big) \leq \pi/2 - 4\qzero,
\end{equation}
a choice made possible by the compactness and convexity of $\wulff$.

Let $c_0 > 0$ is chosen so that, whenever $\bo{x},\bo{y} \in \partial
\wulff$,
$\argu(\bo{x}) < \argu(\bo{y})$ and 
$\ang \big( \bo{x} , \bo{y} \big) \leq 2 c_0$, then
\begin{equation}\label{czercond}
{\rm ang} \big( \bo{x} - \bo{y}, - \bo{y}^{\perp}  \big) \leq
\pi/2 - 3\qzero.  
\end{equation}
This is possible by (\ref{supang}) and Lemma \ref{lemozstr}.
\end{definition}
We also record Lemma 2.5 of \cite{hammondone} regarding control of $\cir$ near regeneration sites:
\begin{lemma}\label{lemdistang}
If $\bo{x},\bo{y} \in \R^2$ satisfy $\ang \big( \bo{x}, \bo{y} \big) \leq c_0$   and $\bo{y} \in \clu{\bo{x}} \cup \clum{\bo{x}}$, then $\vert\vert \bo{y} - \bo{x}  \vert\vert \leq \csc \big( \qzero/2 \big) \vert\vert \bo{x} \vert\vert \ang \big( \bo{x} , \bo{y} \big)$. 
\end{lemma}
\end{subsubsection}
\begin{subsubsection}{Convention regarding constants}
Some constants are fixed in all arguments: notably, $c_0$ and $\qzero$ have just been specified.  A few constants have been or will be fixed in lemmas, these constants carrying letter subscripts that are acronyms evident from the defining context. 
Throughout, constants with a capital $C$ indicate large constants, and, with a lower case, small constants.
 The constants $C$ and $c$ may change from line to line and are used to absorb and simplify expressions involving other constants. 
\end{subsubsection}
\end{subsection}
\end{section}
\begin{section}{Any positive-angle cone contains regeneration sites}\label{secmac}
Our first step on the route to proving Theorem \ref{thmmaxrg} is to show that set of regeneration sites $\reg$ of the conditioned circuit $\cir$ intersect any positive-angle sector. 
\begin{prop}\label{proprgmac}
Let $P = \P_{\beta,q}$ with $\beta < \hat{\beta}_c$ and $q \geq 1$. There exists $c > 0$  
such that, for any $\epsilon \in \big( 0,
c \big)$,
$$
P \Big( 
  \mar > \epsilon  \Big\vert \acon \Big) \leq \exp \big\{ - c \epsilon n \big\}.
$$ 
\end{prop}
We explain the ideas of the proof before beginning formally.
By Proposition \ref{propglobdis}, 
and a union bound, it suffices to show that there exists $c > 0$ such that, for $C_0 > 2c_0$ and for each $\bo{u} \in S^1$,
writing 
$$
a_{\bo{u}} = 
P \Big( \csec \cap \reg = \emptyset \Big\vert \acon, \globdis \leq
  \frac{\epsilon n}{C_0}   \Big),
$$
we have that
\begin{equation}\label{aubd}
a_{\bo{u}} \leq \exp \big\{ - c \epsilon n \big\}.
\end{equation}
We now summarize the approach to proving this assertion.
Suppose given a configuration  $\omega \in \zoz$ realizing 
$\big\{ \acon, \csec \cap \reg = \emptyset, \globdis \leq
  \frac{\epsilon}{C_0} n \big\}$. Let $\bo{x}$ and $\bo{y}$ denote the points of contact of $\cir$ under $\omega$ with the opposite sides of $\csec$. These points may of course not be unique, but we prefer to neglect this problem for now and return to it at the end of the sketch. 
 Set $\tilde\omega = \omega \big\vert_{\csece^c}$.  
It suffices for (\ref{aubd}) that
\begin{equation}\label{sketone}
 P_{\tilde\omega} \Big( \acon , \globdis \leq
  \frac{\epsilon}{C_0} n ,  \csec \cap \reg = \emptyset \Big) \leq \exp \big\{ - c \epsilon n \big\} P \big( \bo{x} \to \bo{y} \big)
\end{equation}
and
\begin{equation}\label{skettwo}
 P_{\tilde\omega} \Big(  \acon,  \globdis \leq
  \frac{\epsilon n}{C_0}    \Big) \geq \exp \big\{ - c' \epsilon n \big\} P \big( \bo{x} \to \bo{y} \big),
\end{equation}
where $c' < c$. 
\begin{figure}\label{figmrsumbox}
\begin{center}
\includegraphics[width=0.7\textwidth]{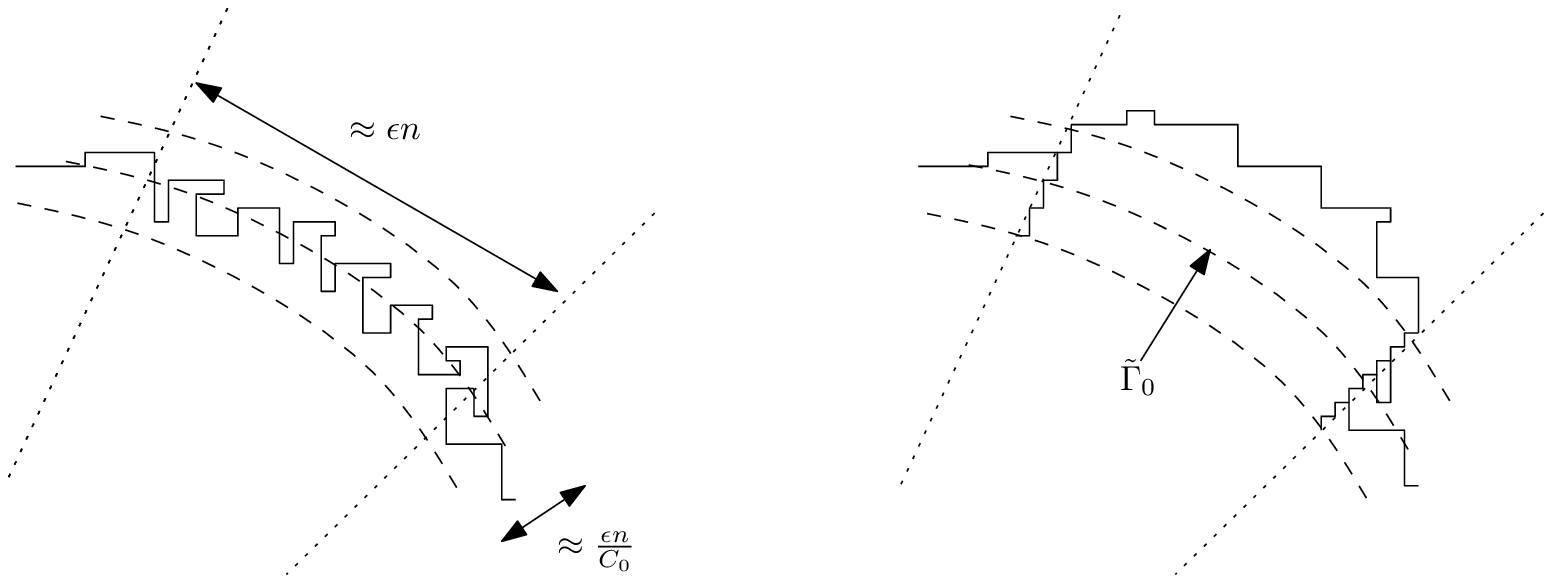} \\
\end{center}
\caption{A pictorial summary of the proof of Proposition \ref{proprgmac}. Given $\acon$ and $\globdis \leq \epsilon n/C_0$, the circuit $\cir$ may realize $\reg \cap \csec = \emptyset$ only by continual local backtracking as it passes through $\csec$, as depicted in the first figure. This peculiarity carries a probabilistic cost of $\exp \big\{ -c \epsilon n\}$ beyond the mere existence of a connecting path between the opposite sides of $\csec$. 
Given that such a circuit segment  $\cir \cap \csec$ is necessarily trapped in a corridor of width $\epsilon n/C_0$ centred on a dilate $\tcir :=  n \partial \wulff$ of the Wulff curve $\partial \wulff$, it is, if $C_0$ is fixed but high, probabilistically cheaper to forge a new circuit, trapping as much area as the first one, by sealing open paths along the two short sides of the corridor and then connecting these seals on the outside. The second figure depicts this alternative. That it is indeed cheaper depends on the differentiability of $\partial \wulff$, which ensures that the corridor is not kinking.}
\end{figure}
Roughly put, (\ref{sketone}) holds because,  given $\acon$ and $\globdis \leq
  \frac{\epsilon n}{C_0}$, the occurrence of 
 $\csec \cap \reg = \emptyset$ entails that the circuit $\reg$ crosses the sector $\csec$ in an inefficient way, with continual local backtracking required so that $\reg \cap \csec = \emptyset$.
 The reason that such backtracking is required is that, if the condition $\bo{v} \in \reg$ is to be violated for all $\bo{v} \in \cir \cap \csec$, then the curve segment $\cir \cap W_{\bo{v} - c_0,\bo{v} + c_0}$ must leave $\clu{\bo{v}} \cup \clum{\bo{v}}$ for all such $\bo{v}$. Our assumption that $\globdis$ is a small multiple of $\epsilon n$ means that $\cir$ is close to the dilate $\tcir = n \partial \wulff$ of the smooth Wulff curve, and this forces $\bo{u} \in \cir$ inside   $\clu{\bo{v}} \cup \clum{\bo{v}}$
provided that $\ang\big( \bo{u},\bo{v} \big)$ is bigger than a small constant multiple of $\epsilon$ (so that the angular change beats the error in the approximation of $\cir$ by $\tcir$) but smaller than $c_0$ (after which, the Wulff curve itself may turn too much).
This means that, if $\bo{v} \not\in \reg$ is to hold, the curve $\cir$ must jump out of  $\clu{\bo{v}} \cup \clum{\bo{v}}$ close to $\bo{v}$, specifically, inside $\csec$, provided that $\bo{v}$ is not close to the boundary of $\csec$. 
Recall Definition \ref{defpathreg}. 
We will show that, for small $\delta > 0$, if $\bo{v}$ is a $\delta$-connection regeneration site for the connection $\bo{x} \to \bo{y}$ given by the path $\cir \cap \csec$, then in fact  $\cir \cap \csec \subseteq \clu{\bo{v}} \cup \clum{\bo{v}}$. So the only way  $\csec \cap \reg = \emptyset$  may occur is if this connection has no connection regeneration sites. As such, Lemma \ref{defpathreg} yields (\ref{sketone}).

The second bound (\ref{skettwo}) makes a crucial use of the differentiability of the Wulff curve. The circuit $\cir$ under the original configuration $\omega$ is approximated with a $O\big( \epsilon n/C_0 \big)$ error by some dilate $\tcir$ of the Wulff shape. Suppose that $\tcir$ cuts the two sides of $\csec$ at $\bo{x'}$ and $\bo{y'}$. The dilate has diameter of order $n$, the sector $\csec$ has an angle of $\epsilon$, and the Wulff curve is second differentiable. These ingredients tell us that, inside $\csec$, $\tcir$ is never at distance more than $\Theta \big( \epsilon^2 n \big)$ from the line segment $\big[\bo{x'},\bo{y'}\big]$. This means that, if under the resampling $\omega'$ in the sector $\csec$, there is a path from $\bo{x}$ to $\bo{y}$ that begins (with a segment $Q_1$) by tracing its way from $\bo{x}$ along one boundary of $\csec$ away from the origin for a distance $2\epsilon n/C_0$, and ends (with a segment $Q_2$) by  tracing its way to $\bo{y}$ along the opposite boundary of $\csec$ towards the origin for the same distance, while crossing $\csec$ in between in a typical fashion, then the resulting full-plane configuration does satisfy $\acon$. This is because a path with such a start and finish has fluctuated outwards enough to trap at least as much area as its counterpart under the original $\omega$ did, since the outward displacement of the new path from $\big[\bo{x'},\bo{y'} \big]$ is greater than both the displacement between      
$\big[\bo{x'},\bo{y'} \big]$ and $\tcir \cap \csec$, and the error arising from the approximation of the curve $\cir$ under $\omega$ by the dilate $\tcir$. The probability of such a path is at least $\exp \big\{ - C \epsilon n/C_0\big\} P \big( \bo{x} \to \bo{y} \big)$, where $C$ is a universal constant, the argument of the exponential corresponding in essence (by means of the bounded energy property of $P$)  to the length of the beginning and end sections of the constructed path. We achieve the required $c' < c$, by fixing $C_0 > 0$ high enough.

Of course, in this discussion, we pretended that we knew that $\cir$ under $\omega$ meets each boundary line segment of $\csec$ just once. This may be untrue. However, 
by $\globdis \leq \epsilon n/C_0$, each of these intersections occurs in an interval of length at most a constant multiple of $\epsilon n/C_0$. The condition on $\omega'$ in the preceding paragraph should thus be altered so that these intervals are sealed by a path of open edges, so that the resulting configuration does not have gaps in the circuit around the boundary of $\csec$. The paths that must be opened are continuations of $Q_1$ and $Q_2$ of comparable length, so the preceding argument may be adapted to deal with this problem. 

In the proof itself, we will make use of the sector storage-replacement operation 
that was introduced in \cite{hammondone}, because it provides a convenient framework 
for extracting peculiarities from a conditioned circuit and renewing the effected section of the circuit. 
While logically self-contained, the formal argument is based on the template provided by the proof of Proposition 2 of \cite{hammondone}. 
There are also a few notational complications involved. As such, the reader is advised to read the overview of the operation provided in Section 1.1 of \cite{hammondone} and to briefly review the proof of Proposition 3 in that paper before proceeding.

To recall, then, the definition of the operation:
\begin{definition}\label{defsrro}
Let $\bo{x},\bo{y} \in \Z^2$, $\argu(\bo{x}) < \argu(\bo{y})$, be given.
Let $P$ be a given measure on configurations $\zoz$. 
The sector storage-replacement operation 
$\sigma_{\bo{x},\bo{y}}$ is a random map 
$$
\sigma_{\bo{x},\bo{y}}: \zoz \to \zoz \times \big\{ 0,1 \big\}^{\axye},
$$
whose law is specified in terms of $\bo{x},\bo{y}$ and $P$.
The output $\sigma_{\bo{x},\bo{y}}(\omega) = \big(\omega_1,\omega_2  \big)$
is given as follows. 

We set $\omega_2 = \omega \big\vert_{E(A_{\bo{x},\bo{y}})}$. We then define
$\omega' \in  \{0,1 \}^{E(A_{\bo{x},\bo{y}})}$ to be a random variable 
having the marginal in $E(A_{\bo{x},\bo{y}})$ of the law $P \big(\cdot \big\vert \omega\vert_{E(\Z^2) \setminus
    E(A_{\bo{x},\bo{y}})} \big)$ (specified in Definition \ref{defmarg}).
 We set 
\begin{eqnarray}
 \omega_1 & = & \omega' \qquad \qquad \textrm{on $E(A_{\bo{x},\bo{y}})$}
\nonumber \\
 & = &  \omega \qquad \qquad \textrm{on $E(\Z^2) \setminus E(A_{\bo{x},\bo{y}})$.}
\nonumber
\end{eqnarray}
That is, in acting on $\omega$, we begin by removing the contents of  $E\big(A_{\bo{x},\bo{y}}\big)$ and storing this information as
$\omega_2$.
We then resample these bonds subject to the untouched
information in the complement  $E(\Z^2) \setminus E(A_{\bo{x},\bo{y}})$. The new configuration,
that coincides with the original one in  $E(\Z^2) \setminus
E(A_{\bo{x},\bo{y}})$, is recorded as $\omega_1$.

We will call $\omega_1 \in \zoz$ the full-plane output, 
$\omega_2 \in \big\{ 0,1 \big\}^{\axye}$ the sector output, and
$\omega_1 \big\vert_{\axye}$ the updated configuration.
\end{definition}
In all its applications, the sector storage-replacement operation will act in the following way.
\begin{definition}\label{defregact}
Let $\bo{x},\bo{y} \in \Z^2$, $\argu(\bo{x}) < \argu(\bo{y})$, be given. 
The sector storage-replacement operation 
$\sigma_{\bo{x},\bo{y}}$ will be said to act regularly if 
\begin{itemize}
 \item the input configuration has the distribution $P$, and
 \item 
the randomness of the action is chosen such that, 
given the input $\omega \big\vert_{E(\Z^2) \setminus \axyoe}$, the updated configuration
$\omega_1 \big\vert_{\axyoe}$ is conditionally independent of the stored sector configuration $\omega_2 \big\vert_{\axye} = \omega \big\vert_{\axye}$.
\end{itemize}
\end{definition}
We require a minor change in definition from applications of the operation elsewhere (in \cite{hammondone} or later in the present paper): 
\begin{definition}
For $\bo{u} \in S^1$ and $\delta > 0$, let $\sigma_{\bo{u},\delta}$
denote the operation given by the definition of the sector
storage-replacement operation $\sigma_{\xo,\yo}$ in which the sector
$A_{\xo,\yo}$ is substituted by $W_{\bo{u},\delta}\big( \bo{0} \big)$. The definition of regular action of $\sigma_{\bo{u},\delta}$ is similarly defined.
\end{definition}
Before beginning the formal argument, we make some further definitions, and recall from \cite{hammondone} the notion of {\it good area capture} and the associated Lemma 3.1:
\begin{definition}\label{defgagab}
Let $\bo{x},\bo{y} \in \Z^2$ satisfy $\argu(\bo{x}) < \argu(\bo{y})$.
Let $\omega \in \zoz$ realize the events that $\bo{x} \build\leftrightarrow_{}^{A_{\bo{x},\bo{y}}} \bo{y}$ and that the common cluster of $\bo{x}$ and $\bo{y}$ in $A_{\bo{x},\bo{y}}$ is finite. We define the outermost open path from $\bo{x}$ to $\bo{y}$ in $A_{\bo{x},\bo{y}}$ to be  the open path $\sentier$ 
from $\bo{x}$ to $\bo{y}$ in $A_{\bo{x},\bo{y}}$ such that 
the bounded component of 
$\axy \setminus \sentier$ is maximal. We denote this path by %\hfff{gamoop} 
$\gamma_{\bo{x},\bo{y}} = \gamma_{\bo{x},\bo{y}}(\omega)$. (Note that there is only one bounded component, because a path is, by our definition, self-avoiding.)
We further write  %\hfff{ogamoop} 
$\overline{\gamma}_{\bo{x},\bo{y}} =  \overline{\gamma}_{\bo{x},\bo{y}}(\omega)$ for the common $\omega$-open cluster of $\bo{x}$ and $\bo{y}$ in $\axy$. 
\end{definition}
\noindent{\bf Remark.} Note that $\gamma_{\bo{x},\bo{y}}$ is almost surely well-defined under $P$ given
  $\bo{x} \build\leftrightarrow_{}^{A_{\bo{x},\bo{y}}} \bo{y}$, with $P = \P_{\beta,q}$, $\beta < \hat{\beta}_c$, since all open clusters are finite $P$-a.s. (as implied by the remark after Definition \ref{defcirsur}).
Note also that there is no difficulty in regard to the uniqueness of $\gamma_{\bo{x},\bo{y}}$, 
because we have decided to identify any path with the subset of $\R^2$ given by the union of the edges contained in the corresponding discrete path.
\begin{definition}\label{defgac}
Let $\bo{x},\bo{y} \in \Z^2$, and let $\gamma \subseteq \axy$ denote a path from $\bo{x}$
to $\bo{y}$. 
Let $I_{\bo{x},\bo{y}}\big( \gamma \big) \subseteq \R^2$ denote the bounded component of $\axy  \setminus \gamma$. The path $\gamma$ is said to have
$\epsilon$-{\rm good area capture} in $\axy$ if ${\rm diam} \big( \gamma \big) \leq \cgac \vert\vert
\bo{x} - \bo{y} \vert\vert$ and
\begin{equation}\label{eqgac}
 \Big\vert I_{\axy} \big( \gamma \big) \Big\vert
  \geq \big\vert  T_{\bo{0},\bo{x},\bo{y}} \big\vert + \epsilon
  \vert\vert \bo{x} - \bo{y} \vert\vert^{3/2} \big( \log \vert\vert \bo{x} -
  \bo{y} \vert\vert \big)^{1/2}.
\end{equation}
For $\epsilon > 0$, 
we write ${\rm GAC}\big(\bo{x},\bo{y}, \epsilon \big)$ for the subset of 
configurations $\omega \in \{ 0,1 \}^{\axye}$
such that
$\bo{x}  \build\leftrightarrow_{}^{\axy} \bo{y}$ under $\omega$,
 and for which 
$\gamma_{\bo{x},\bo{y}}$  has $\epsilon$-good area capture in $\axy$.

The path $\gamma$ is said to have good area capture if it has $1/10$-good area capture. We set 
${\rm GAC}\big(\bo{x},\bo{y} \big) = {\rm GAC}\big(\bo{x},\bo{y}, 1/10 \big)$.
\end{definition}
\begin{lemma}\label{lemgac}
Let $P = \P_{\beta,q}$, with $\beta < \hat{\beta}_c$  and $q \geq 1$.  For any sufficiently high constant $\clemgac > 0$ the following holds. Let $\bo{x},\bo{y} \in \Z^2$ 
satisfy $\arg(\bo{x}) < \arg(\bo{y})$, $\vert\vert \bo{x} \vert\vert, \vert\vert \bo{y} \vert\vert \leq \clemgac n$, $\vert\vert \bo{y} - \bo{x} \vert\vert \geq 4 \clemgac \log n$, 
$\bo{y} \in  C_{\pi/2 - \qzero}^F \big( \bo{x} \big)$ and 
$\bo{x} \in  C_{\pi/2 - \qzero}^B \big( \bo{y} \big)$ .
Let $\omega \in \{0,1\}^{E(\Z^2)
  \setminus \axye}$ be arbitrary. 
Then, for $\epsilon > 0$,
$$
P\Big(  {\rm GAC}\big( \bo{x},\bo{y},\epsilon \big) \Big\vert \omega\big\vert_{E(\Z^2) \setminus E(\axy)} \Big)
 \geq n^{-\clemgac \epsilon} P \Big(  \bo{x} \build\leftrightarrow_{}^{\axy} \bo{y} \Big).
$$
\qed
\end{lemma}
\noindent{\it Remark.} The precise definition of good area capture is being used purely for convenience, in the sense that it has been proved in \cite{hammondone} in this form. We have no need for the second term on the right-hand-side of (\ref{eqgac}) in the applications in the present paper. \\
\noindent{\bf Proof of Proposition \ref{proprgmac}.}
We will consider the regular action of the  sector storage-replacement operation
$\sigma_{(u,\epsilon/2)}$. \\
\noindent{\bf Definition of satisfactory input.}
The input $\omega \in \zoz$ will be defined to be satisfactory if it
realizes the event 
$$
 {\rm SAT}  : = 
 \Big\{  \acon, \globdis \leq
  \frac{\epsilon n}{C_0} , \cir \subseteq B_{\cctwo n}  \setminus  B_{\ccone n},
  \csec \cap \reg = \emptyset \Big\}, 
$$
where $C_0$ is a large constant, to be specified later. Throughout the proof, we write 
$\cir = \cir(\omega)$ as specified in Definition \ref{defncentr}. Since we will consider only 
$\omega \in {\rm SAT}$, $\cir$ is non-empty and $\omega$-open. \\ 
\noindent{\bf Definition of successful action.} We require:
\begin{definition}
Let $\ell$ denote a semi-infinite planar line segment that has $\bo{0}$
as an endpoint. We introduce two semi-infinite simple lattice paths
associated with $\ell$: the clockwise and counterclockwise boundary
segments of $\ell$, $\partial^- \ell$ and $\partial^+ \ell$. 

To define $\partial^- \ell$, suppose firstly that $\bo{u} \in \ell$
with $\arg(\bo{u}) \in [ 0,\pi/2)$.
The starting point of
$\partial^- \ell = \big\{ \bo{x_0},\bo{x_1},\ldots \big\}$ is taken to be
$\bo{x_0} = \bo{0}$.
The simple path $\big\{ \bo{x_0},\bo{x_1},\ldots \big\}$
will be constructed such that $\arg(\bo{x_i}) \leq \arg (\bo{u})$ for each $i \in \N$.
Given $\big\{ \bo{x_0},\ldots,\bo{x_{i-1}} \big\}$, we set
$\bo{x_i} = \bo{x_{i-1}} + (0,1)$ provided that  
  $\arg \big( \bo{x_{i-1}} + (0,1)  \big) \leq \arg (\bo{u})$.
Otherwise, we set $\bo{x_i} = \bo{x_{i-1}} + (1,0)$.

The path $\partial^+ \ell = \big\{ \bo{y_0}, \bo{y_1}, \ldots  \big\}$ is constructed
so that $\arg \big( \bo{y_i} \big) \geq \arg \big( \bo{u} \big)$,
with the choice 
$\bo{y_i} = \bo{y_{i-1}} + (0,1)$ being made if   
  $\arg \big( \bo{y_{i-1}} + (0,1) \big) \geq \arg (\bo{u})$, and, otherwise,
$\bo{y_i} = \bo{y_{i-1}} + (1,0)$.

A similar definition is used in the case of general $\arg \big( \bo{u}
\big)$.
To be explicit, let $I: \R^2 \to \R^2$
denote a rotation by a multiple of $\pi/2$ such that $I(\ell)$ lies in
the first quadrant. We set $\partial^+ \ell = I^{-1} \big( \partial^+
I(\ell) \big)$ and $\partial^- \ell = I^{-1} \big( \partial^- I(\ell) \big)$.

Finally, whenever $\bo{u},\bo{v} \in \Z^2$ satisfy
$\bo{u},\bo{v} \in \partial^+ \ell$ with $\bo{u}$ encountered before
$\bo{v}$ in the path $\partial^+ \ell$, we write 
$\partial_{\bo{u},\bo{v}}^+ \ell$ for the subpath of 
$\partial^+ \ell$ that starts at $\bo{u}$ and ends at $\bo{v}$.
Similarly, we define the path $\partial_{\bo{u},\bo{v}}^- \ell$.
\end{definition}
To define the successful action of $\sigma_{(\bo{u},\epsilon/2)}$,
note that the boundary $\partial W_{\bo{u},\epsilon/2}$ takes the form
 $\partial W_{\bo{u},\epsilon/2} = \ell_1 \cup \ell_2$, where
\begin{equation}\label{ellonedef}
\ell_1 = \Big\{ \bo{x} \in \R^2: \arg \big( \bo{x} \big)
  = \arg \big( \bo{u} \big) - \epsilon/2  \Big\}
\end{equation}
and
\begin{equation}\label{elltwodef}
\ell_2 = \Big\{ \bo{x} \in \R^2: \arg \big( \bo{x} \big)
  = \arg \big( \bo{u} \big) + \epsilon/2  \Big\}.
\end{equation}
Let $\xo \in \partial^+ \ell_1$ and 
$\yo \in \partial^- \ell_2$ respectively attain the minimal norm of the
sets $\partial^+ \ell_1 \cap \vcir$ and
 $\partial^- \ell_2 \cap \vcir$. Let $\bo{x_1} \in \partial^+ \ell_1$
have minimal norm subject to
\begin{equation}\label{defxone}
\vert\vert  \bo{x_1}  \vert\vert
 \geq \sup \Big\{ \vert\vert  \bo{x} \vert\vert: 
  \bo{x} \in \partial^+ \ell_1 \cap \vcir \Big\} + 2 \clemlu \epsilon n/C_0
\end{equation}
and let
$\bo{y_1} \in \partial^- \ell_2$
have minimal norm subject to
\begin{equation}\label{defxtwo}
\vert\vert  \bo{y_1}  \vert\vert
 \geq \sup \Big\{ \vert\vert  \bo{x} \vert\vert: 
  \bo{x} \in \partial^- \ell_2 \cap \vcir \Big\} + 2 \clemlu \epsilon n/C_0,
\end{equation}
where the constant $\clemlu$ will be specified in the upcoming Lemma \ref{lemlu}. 
\begin{definition}
Let $\bo{r} = \bo{r}_{\bo{u},\epsilon/2} \in \Z^2$ denote the  point of intersection of 
$\partial^+ \ell_1$ and $\partial^- \ell_2$ of maximal norm. Let $\Delta_{\bo{u},\epsilon/2} \subseteq \R^2$ denote the open region inside $W_{\bo{u},\epsilon/2}$ satisfying 
$\partial \Delta_{\bo{u},\epsilon/2} = \partial_{\bo{r},\infty}^+ \ell_1 \cup 
 \partial_{\bo{r},\infty}^- \ell_2$.
\end{definition}
Recall the quantity ${\rm fluc}_{\bo{x},\bo{y}}\big( \cdot \big)$ from Definition \ref{deffluc}. Let $\ccthr$ be a constant satisfying $\ccthr \geq 4 \clemlu^{-1} \csc \big( \qzero/2 \big)$. The operation is said to act successfully if
the updated 
configuration $\omega_1 \big\vert_{\csece}$
realizes the event
\begin{eqnarray}
 & & \Big\{ \partial^+_{\bo{x_0},\bo{x_1}}  \ell_1  \cup \partial^-_{\bo{y_0},\bo{y_1}}  \ell_2  \, \, \textrm{is open}
 \Big\} \cap \Big\{ \bo{x_1} \build\leftrightarrow_{}^{\Delta_{\bo{u},\epsilon/2} \cup \{ \bo{x_1} \} \cup \{ \bo{y_1} \}} \bo{y_1} \Big\} \nonumber \\
  & &  \cap \Big\{ {\rm fluc}_{\bo{x_1},\bo{y_1}} \big( \overline\gamma_{\bo{x_1},\bo{y_1}} \big) \leq \frac{\epsilon n}{\ccthr C_0} \Big\} 
       \cap \Big\{   \overline\gamma_{\bo{x_1},\bo{y_1}} \subseteq \clu{\bo{x_1}} \cap \clum{\bo{y_1}} \Big\},  \label{rgsa}
\end{eqnarray}
where $\overline\gamma_{\bo{x_1},\bo{y_1}}$ denotes the common $\omega_1$-open component of $\bo{x_1}$ and $\bo{y_1}$ in $\Delta_{\bo{u},\epsilon/2} \cup \big\{ \bo{x_1} \big\} \cup \big\{ \bo{y_1} \big\}$.
(Note that, since the open set $\Delta_{\bo{u},\epsilon/2}$ is disjoint from   
the boundary $\partial^+ \ell_1 \cup \partial^- \ell_2$, $\overline\gamma_{\bo{x_1},\bo{y_1}}$ is a connected set  disjoint from this boundary except at $\bo{x_1}$ and $\bo{y_1}$. 
This definition of  $\overline\gamma_{\bo{x_1},\bo{y_1}}$  is adopted in order that, in the definition of successful action, 
the opening of  $\partial^+_{\bo{x_0},\bo{x_1}} \big( \ell_1 \big) \cup \partial^-_{\bo{y_0},\bo{y_1}} \big( \ell_2 \big)$ does not affect the value of the fluctuation ${\rm fluc}_{\bo{x_1},\bo{y_1}} \big(  \overline\gamma_{\bo{x_1},\bo{y_1}}  \big)$.)  See the second picture in Figure 4.
%\ref{figmacregbox}.
\begin{definition}\label{defepsstar}
Let $\epsilon^* = \epsilon - \conepsstar n^{-1} \log n$, where $\conepsstar > 0$ is a constant to be specified later. 
Set
$$
\ell^*_1 = \Big\{ \bo{x} \in \R^2: \arg \big( \bo{x} \big)
  = \arg \big( \bo{u} \big) + \epsilon^*/2  \Big\}
$$
and
$$
\ell^*_2 = \Big\{ \bo{x} \in \R^2: \arg \big( \bo{x} \big)
  = \arg \big( \bo{u} \big) - \epsilon^*/2  \Big\}.
$$
Let $\bo{u}_0, \bo{u}_1$ attain the minimal and maximal norms among
elements of the set 
$\partial^- \ell_1^* \cap \vcir$. Let $\bo{v_0}$ and $\bo{v_1}$
play the analogous role for the set $\partial^+ \ell_2^* \cap \vcir$.
Set $P_1 = \partial_{\bo{u_0},\bo{u_1}}^- \ell_1^*$ and
 $P_2 = \partial_{\bo{v_0},\bo{v_1}}^+ \ell_2^*$.  
\end{definition}
\noindent{\bf Output properties.}
Recall the notation  $E^*(A)$ for $A \subseteq \R^2$
from Definition \ref{defpathedge} in Section \ref{secnot}.
Fix $\phi \subseteq E^*(B_K)$
that contains a path from $\partial B_K \cap W_{-(\bo{v_0} - \bo{u_0}),\qzero/2}\big( \bo{0} \big)$ to  $\partial B_K \cap W_{\bo{v_0} - \bo{u_0},\qzero/2}\big( \bo{0} \big)$ and satisfies 
$\phi  \cap \big( W_{\bo{u_0},2\qzero}\big( \bo{0} \big) \cup  W_{- \bo{u_0}, 2\qzero} \big( \bo{0} \big)  \big) = \emptyset$. 
That such a $\phi$ exists is ensured by 
\begin{equation}\label{angvouo}
  \ang \big( \bo{v_0} - \bo{u_0} , \perpu{\bo{u_0}} \big) < \pi/2 - 3\qzero/2.
\end{equation}
We defer verifying (\ref{angvouo}) 
for now, since we will anyway obtain it in due course.

We make the claim that, if the input configuration $\omega$ is
satisfactory,
 and if the operation $\sigmac$ acts successfully, then the output 
$\sigmac(\omega) = \big( \omega_1,\omega_2 \big)$ enjoys the following two 
properties:
\begin{enumerate}
\item {\bf Full-plane circuit property:} the full-plane configuration $\omega_1 \in \zoz$ contains an open circuit $\Gamma$ satisfying 
$\Gamma \subseteq B_{\cfpo n}$, with $\cfpo = \csc \big( \qzero/2 \big) c_0 (\cctwo + 1)$,  and 
 $\big\vert  {\rm INT}  \big( \Gamma \big) \big\vert \geq n^2$,
\item {\bf Sector open-set property:}
and the sector configuration $\omega_2 \in \big\{ 0,1
  \big\}^{\csece}$ realizes the event, to be denoted by ${\rm SOSP}$, that
there exists a set $\gamma$ that is a union of $\omega_2$-open edges, with 
$E(\gamma) \subseteq E^* \big(  W_{\bo{u},\epsilon^*/2}(\bo{0}) \big)$ and $\gamma \subseteq
B_{\cctwo n} \setminus B_{\ccone n}$, and is
such that the set $\gamma' : = \gamma \cup P_1 \cup P_2$
lies in $W_{\bo{u},\epsilon/2} (\bo{0})$, is connected, 
satisfies $\{ \bo{u_0} \} \cup \{ \bo{v_0} \} \subseteq V \big( \gamma'  \big)$
and 
\begin{equation}\label{rgcemp}
{\rm CRG}_{\qzero/2,K,\phi}^{\bo{u_0},\bo{v_0}} \big(
 \gamma'  \big) \cap W_{\bo{u},\epsilon/4} \big( \bo{0} \big)
 = \emptyset.
\end{equation}
\end{enumerate}
\begin{figure}\label{figmacregbox}
\begin{center}
\includegraphics[width=0.6\textwidth]{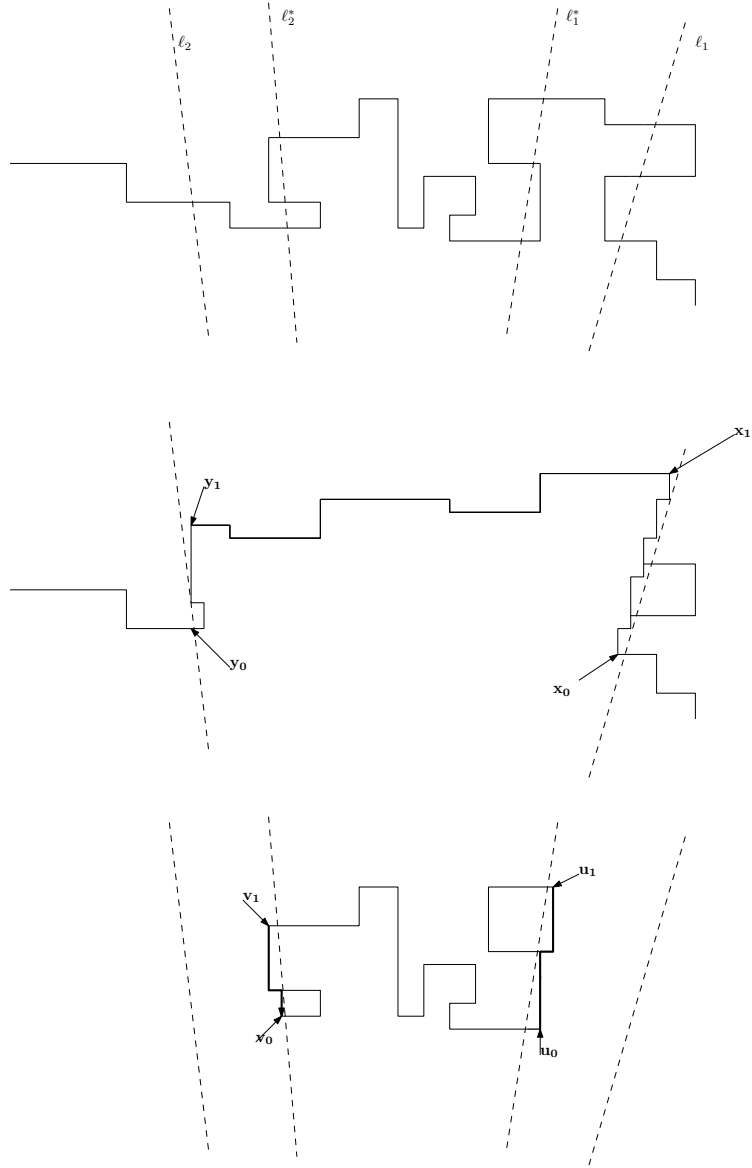} \\
\end{center}
\caption{The first figure depicts $\cir$ in a neighbourhood of $\csec$ from a satisfactory input. In the second, the same region of the full-plane output is depicted, after successful action. Note how the opening of $\partial_{\xo,\bo{x_1}}^+ \big( \ell_1 \big) \cup \partial_{\yo,\bo{y_1}}^- \big( \ell_2 \big)$ glues the new crossing path to the existing circuit. In the third figure, the set $\gamma$ in the sector output is shown, with the paths $P_1$ and $P_2$ (depicted in bold) joining up its components to form the connected set $\gamma'$.}
\end{figure}
\noindent{\bf Proof of the sector open-set property.} 
Note that 
$\omega_2 \in \{0,1\}^{\csece}$ 
satisfies
$\omega_2 = \omega \big\vert_{\csece}$, where the input $\omega$ realizes ${\rm SAT}$.
We take $\gamma$ so that its edge-set is given by
$$
E\big( \gamma \big)  = E\big(\cir\big) \cap E^* \Big( W_{\bo{u},\epsilon^*/2}\big( \bo{0} \big) \Big).
$$
By $\gamma \subseteq \cir$, we certainly have $\gamma \subseteq B_{\cctwo n} \setminus B_{\ccone n}$.
To see that $\gamma'$ is connected,
note firstly that:
\begin{lemma}\label{lemtopl}
If $\bo{v} \in \vcir$ satisfies $\arg(\bo{v}) < \arg(\bo{u}) -
\epsilon^*/2$, and there exists $\bo{w} \in \Z^2$, $\vert\vert \bo{w} -
\bo{v} \vert\vert = 1$ such that $\arg(\bo{w}) \geq \arg(\bo{u}) -
\epsilon^*/2$, then $\bo{v} \in P_1$.
\end{lemma}
\noindent{\bf Proof.}
Note that if a semi-infinite line $\ell$ (emanating from $\bo{0}$) and a nearest-neighbour edge 
$[\bo{w},\bo{w'}]$ intersect, then $\bo{w} \in \partial^+ \ell \cup \partial^- \ell$.
From this, $P_1
= \partial_{\bo{u_0},\bo{u_1}}^- \ell_1^*$, and the definition of $\bo{u_0}$ and
$\bo{u_1}$, follows the statement. \qed

Note now that $\bo{0} \in {\rm INT} \big( \cir \big)$
implies that $\bo{0}$ is disconnected in $W_{\bo{u},\epsilon^*/2}\big( \bo{0}
\big)$   from $\infty$
by $\cir$. However,
$$
  \cir \cap W_{\bo{u},\epsilon^*/2} \big( \bo{0}
 \big) \subseteq \gamma,
$$
so that $\bo{0}$ is disconnected from $\infty$
by $\gamma$ in $W_{\bo{u},\epsilon^*/2}\big( \bo{0}
\big)$. This forces $\gamma$ to contain a path connecting the
two lines $\ell_1^*$ and $\ell_2^*$ comprising $\partial W_{\bo{u},\epsilon^*/2}\big( \bo{0} \big)$.
The minimal nearest-neighbour path in which this path is contained necessarily
connects $P_1$ and $P_2$, by Lemma \ref{lemtopl} and its counterpart for $P_2$.
Thus, a path in $\gamma'$ runs between the connected sets $P_1,P_2 \subseteq \gamma'$. To establish that $\gamma'$ is connected, it remains to show that, from any $\bo{x} \in \gamma$ runs a path in $\gamma$ to either $P_1$ or $P_2$. However, since $\bo{x} \in \cir$ and $\bo{0} \in \intg$, there exists a nearest-neighbour path with edge-set contained in  $E\big(\cir\big) \cap E^* \big( W_{\bo{u},\epsilon/2} \big)$ that begins at 
$\bo{x}$ and whose final edge (but no other) crosses $\partial W_{\bo{u},\epsilon/2}$. 
By the definition of $\gamma$, this path lies in $\gamma$; and its final vertex lies on $P_1 \cup P_2$, by 
Lemma \ref{lemtopl}. 

It remains to verify (\ref{rgcemp}). Recall that $\omega \in {\rm SAT}$ implies that
$\csec \cap \reg = \emptyset$. It suffices then to establish that
\begin{equation}\label{rwrg}
 {\rm CRG}_{\qzero/2,K,\phi}^{\bo{u_0},\bo{v_0}} \big(
  \gamma'   \big) \cap W_{\bo{u},\epsilon/4} \big( \bo{0} \big) 
 \subseteq \reg,
\end{equation}
where we identify the space $\csec$ in which $\gamma'$ is contained with the corresponding subset in the copy of $\R^2$ that contains $\cir$. In deriving (\ref{rwrg}), we will also obtain (\ref{angvouo}). The derivation takes some time, and we prefer to give it at the end of the proof. 
Admitting (\ref{rwrg}) for now, the output indeed satisfies the sector open-set property. \\
\noindent{\bf Proof of the full-plane circuit property.}
%Let $\gamma$ denote the connected component of $\bo{x_1}$
%in $\csec$ of the replacing configuration $\sigma' \in \{ 0,1 \}^{\csece}$.
Recall the set  $\overline\gamma_{\bo{x_1},\bo{y_1}}$ defined after (\ref{rgsa}).
We define a set $\tilde\Gamma$ so that its edge-set is given by 
$$
E \big( \tilde\Gamma \big) =  \Big( E \big( \Gamma_0 \big) \cap \csece^c  \Big) 
\cup E \big( \overline\gamma_{\bo{x_1},\bo{y_1}} \big) \cup E \big( \partial_{\xo,\bo{x_1}}^+(\ell_1) \big) \cup E \big( \partial_{\yo,\bo{y_1}}^-(\ell_2) \big).
$$
Let $\Gamma$ denote the outermost circuit of $\tilde\Gamma$. Note that
$\big\vert {\rm INT} \big( \Gamma \big) \big\vert = \big\vert {\rm INT} \big(  \tilde\Gamma \big) \big\vert$. 

To show that $\tilde{\Gamma} \subseteq B_{\cfpo n}$, note that, by $\cir \subseteq B_{\cctwo n}$, 
$\partial_{\xo,\bo{x_1}}^+(\ell_1) \subseteq B_{\vert\vert \bo{x_1} \vert\vert}$, $\partial_{\yo,\bo{y_1}}^+(\ell_1) \subseteq B_{\vert\vert \bo{y_1} \vert\vert}$ and
$\bo{x_1}, \bo{y_1} \in \overline\gamma_{\bo{x_1},\bo{y_1}}$,   
we need only confirm that $\overline\gamma_{\bo{x_1},\bo{y_1}}  \subseteq B_{\cfpo n}$. 
Recall that $\overline\gamma_{\bo{x_1},\bo{y_1}}    \subseteq \clu{\bo{x_1}} \cap \csec \subseteq \clu{\bo{x_1}} \cap W^-_{\bo{x_1},\epsilon}$. Imposing $\epsilon < c_0$, Lemma \ref{lemdistang} implies that the maximal norm of a point in $\clu{\bo{x_1}} \cap  W^+_{\bo{x_1},\epsilon}$ is at most $\csc\big( \qzero/2 \big) \vert\vert \bo{x_1} \vert\vert c_0$. 
The definition of $\bo{x_1}$, and $\cir \subseteq B_{\cctwo n}$, entail that $\vert\vert \bo{x_1} \vert\vert \leq \cctwo n + \epsilon n/C_0 + 1$, so that $\tilde{\Gamma} \subseteq B_{\cfpo n}$ 
follows from our assumption that $C_0 > 2 c_0$.
%We replace $\tilde\Gamma$ by the outermost dual circuit that it contains. [Check that $V(\tilde\Gamma) \subseteq B_{2Cn}$.]]

%Let $\Delta$ denote the infinite component of $\R^2 \setminus 
%\big( \partial^+ \ell_1 \cup \partial^- \ell_2 \big)$
%that is contained in $W_{\bo{u},\epsilon/2}$.
%Let $\tilde{T}_{\bo{0},\bo{x_1},\bo{y_1}}$ denote the bounded component of $\Delta \setminus \big[ \bo{x_1},\bo{y_1} \big]$.

Let $\ell'$ denote the line parallel to $\ell_{\bo{x_1},\bo{y_1}}$
that has a displacement from $\ell_{\bo{x_1},\bo{y_1}}$ by $\frac{\epsilon n}{\ccthr C_0}$, in the direction of the origin.  Set $\Delta =\Delta_{\bo{u},\epsilon/2}$ and 
let $\tilde{T}$ denote the subset of $\Delta$ lying on the same side of $\ell'$ as does $\bo{0}$.

We now show that
$$
 \Big\vert {\rm INT} (\tilde\Gamma) \Big\vert = 
 \Big\vert {\rm INT} (\tilde\Gamma) \cap \Delta \Big\vert
 +   \Big\vert {\rm INT} (\tilde\Gamma) \cap \Delta^c \Big\vert, \quad \textrm{and likewise for $\cir$,}
$$
\begin{equation}\label{intgr}
 {\rm INT} \big( \cir  \big) \cap \Delta^c \subseteq   
 {\rm INT} \big( \tilde\Gamma \big) \cap \Delta^c,
\end{equation}
\begin{equation}\label{inttbd}
 \tilde{T} \subseteq  {\rm INT} \big( \tilde\Gamma \big) \cap \Delta,
\end{equation}
\begin{equation}\label{opintbd}
  {\rm INT} \big( \cir \big) \cap \Delta 
 \subseteq  \tilde{T},
\end{equation}
whence
$$
 \big\vert {\rm INT} \big( \tilde\Gamma \big) \big\vert
 \geq \big\vert {\rm INT} \big( \cir \big) \big\vert \geq n^2,
$$
as required for the full-plane circuit property. 
%We establish (\ref{intgr}) similarly to (3.24) in 
%intinc
%\cite{hammondone}. 
The main ingredients for establishing (\ref{intgr}) are the assertions $\cir \cap \overline{\Delta}^c = \tilde\Gamma \cap \overline{\Delta}^c$, (where $\overline{\Delta}$ denotes the closure of $\Delta$), 
\begin{equation}\label{auxeq}
{\rm INT} \big( \tilde\Gamma \big) \cap \big( \partial^+ \ell_1 \cup \partial^- \ell_2 \big) \supset \partial_{\bo{0},\bo{x_1}}^+ \ell_1 \cup \partial_{\bo{0},\bo{y_1}}^- \ell_2 
\end{equation}
and  
${\rm INT} \big( \cir \big) \cap \big( \partial^+ \ell_1 \cup \partial^- \ell_2 \big) \subset \partial_{\bo{0},\bo{x_1}}^+ \ell_1 \cup \partial_{\bo{0},\bo{y_1}}^- \ell_2$. 
From these, (\ref{intgr}) follows readily. 
%(For more detail on this point, see the derivation of (3.24) in \cite{hammondone}.)
To establish (\ref{auxeq}), note that $\tilde\Gamma$ separates 
$\partial^+_{\bo{0},\bo{x_0}} \ell_1 \cup \partial^-_{\bo{0},\bo{y_0}} \ell_2$
and
$\partial^+_{\bo{x_1},\infty} \ell_1 \cup \partial^-_{\bo{y_1},\infty} \ell_2$
both in $\Delta$ and in $\Delta^c$.
This implies that 
$\partial^+_{\bo{0},\bo{x_0}} \ell_1 \cup \partial^-_{\bo{0},\bo{y_0}} \ell_2 \subseteq {\rm INT} \big( \tilde\Gamma \big)$. We use this alongside 
$\partial^+_{\bo{x_0},\bo{x_1}} \ell_1 \cup \partial^-_{\bo{y_0},\bo{y_1}} \ell_2  \subseteq \tilde\Gamma$ to obtain (\ref{auxeq}).
 
%To see (\ref{intgr}), suppose that $\bo{x} \in \big( \R^2 \setminus \csec
%\big) \setminus {\rm INT} \big( \tilde\Gamma \big)$. We must show that 
%$\bo{x} \not\in \intg$. Let $\phi:[0,\infty) \to \R^2$ satisfy $\phi(\bo{0}) =
%\bo{x}$, $\vert\vert \phi(t) \vert\vert \to \infty$ as $t \to \infty$, and $\phi[0,\infty) \cap \tilde\Gamma = \emptyset$. If
%$\phi[0,\infty) \subseteq \R^2 \setminus \csec$, then $\phi [0,\infty) \cap
%\cir = \emptyset$, since $\tilde\Gamma$ and $\Gamma_0$ 
%coincide in $\R^2 \setminus \csec$. Note that ${\rm INT} \big(
%\tilde\Gamma \big) \cap \partial \csec$ contains 
%$\big[ \bo{0} , \bo{r_1} \big] \cup \big[ \bo{0} , \bo{r_2} \big]$,
%where $\bo{r_1}$ and $\bo{r_2}$ are the points of intersection of $E \big(
%\Gamma_0 \big)$ with $\ell_1$ and $\ell_2$ of maximal norm. [Justify.]
%Hence, if $\phi[0,\infty) \cap \big( \R^2 \setminus \csec \big) \not=
%\emptyset$, then
%$$
%\phi[0,\infty) \cap \partial \csec \subseteq \big\{ t \bo{r_1}: t \geq 1
%\big\} \cup \big\{ t \bo{r_2}: t \geq 1  \big\}, 
%$$
%since $\phi[0,\infty) \cap \tilde\Gamma = \emptyset$. However,
%$\big\{ t \bo{r_1}: t \geq 1  \big\} \cup \big\{ t \bo{r_2}: t \geq 1
%\big\} \subseteq \intg^c$, so that $\phi$ may be modified to demonstrate
%that $\bo{x} \not\in \intg$, as required for (\ref{intgr}).

To derive (\ref{inttbd}), note that ${\rm INT} \big( \tilde\Gamma \big) \cap \Delta = {\rm INT} \big( \overline{\gamma}_{\bo{x_1},\bo{y_1}}  \cup \partial^+_{\bo{0},\bo{x_1}} \ell_1
\cup  \partial^-_{\bo{0},\bo{y_1}} \ell_2 \big)$ by (\ref{auxeq}). 
Hence, (\ref{inttbd}) follows from 
$\overline\gamma_{\bo{x_1},\bo{y_1}} \cap \ell' = \emptyset$, which is implied by 
${\rm fluc}_{\bo{x_1},\bo{y_1}} \big( \overline\gamma_{\bo{x_1},\bo{y_1}} \big) \leq \frac{\epsilon n}{\ccthr C_0}$.

To justify (\ref{opintbd}), we require the following lemma.
\begin{lemma}\label{lemlu}
Recall that $\tcir = n \partial \wulff$ denotes the dilation of the Wulff curve attaining $\globdis$. For each $\bo{u} \in S^1$, write 
$L_{\bo{u}} = \ell_{\bo{0},\bo{u}}^+ \cap \big( \tcir + \globdis \big)$.
There exists $\clemlu > 0$ such that, for all $\bo{u} \in S^1$, 
$\big\vert L_{\bo{u}} \big\vert \leq \clemlu \globdis$.
\end{lemma}
In essence, Lemma \ref{lemlu} follows from the tangent vectors on the Wulff shape boundary never pointing too directly towards the origin (\ref{supang}). The details are perhaps a little distracting at this stage, so that we give the proof of Lemma \ref{lemlu} after the end of the proof. 

For (\ref{opintbd}), note that it suffices to show that $\cir \cap \csec \subseteq H$, where $H$ denotes the half-plane with $\partial H = \ell_1$ that contains $\bo{0}$. Note that $\tcir = n \partial \wulff$ satisfies 
\begin{equation}\label{eqcirtcir}
\cir \subseteq \tcir + B_{\globdis}, 
\end{equation}
%a choice possible by the Definition \ref{defglobdis} of $\globdis$ and 
since $\centre(\cir) = \bo{0}$. 
Let $\mathcal{C}^*$ denote the intersection of the outer boundary of 
$\tcir + B_{\globdis}$ with $\csec$, and write $\bo{x^*} = \mathcal{C}^* \cap \ell_1$ and $\bo{y^*} = \mathcal{C}^* \cap \ell_2$ for the endpoints of $\mathcal{C}^*$. For (\ref{eqcirtcir}), it suffices to show that $\mathcal{C}^* \subseteq H$. 
We parameterize the lines $\ell_{\bo{x^*},\bo{y^*}}$,  $\ell_{\bo{x_1},\bo{y_1}}$, $\ell'$ and $\mathcal{C}^*$ according to the argument-value $\theta \in \big[ \argu(\bo{u})- \epsilon/2 , \argu(\bo{u}) + \epsilon/2  \big]$. (For the rest of the derivation of (\ref{opintbd}), statements containing $\theta$ are asserted for all $\theta$ in this interval.)  It is our aim to show that $\mathcal{C}^*(\theta) \leq \ell'(\theta)$.
This will follow from the bounds
\begin{equation}\label{eqellone}
  \ell_{\bo{x^*},\bo{y^*}} \big( \theta \big) \leq \ell_{\bo{x_1},\bo{y_1}} \big( \theta \big)
 -  \frac{\clemlu}{2 C_0} \epsilon n,
\end{equation}
\begin{equation}\label{eqelltwo}
  \mathcal{C}^*\big( \theta \big) \leq \ell_{\bo{x^*},\bo{y^*}} \big( \theta \big)
 +  \ccfour \epsilon^2 n,
\end{equation}
where $\ccfour$ is a constant depending only on $\partial \wulff$, and
\begin{equation}\label{eqellthr}
  \ell' \big( \theta \big) \leq \ell_{\bo{x_1},\bo{y_1}} \big( \theta \big)
 -  \csc\big( \qzero/2 \big) \frac{\epsilon n}{C_0}.
\end{equation}
To derive (\ref{eqellone}), note that
$\max \big\{ d(\xo,\bo{x^*}),d(\yo,\bo{y^*})  \big\} \leq \clemlu C_0^{-1} \epsilon$
follows from $\globdis \leq \epsilon n/C_0$ and Lemma \ref{lemlu}.
We have that 
$\min \big\{ d(\bo{x_0},\bo{x_1}),d(\bo{y_0},\bo{y_1})  \big\} \geq 2 \clemlu C_0^{-1} \epsilon n$
by the definition of the points in question. From these inequalities, it follows that 
$\min \big\{ d(\bo{x^*},\bo{x_1}),d(\bo{y^*},\bo{y_1})  \big\} \geq 2 \clemlu C_0^{-1} \epsilon n$,
which provides (\ref{eqellone}) for the endpoint choices of $\theta$. It is easily seen that the quantity 
$\ell_{\bo{x_1},\bo{y_1}}(\theta) - \ell_{\bo{x^*},\bo{y^*}}(\theta)$ is at least one-half of the minimum of its endpoint values (since $\epsilon$ is taken to be small). Hence, we obtain (\ref{eqellone}).

%To obtain (\ref{eqelltwo}), note that the dilation factor $t$ in $\tcir$ satisfies $t \leq 2 \cctwo \cwinf n$,
%where $\cwinf = \inf \big\{ \vert\vert \bo{x} \vert\vert: \bo{x} \in \partial W \big\}$.
%This is a consequence of $\cir \subseteq B_{\cctwo n}$ and $\globdis \leq \epsilon n/C_0 \leq \cctwo n$. By this bound on $t$, and the second-order differentiability of $\partial \wulff$ (by Lemma \ref{lemozstr}), (\ref{eqelltwo}) readily follows. 

To obtain (\ref{eqelltwo}), note that $\tcir = n \partial \wulff$.  By the second-order differentiability of $\partial \wulff$ (from Lemma \ref{lemozstr}), (\ref{eqelltwo}) readily follows. 

For (\ref{eqellthr}), set $\psi = \ang \big( \bo{y_1},\bo{x_1}, \bo{x_1^\perp} \big)$. It is easy to see that
$\ell_{\bo{x_1},\bo{y_1}}(\theta) - \ell'(\theta) \geq \sec \big( \psi + \epsilon \big) d \big( \ell', \ell_{\bo{x_1},\bo{y_1}} \big)$, so that (\ref{eqellthr}) follows from 
\begin{equation}\label{xyyxinc}
\bo{y_1} \in \clu{\bo{x_1}} \qquad \textrm{and} \qquad \bo{x_1} \in \clum{\bo{y_1}},
\end{equation}
$\epsilon < \qzero/2$, $d \big( \ell', \ell_{\bo{x_1},\bo{y_1}} \big) = \frac{\epsilon n}{\ccthr C_0}$ and $\ccthr \geq 4 \clemlu^{-1} \csc \big( \qzero/2 \big)$.

It remains only to verify (\ref{xyyxinc}). By minor adjustments to the derivation of the upcoming (\ref{angcon}), we find that  $\ang \big( \yo - \xo , \perpu{\xo} \big) \leq \pi/2 - 3 \qzero/2$. By $\globdis \leq \epsilon n/C_0$, Lemma \ref{lemlu} and the definitions (\ref{defxone}) and (\ref{defxtwo}) of $\bo{x_1}$ and $\bo{y_1}$, we have that
$\vert\vert \bo{x_1} - \xo \vert\vert, \vert\vert \bo{y_1} - \yo \vert\vert \leq \clemlu \epsilon n/c_0 + \clemlu \epsilon n/c_0  = 2 \clemlu \epsilon n/c_0$. From $\xo,\yo \in \cir$, $\cir \cap B_{\ccone n} = \emptyset$,
and $\ang(\xo,\yo)  \geq \epsilon/2$, we obtain  $\vert\vert \bo{x_0} - \yo \vert\vert \geq  \pi^{-1} c_1 n \epsilon$. Hence, for $C_0$ a high enough constant,  
$\ang \big( \bo{y_1} - \bo{x_1} , \perpu{\bo{x_1}} \big) \leq \pi/2 - \qzero$. 
Hence, we obtain (\ref{xyyxinc}),
since its second inclusion is similarly derived. This completes the proof of (\ref{opintbd}). \\
\noindent{\bf The upper bound on the probability of the two output properties.}
Note that, since the input configuration has law $P$, the full-plane configuration $\omega_1 \in \zoz$ in the output also has this law. Thus, the full-plane circuit property is satisfied by the output with probability equal to 
 \begin{eqnarray}
 & & P \Big( \exists \, \textrm{an open circuit} \, \, \Gamma: \Gamma \subseteq B_{\cfpo n}, 
  \big\vert {\rm INT} \big( \Gamma \big) \big\vert \geq n^2 \Big) \nonumber \\
   & \leq & \cpi \cfpo^2 n^2    P \Big( 
     \acon  \Big), \nonumber
 \end{eqnarray}
the inequality due to the following general device for centring circuits, (Lemma 3.2 in \cite{hammondone}):
\begin{lemma}\label{lemcirprop}
Let $A \in \N$. 
%Set $\cpi = 20\pi$. 
For any constant $\cgenbig > 0$,
$$
P \Big( \exists \,  \textrm{an open circuit} \, \, \Gamma:  \Gamma \subseteq B_{\cgenbig n},  
   \big\vert {\rm INT}(\Gamma) \big\vert \geq A  \Big) \leq
\cpi  \cgenbig^2 n^2 P \big(  \area{A} \big).
$$
\end{lemma}
Recall that we apply $\sigma_{(\bo{u},\epsilon/2)}$ so that it acts regularly. Hence,
%(\ref{ptomf}) 
the conditional probability of the sector open-set property, given the full-plane circuit property, is at most
\begin{equation}\label{omcsece}
 \sup \Big\{ P_{\tilde\omega} \big( {\rm SOSP} \big) : \tilde\omega \in \{ 0,1 \}^{E(\Z^2) \setminus \csece} \Big\}.
\end{equation}
(The paragraph leading to (3.45) in \cite{hammondone} explains this point.)
Note that
$$
{\rm SOSP} \in \sigma \Big\{ E^* \Big( W_{\bo{u},\epsilon^*/2} \big( \bo{0} \big) \Big)
\cap E \Big( B_{2Cn} \setminus B_{\ccone n} \Big) \Big\}.
$$
The edge-set on the right-hand-side is at distance at least $\pi^{-1} \ccone \conepsstar \log n$ from $E\big(\Z^2\big) \setminus \csece$, 
where the constant $\conepsstar$ appears in the
Definition \ref{defepsstar} of $\epsilon^*$. 
By tuning this constant to be high enough, this edge-set and 
$E\big(\Z^2\big) \setminus \csece$ are $(\clemkam,0)$-well separated, 
in the sense of Lemma \ref{lemkapab}. 
Applying that lemma, 
we find that (\ref{omcsece}) is at most $\conka P \big( R \big)$.

Note that
the occurrence of ${\rm SOSP} \cap \big\{ E(P_1) \cup E(P_2) \, \, \textrm{is open}
\big\}$ entails the presence of a connected open set $\gamma' \subseteq B_{\ccone n}^c$
 for which $\bo{u_0},\bo{v_0} \in V(\gamma')$ and
 $$
 {\rm CRG}_{\qzero/2,K,\phi}^{\bo{u_0},\bo{v_0}} \big( \gamma' \big) \cap \csech = \emptyset.
 $$
The two straight boundary segments of $\csec \cap B_{\ccone n}^c$ being at distance at least $\epsilon \ccone \pi^{-1} n$, we see that ${\rm MAXREG}_{\qzero/2,K,\phi}^{\bo{u_0},\bo{v_0}} \big( \gamma' \big) \geq \epsilon \ccone \pi^{-1} n$.
 By Lemma \ref{lemmaxreg}, we obtain that the probability under $P$ that such a
 $\gamma'$ exists
 is at most 
$P \big( \bo{u_0} \leftrightarrow \bo{v_0} \big) \exp \big\{ - c \ccone \pi^{-1} \epsilon n \big\}$, so that
 $$
 P \Big( 
{\rm SOSP} \cap \Big\{ E(P_1) \cup E(P_2) \, \, \textrm{is open}
\Big\} \Big)
 \leq P \Big( \bo{u_0} \leftrightarrow \bo{v_0} \Big) \exp \big\{ - c \ccone \pi^{-1} \epsilon n \big\}.
 $$
 However, the conditional probability that the set 
  $E(P_1) \cup E(P_2)$ is open, given that ${\rm SOSP}$ occurs, is at least $\cposen^{\vert E(P_1) \vert + \vert E(P_2) \vert}$, by conditioning on the set $\gamma$ appearing in the definition of ${\rm SOSP}$, and by the bounded energy property of $P$. Hence, we have use for:
  \begin{lemma}\label{lemedgep}
   $$
    \big\vert E(P_1) \big\vert + \big\vert E(P_2) \big\vert \leq \frac{5 \clemlu \epsilon n}{C_0}.
   $$
  \end{lemma}
This lemma is very similar to Lemma \ref{lemlu}, and its proof appears at the end of this section.
By applying it, we find that
$$
P \Big( {\rm SOSP}  \Big) \leq  \cposen^{- \frac{5 \clemlu}{C_0} \epsilon n} P \Big( \bo{u_0} \leftrightarrow \bo{v_0} \Big)
 \exp \Big\{ - c  \ccone \pi^{-1} \epsilon n \Big\}. 
$$
By fixing $C_0 > 0$ high enough, then, 
$P \big( {\rm SOSP}  \big) \leq  P \Big( \bo{u_0} \leftrightarrow \bo{v_0} \Big)
 \exp \big\{ -  (c/2) \ccone \pi^{-1}  \epsilon n \big\}$. 

In summary, in acting on an input with law $P$, $\sigma_{(\bo{u},\epsilon/2)}$ will return an output having the full-plane circuit and sector open-set properties with probability at most 
\begin{equation}\label{clsatbd}
 \cpi \cfpo^2 n^2   
 P \Big( \acon \Big)
  \conka P \Big( \bo{u_0} \leftrightarrow  \bo{v_0} \Big)  \exp \Big\{ -  (c/2) \ccone \pi^{-1} 
  \epsilon n  \Big\}.
\end{equation}
\noindent{\bf The lower bound on the probability of satisfactory input and successful action.}
The operation $\sigmac$ acts on an input $\omega \in \zoz$
having the law $P$. 
Analogously to (3.37) in \cite{hammondone},
%(\ref{bdsatcn}), 
the input is satisfactory and the operation acts successfully with probability at least
$$
 P \Big( {\rm SAT} \Big) \inf \Big\{ P_{\tilde\omega}\big( \ref{rgsa} \big): \tilde\omega \in \{ 0,1\}^{E(\Z^2 \setminus \csece)} \Big\}.
$$ 
We will use:
\begin{lemma}\label{lemgaco}
Let $\bo{u} \in S^1$, $\epsilon > 0$, 
$\bo{z_1} \in \partial^+ \ell_1$ and 
$\bo{z_2} \in \partial^- \ell_2$, where $\ell_1$ and $\ell_2$ are defined in (\ref{ellonedef}) and (\ref{elltwodef}). Assume that $\bo{z_1}$ and $\bo{z_2}$ satisfy the hypotheses of Lemma \ref{lemgac}. Let $\omega \in \{ 0,1 \}^{E(\Z^2) \setminus E(W_{\bo{u},\epsilon/2})}$. Then, for $c > 0$ sufficiently small,
\begin{eqnarray}
 & & P_\omega \Big( \Big\{ \bo{z_1} 
\build\leftrightarrow_{}^{\Delta_{\bo{u},\epsilon/2} \cup \{ \bo{z_1} \} \cup \{ \bo{z_2} \} }  \bo{z_2} \Big\} \cap 
 \Big\{ {\rm fluc}_{\bo{z_1},\bo{z_2}} \big( \overline\gamma_{\bo{z_1},\bo{z_2}} \big) \geq c \vert\vert  \bo{z_1} - \bo{z_2}  \vert\vert \Big\}
   \label{deloutfl} \\
 & &  \qquad \qquad \cap \Big\{  \overline\gamma_{\bo{z_1},\bo{z_2}}  \subseteq \clu{\bo{z_1}} \cap \clum{\bo{z_2}} \Big\} \Big) 
 \leq \exp \Big\{ - c \vert\vert    \bo{z_1} - \bo{z_2} \vert\vert \Big\} P  \Big(  \bo{z_1} \build\leftrightarrow_{}^{}  \bo{z_2} \Big), \nonumber
 \end{eqnarray}
 %{\rm GAC} \big( \bo{u}, \epsilon, \bo{z_1},\bo{z_2}\big) \Big) \geq n^{-C} P \Big( \bo{z_1} \build\leftrightarrow_{}^{\Delta_{\bo{u},\epsilon/2}}  \bo{z_2}  \Big).
%$$
where recall that $\overline\gamma_{\bo{z_1},\bo{z_2}}$ denotes the common open cluster of $\bo{z_1}$ and $\bo{z_2}$ in $\Delta_{\bo{u},\epsilon/2}  \cup \{ \bo{z_1} \} \cup \{ \bo{z_2} \}$.
\end{lemma}
\noindent{\bf Proof.} 
Note that the event on the left-hand-side of (\ref{deloutfl}) is measurable with respect to 
$\sigma \big\{ E \big( \clu{\bo{z_1}} \cap \clum{\bo{z_2}}  \big) \cap E \big( W_{\bo{u},\epsilon/2} \big) \big\}$. 
It is easy to see that there exist $\clemkam,\clemkac > 0$ such that the pair of edge-sets 
$E^* \big( \club{\bo{z_1}} \cap \clumb{\bo{z_2}}  \big) \cap 
E \big( W_{\bo{u},\epsilon/2} \big)$ 
and 
$E \big( \Z^2 \big) \setminus E \big( W_{\bo{u},\epsilon/2} \big)$ 
are 
$(\clemkam,\clemkac)$-well separated.
 Hence, by Lemma \ref{lemkapab}, it suffices to prove the variant of (\ref{deloutfl}) in which the measure $P_{\omega}$ on the left-hand-side is replaced by the unconditioned measure $P$. 
The required statement follows from Lemma \ref{lemmdf}. \qed
Note that, for any
$\omega \in \{ 0,1 \}^{E(\Z^2) \setminus \csece}$,
\begin{eqnarray}
 & & P_\omega \Big( \Big\{ \bo{x_1} \build\leftrightarrow_{}^{\Delta_{\bo{u},\epsilon/2}
  \cup \{ \bo{x_1} \} \cup \{ \bo{y_1} \} }  \bo{y_1} \Big\} 
  \cap \Big\{  \overline\gamma_{\bo{x_1},\bo{x_2}}  \subseteq \club{\bo{x_1}} \cap \clumb{\bo{y_1}} \Big\} \Big) \nonumber \\
 & \geq & c P \Big(  \bo{x_1} \build\leftrightarrow_{}^{\Delta_{\bo{u},\epsilon/2}
   \cup \{ \bo{x_1} \} \cup \{ \bo{y_1} \} }  \bo{y_1} \Big), \label{eqxoney}
\end{eqnarray}
the inequality due to the same of pair of edge-sets as in the preceding paragraph
being $(\clemkam,\clemkac)$-well separated,
and Lemma \ref{lemkapab} (to replace $P_\omega$ by $P$) and a use of Lemma \ref{lemoznor}. 
Note also that 
\begin{equation}\label{eqconest}
 \frac{P \big( \bo{x_1} \build\leftrightarrow_{}^{\Delta_{\bo{u},\epsilon/2}   \cup \{ \bo{x_1} \} \cup \{ \bo{y_1} \}}  \bo{y_1} \big)}{P \big( \bo{x_1} \leftrightarrow  \bo{y_1} \big)} \geq c. 
\end{equation}
Indeed, (\ref{xyyxinc}) implies that $W_{\bo{x_1} - \bo{y_1},\qzero}\big( \bo{y_1} \big) \cap W_{\bo{y_1} - \bo{x_1},\qzero}\big( \bo{x_1} \big) \subseteq A_{\bo{x_1},\bo{y_1}}$. Hence, any intersection between $W_{\bo{x_1} - \bo{y_1},\qzero}\big( \bo{y_1} \big) \cap W_{\bo{y_1} - \bo{x_1},\qzero}\big( \bo{x_1} \big)$
and $\Delta_{\bo{u},\epsilon}$ must occur within a bounded distance of 
$\{ \bo{x_1} \} \cup \{ \bo{y_1} \}$. From Lemma \ref{lemoznor}, we easily obtain (\ref{eqconest}).

By Lemma \ref{lemgaco}, (\ref{eqxoney}) and (\ref{eqconest}), 
we obtain
\begin{eqnarray}
 & & P_\omega \Big( \Big\{ \bo{x_1} \build\leftrightarrow_{}^{\Delta_{\bo{u},\epsilon/2} \cup \{ \bo{x_1} \} \cup \{ \bo{y_1} \} }  \bo{y_1} \Big\} \cap 
 \Big\{ {\rm fluc}_{\bo{x_1},\bo{y_1}} \big( \overline\gamma_{\bo{x_1},\bo{y_1}} \big) \leq \frac{\epsilon n}{\cthr C_0} \Big\} \nonumber \\
 & & \qquad \qquad \qquad 
  \cap \Big\{    \overline\gamma_{\bo{x_1},\bo{y_1}} \subseteq \clu{\bo{x_1}} \cap \clum{\bo{y_1}} \Big\} \Big) 
    \geq    \frac{c^2}{2} P \Big(  \bo{x_1} \build\leftrightarrow_{}^{\Delta_{\bo{u},\epsilon/2}  \cup \{ \bo{x_1} \} \cup \{ \bo{y_1} \}}  \bo{y_1} \Big). \nonumber
\end{eqnarray}
The remaining condition 
required for successful action, as stated in (\ref{rgsa}), namely that  $\partial^+_{\bo{x_0},\bo{x_1}} \big( \ell_1 \big) \cup \partial^-_{\bo{y_0},\bo{y_1}} \big( \ell_2 \big)$ be open, 
has conditional probability at least $\cposen^{\frac{6\clemlu \epsilon n}{C_0}}$, by 
the bounded energy property of $P$ and the bounds $\big\vert \partial^-_{\bo{x_0},\bo{x_1}} \big( \ell_1 \big) \big\vert, \big\vert \partial^-_{\bo{y_0},\bo{y_1}} \big( \ell_2 \big) \big\vert \leq 4 \clemlu \epsilon n/C_0 + 5 \clemlu \epsilon n/C_0$. 
To prove these upper bounds on 
$\big\vert \partial^-_{\bo{x_0},\bo{x_1}} \big( \ell_1 \big) \big\vert$ 
and  
$\big\vert \partial^-_{\bo{y_0},\bo{y_1}} \big( \ell_2 \big) \big\vert$, 
let 
$\bo{x_2}$ 
attain 
$\sup \big\{  \vert\vert  \bo{x} \vert\vert: \bo{x} \in \partial^- \big( \ell_1 \big) \cap \vcir \big\}$. 
Then 
$\big\vert \partial^-_{\bo{x_0},\bo{x_1}} \big( \ell_1 \big) \big\vert 
= \big\vert \partial^-_{\bo{x_0},\bo{x_2}} \big( \ell_1 \big) \big\vert  + \big\vert \partial^-_{\bo{x_2},\bo{x_1}} \big( \ell_1 \big) \big\vert$. A near-verbatim argument to the proof of  Lemma \ref{lemedgep} shows that  
$\big\vert \partial^-_{\bo{x_0},\bo{x_2}} \big( \ell_1 \big) \big\vert \leq 5 \clemlu \epsilon n/C_0$, and that 
$\big\vert \partial^-_{\bo{x_2},\bo{x_1}} \big( \ell_1 \big) \big\vert \leq 4 \clemlu \epsilon n/C_0$.
The bound on  $\big\vert \partial^+_{\bo{y_0},\bo{y_1}} \big( \ell_2 \big) \big\vert$ is identical.

Hence,
the input is satisfactory, and the operation acts
successfully, with probability at least 
\begin{equation}\label{csatbd}
   P \Big( {\rm SAT} \Big)  c \cdot \cposen^{6 \clemlu  \epsilon n/C_0}
  P \Big(  \bo{x_1} \build\leftrightarrow_{}^{\Delta_{\bo{u},\epsilon/2} \cup \{ \bo{x_1} \} \cup \{ \bo{y_1} \} }  \bo{y_1} \Big). 
 \end{equation}
\noindent{\bf Conclusion by comparison of the obtained bounds.} However, we have seen that circumstances arise in which such an output will definitely be produced whose probability is at least the quantity given in (\ref{clsatbd}). Thus, the quantity in (\ref{csatbd}) is at most that in (\ref{clsatbd}). That is,
\begin{eqnarray}
 P \Big( {\rm SAT} \Big) 
 & \leq  &   
%16 c^{-1} \pi C^2 n^2 C_1  
C
\cposen^{-6 \clemlu \epsilon n/C_0} \frac{P \big(  \bo{u_0} \leftrightarrow  \bo{v_0} \big)}{P \big(  \bo{x_1 } \build\leftrightarrow_{}^{\Delta_{\bo{u},\epsilon/2}  \cup \{ \bo{x_1} \} \cup \{ \bo{y_1} \} }  \bo{y_1} \big)}  \nonumber \\
 & &  \qquad 
 P \Big( \acon  \Big)   \exp \Big\{ -  (c/2) \ccone \pi^{-1} \epsilon n \Big\} .
\end{eqnarray}
 Noting that $d \big( \bo{u_0}, \bo{x_1} \big), d \big( \bo{v_0}, \bo{y_1} \big) \leq C \log n + 3 \clemlu \epsilon n/C_0$, the bounded energy property of $P$ implies that 
$$
 \frac{P \big(
   \bo{u_0} \leftrightarrow \bo{v_0} \big)}{P \big(
   \bo{x_1} \leftrightarrow  \bo{y_1} \big) } \leq
 \cposen^{-\frac{4 \clemlu \epsilon n}{C_0}}.
$$
Setting $C_0$ to be a sufficiently high constant, we obtain, from (\ref{eqconest}) and Lemma \ref{lemmac},
that
$a_{\bo{u}} \leq \exp \big\{ - c \epsilon n \big\}$, as we sought (in (\ref{aubd})). \\
\noindent{\bf The derivation of (\ref{rwrg}) and (\ref{angvouo}).}
Let $\bo{v}$ be an element on the left-hand-side of (\ref{rwrg}). 
We must verify that
\begin{equation}\label{vciinc}
 \cir \cap \Big( W_{\bo{v},c_0}\big( \bo{0} \big) \setminus B_K \big( \bo{v} \big)  \Big)
\subseteq \clus{\bo{v}}{\bo{v}},
\end{equation}
and that
\begin{equation}\label{gamclus}
 \cir \cap B_K \big( \bo{v} \big) \subseteq  \clus{\bo{v}}{\bo{v}}.
\end{equation}
To derive (\ref{gamclus}), note that $\big( P_1 \cup P_2 \big) \cap W_{\bo{u},\epsilon^*/2} = \emptyset$ implies that  $\cir \cap W_{\bo{u},\epsilon^*/2}\big( \bo{0}\big) = \gamma' \cap W_{\bo{u},\epsilon^*/2}\big( \bo{0}\big)$. Using this, as well as  $B_K\big( \bo{v} \big) \subseteq W_{\bo{u},\epsilon^*/2}\big( \bo{0}\big)$ and $\bo{v} \in {\rm RG}^{\bo{u_0},\bo{v_0}}_{\qzero/2,K,\phi} \big( \gamma' \big)$, we find that
 $E(\cir) \cap E^* \big( B_K(\bo{v}) \big) = \bo{v} + \phi$. Hence,
 $$
  \cir \cap B_K \big( \bo{v} \big) \subseteq \Big( W_{\bo{u_0},2\qzero}\big( \bo{v} \big) \cup  W_{- \bo{u_0},2\qzero}\big( \bo{v} \big)  \Big)^c,
 $$
by the choice of $\phi$ and since $\cir \cap B_K(\bo{v})$ is contained in the union of the members of the edge-set $E(\cir) \cap E^* \big( B_K(\bo{v}) \big)$. To derive (\ref{gamclus}), it suffices then to show that
 $$
  \bigg(\Big( W_{\bo{u_0},2 \qzero}\big( \bo{v} \big) \cup  W_{- \bo{u_0},2 \qzero}\big( \bo{v} \big) \Big) \bigg)^c
  \subseteq \clus{\bo{v}}{\bo{v}}.
 $$
 For this, we must show that $\ang \big( \bo{u_0}, \bo{v} \big) \leq \qzero$. Note that $\argu(\bo{u_0})$ differs from $\argu(\bo{u}) + \epsilon^*/2$ by at most $\sin^{-1} \big( 1/(c_1 n) \big)$, because $\bo{u_0}$ is a distance at most one from $\ell_1^*$, and $\vert\vert \bo{u_0} \vert\vert \geq c_1 n$, (since $\bo{u_0} \in \vcir$). From $\bo{v} \in W_{\bo{u},\epsilon/4}\big( \bo{0} \big)$, we find that $\big\vert \argu\big( \bo{v} \big) - \argu\big( \bo{u} \big) \big\vert \leq \epsilon/4$, so that $\big\vert \argu(\bo{v}) - \argu(\bo{u_0}) \big\vert \leq \epsilon/4 + \epsilon^*/2 + \sin^{-1} \big( 1/(c_1 n) \big) \leq \epsilon$, for $n$ high. We see that $\epsilon < \qzero$ completes the derivation of (\ref{gamclus}). 
 
Turning to the derivation of (\ref{vciinc}), note that  the set 
$$
 \cir  \cap \Big( W_{\bo{v},c_0}\big( \bo{0} \big) \setminus B_K \big( \bo{v} \big)  \Big)
$$
is contained in the disjoint union
$$
 \Big( \cir  \cap A_{\argu(\bo{v}) - c_0,\argu(\bo{u}) - \epsilon^*/2} \Big) \, \cup \,
 \Big( \cir \cap A_{\argu(\bo{u}) + \epsilon^*/2,\argu(\bo{v}) + c_0} \Big) \, \cup \,
  \Big(  \gamma  \setminus B_K(\bo{v})  \Big),
$$
because  $\cir  \cap W_{\bo{u},\epsilon^*/2} \big( \bo{0}  \big)
 \subseteq   \gamma$.
We note that
\begin{equation}\label{gabcfcb}
 \gamma  \setminus  B_K \big( \bo{v} \big) \subseteq 
 \clus{\bo{v}}{\bo{v}}
\end{equation}
follows from the condition (which is implied by $\bo{v} \in 
{\rm CRG}_{\qzero/2,K,\phi}^{\bo{u_0},\bo{v_0}} \big( \gamma' \big)$)
\begin{equation}\label{cxyinct}
 \gamma'  \setminus  B_K \big( \bo{v} \big) \subseteq 
 W_{-(\bo{v_0} - \bo{u_0}),\qzero/2} \big( \bo{v} \big) \cup  W_{\bo{v_0} - \bo{u_0},\qzero/2} \big( \bo{v} \big),
\end{equation}
along with $\gamma \subseteq \gamma'$,
$W_{-(\bo{v_0} - \bo{u_0}),\qzero/2} \big( \bo{v} \big) \subseteq \clum{\bo{v}}$
and
$W_{\bo{v_0} - \bo{u_0},\qzero/2} \big( \bo{v} \big) \subseteq \clu{\bo{v}}$.
These last two inclusions are equivalent to
\begin{equation}\label{angcon}
 {\rm ang} \big( \bo{v_0} - \bo{u_0},  \bo{v}^{\perp} \big)
\leq \pi/2 - 3\qzero/2.
\end{equation}
We now derive this, obtaining as we do so the promised (\ref{angvouo}).
 Note that $\bo{u_0},\bo{v_0} \in \vcir$ for a
circuit $\cir$
 realizing $\acon$ and $\globdis \leq \frac{\epsilon n}{C_0}$. 
Recall that  $\tcir = n \partial \wulff$ satisfies $\cir \subseteq \tcir + B_{\globdis}$.
Let $\bo{u},\bo{v} \in \tcir$
attain the infimum on $\tcir$ of the distances to $\bo{u_0}$ and $\bo{v_0}$. Note that $\vert \vert \bo{u} - \bo{u_0} \vert\vert, \vert \vert \bo{v} - \bo{v_0} \vert\vert \leq \epsilon n/C_0$. 
Recall from Definition \ref{defrg} that, for $\bo{e} \in S^1$, $w_{\bo{e}}$ denotes the  clockwise-oriented  unit tangent vector
to $\bo{e}$ on $\partial \wulff$, and that  
$\qzero > 0$ is chosen to satisfy (\ref{supang}). Note that
\begin{equation}\label{supangfe}
 \sup \Big\{   \ang \big( \bo{f} - \bo{e}, w_{\bo{e}} \big): \bo{e},\bo{f} \in \tcir, \ang \big( \bo{e},\bo{f} \big)  \leq r \Big\} \to 0
\end{equation}
as $r \to 0$, since $\tcir$ is a dilation of the compact and differentiable curve $\partial \wulff$ (as known from Lemma \ref{lemozstr}). 

We now bound the angle between $\bo{u}$ and $\bo{v}$. Note that $\ang \big( \bo{v_0},\bo{u_0} \big) \geq \epsilon^*$, because the sector $W_{\bo{u},\epsilon^*/2}$ separates $\bo{u_0}$ and $\bo{v_0}$. Note also that
\begin{equation}\label{vouoineq}
\vert\vert \bo{v_0} - \bo{u_0} \vert\vert \geq 2 \pi^{-1} \min \Big\{ \vert\vert \bo{u_0} \vert\vert, \vert\vert \bo{v_0} \vert\vert  \Big\} \ang \big( \bo{u_0}, \bo{v_0} \big) \geq 2 \pi^{-1} c_1 n \epsilon^*,
\end{equation}
where we used $\bo{u_0},\bo{v_0} \in \cir$ and $\cir \cap B_{\ccone n} = \emptyset$ in the second inequality. 

From $\vert\vert \bo{u} - \bo{u_0} \vert\vert, \vert\vert \bo{v} - \bo{v_0} \vert\vert \leq \epsilon n/C_0$ and $\vert\vert \bo{v_0} - \bo{u_0} \vert\vert \geq \ccone n \epsilon^*$, as well as $\ccone n \epsilon^* \geq 2 \epsilon n/C_0$ (which follows if we set $C_0 \geq 4/c_1$), it readily follows that 
\begin{equation}\label{angpi}
 \ang \big( \bo{v_0} - \bo{u_0}, \bo{v} - \bo{u} \big) \leq \pi \frac{\epsilon n/C_0 + \epsilon n/C_0}{\ccone n \epsilon^*} = \frac{\pi \epsilon}{C_0 \ccone \epsilon^*} \leq \frac{2 \pi}{C_0 \ccone},
\end{equation}
where $\epsilon \leq 2 \epsilon^*$ in the second inequality holds provided that $n$ is high.
Note that 
\begin{equation}\label{angthrang}
 \ang \big( \bo{u}, \bo{v} \big) \leq  \ang \big( \bo{u}, \bo{u_0} \big) + \ang \big( \bo{u_0}, \bo{v_0} \big) +  \ang \big( \bo{v_0}, \bo{v} \big). 
\end{equation}
Now, $\vert\vert \bo{u} - \bo{u_0} \vert\vert \leq \epsilon n/C_0$ and $\vert\vert \bo{u_0} \vert\vert \geq \ccone n$ imply that $\ang \big( \bo{u},\bo{u_0} \big) \leq \frac{\pi \epsilon}{2 C_0 \ccone}$. The same bound holds for $\ang \big( \bo{v},\bo{v_0} \big)$. 
Note that $\bo{u_0}$ and $\bo{v_0}$ lie on opposite sides of the sector $W_{\bo{u},\epsilon^*/2}$, and that each has distance at most one from its boundary. By $\vert\vert \bo{u_0} \vert\vert, \vert\vert \bo{v_0} \vert\vert \geq \ccone n$, we find that 
\begin{equation}\label{anguveps}
 \ang \big( \bo{u_0}, \bo{v_0} \big) \leq \epsilon^* + \pi/2 \cdot 2/(c_1 n) \leq 2 \epsilon.
\end{equation}
Returning to (\ref{angthrang}), we obtain $\ang \big( \bo{u}, \bo{v} \big) \leq 2 \epsilon + \frac{\pi \epsilon}{C_0 \ccone}$. Given any $c' > 0$, we may choose $\epsilon > 0$ small enough that, by use of (\ref{supangfe}), $\ang \big( \bo{v} - \bo{u}, w_{\bo{u}} \big) \leq c'$. Note that 
$\ang \big( \perpu{\bo{v}}, \bo{v_0} - \bo{u_0} \big)
  \leq  \ang \big( \bo{v}, \bo{u_0} \big) + \ang \big( \perpu{\bo{u_0}}, \bo{v_0} - \bo{u_0} \big)$. We have also that
\begin{eqnarray}
&  & \ang \big( \perpu{\bo{u_0}}, \bo{v_0} - \bo{u_0} \big) \nonumber \\
 & \leq &
 \ang \big( \perpu{\bo{u_0}}, \perpu{\bo{u}}  \big) + \ang \big( \perpu{\bo{u}}, w_{\bo{u}} \big) \nonumber \\
 & &   +  \ang \big( w_{\bo{u}},  \bo{v} - \bo{u} \big) + \ang \big( \bo{v} - \bo{u}, \bo{v_0} - \bo{u_0} \big) \leq  \frac{\pi \epsilon}{2 C_0 \ccone}  + \Big( \pi/2 - 3\qzero \Big) + c' + \frac{2\pi}{C_0 \ccone}, \nonumber 
\end{eqnarray}
the second term in the final inequality bounded by (\ref{supang}) and the fourth term by (\ref{angthrang}).
We have that 
$\ang \big( \bo{v}, \bo{u_0} \big) \leq \ang 
\big( \bo{v_0}, \bo{u_0} \big) \leq 2 \epsilon$, the first inequality by $\bo{v} \in W_{\bo{u},\epsilon^*/2}$ and $\bo{u_0}$ and $\bo{v_0}$ being outside the opposite sides of $C_{\bo{u},\epsilon^*/2}$, and the second inequality by (\ref{anguveps}). 
By choosing $c' = \qzero/2$, $\epsilon < \frac{2\qzero}{4 + \pi/(C_0 \ccone)}$ and $C_0 \geq \frac{\pi}{\ccone} \cdot \frac{4}{\qzero}$, we obtain 
$\ang \big( \perpu{\bo{v}}, \bo{v_0} - \bo{u_0} \big) \leq \pi/2 - \qzero$, which is (\ref{angcon}). (We also obtain (\ref{angvouo}), as promised.)
This completes the derivation of (\ref{gabcfcb}).

To verify (\ref{vciinc}), it remains to show that
\begin{equation}\label{wvinc}
 \cir \cap A_{\argu(\bo{v}) - c_0,\argu(\bo{u}) -
   \epsilon^*/2}\big( \bo{0} \big) \subseteq \clum{\bo{v}}
\end{equation}
and that $\cir \cap A_{\argu(\bo{u}) + \epsilon^*/2,\argu(\bo{v}) +
   c_0}\big( \bo{0} \big) \subseteq \clu{\bo{v}}$
We check only (\ref{wvinc}), the other inclusion being similar. Let $\bo{w}$ belong to the left-hand-side of (\ref{wvinc}).
%[See sketch!]
We must verify that
${\rm ang} \big( \bo{w} - \bo{v}, - \bo{v}^{\perp}  \big) \leq
\pi/2 - \qzero$.
As before,
let $\tcir = n \partial \wulff$, so that $\cir \subseteq \tcir + B_{\globdis}$. Let $\bo{\tilde{w}}$ and $\bo{\tilde{v}}$ be points of
$\tcir$ for which 
$d \big( \bo{w} , \tcir \big) = d \big( \bo{w} , \bo{\tilde{w}} \big)$ 
and
$d \big( \bo{v} , \tcir \big) = d \big( \bo{v} , \bo{\tilde{v}}
\big)$.
Recall the requirement (\ref{czercond}) on $c_0 > 0$.
Thus,
\begin{equation}\label{angwv}
{\rm ang} \Big( \bo{\tilde{w}} - \bo{\tilde{v}}  ,
- \bo{\tilde{v}}^{\perp} \Big) \leq \pi/2 - 2\qzero, 
\end{equation}
provided that
$\big\vert  \argu(\bo{\tilde{w}})  -  \argu(\bo{\tilde{v}})  \big\vert \leq
2 c_0$. But this follows from  
$\big\vert  \argu(\bo{w}) -  \argu(\bo{v})  \big\vert \leq
c_0$,
\begin{equation}\label{twoineq}
\big\vert  \argu(\bo{\tilde{w}}) -  \argu(\bo{w})  \big\vert \leq
c_0/2 \qquad
\textrm{and} \qquad
\big\vert  \argu(\bo{v}) -  \argu(\bo{\tilde{v}})  \big\vert \leq
c_0/2.
\end{equation}
In this regard, note that
$$
d \big( \bo{w} , \bo{\tilde{w}} \big) = d \big( \bo{w} , \tcir \big)
 \leq \globdis \leq \frac{\epsilon n}{C_0},
$$
and $\vert\vert  \bo{w} \vert\vert  \geq c_1 n$ (since $\bo{w} \in \cir$). From these, we obtain 
$\ang \big( \bo{\tilde{w}} , \bo{w} \big) \leq \sin^{-1} \big( \epsilon C_0^{-1} c_1^{-1} \big) \leq \frac{\pi \epsilon}{2 C_0 c_1}$, so that the first inequality of (\ref{twoineq}) follows by imposing $\epsilon < C_0 c_1 c_0 \pi^{-1}$. The second arises identically. 

We have that
\begin{equation}\label{dwveq}
d \big( \bo{w}, \bo{v} \big) \geq \frac{\epsilon \ccone}{6 \pi} n.
\end{equation}
Indeed, $\cir \cap B_{\ccone n} = \emptyset$ implies that
$\vert\vert \bo{w} \vert\vert, \vert\vert \bo{v} \vert\vert \geq \ccone n$,
which, alongside $\ang \big( \bo{w}, \bo{v} \big) \geq \epsilon/12$ (a consequence of $\big\vert \argu(\bo{v}) - \argu(\bo{u})  \big\vert \leq \epsilon/4$ and  $\big\vert \argu(\bo{w}) - \argu(\bo{u})  \big\vert \geq \epsilon^*/2 \geq \epsilon/3$)
and  
$d\big(\bo{w},\bo{v}\big) \geq  2 \pi^{-1}  {\rm ang} \big( \bo{w} , \bo{v} \big) 
 \min \big\{ \vert\vert \bo{w} \vert\vert,\vert\vert \bo{v}
 \vert\vert\big\}$,
yields (\ref{dwveq}). 

From (\ref{dwveq}), $\vert\vert \bo{w} - \bo{\tilde{w}} \vert\vert, \vert\vert \bo{v} - \bo{\tilde{v}} \vert\vert \leq \epsilon n/C_0$,  it follows by imposing that $C_0 \ge 12 \pi c_1^{-1}$, analogously to (\ref{angpi}), that 
\begin{equation}\label{angwvt}
{\rm ang} \Big( \bo{\tilde{w}} - \bo{\tilde{v}}  ,
    \bo{w} - \bo{v}   \Big) \leq \frac{12 \pi^2}{\ccone C_0}.
\end{equation}
Note that 
\begin{equation}\label{angvv}
{\rm ang} \Big( 
 - \bo{v}^{\perp}  ,
- \bo{\tilde{v}}^{\perp}
       \Big) 
 = 
{\rm ang} \big( 
\bo{v}  ,
\bo{\tilde{v}}
       \big) 
\leq \frac{\pi \epsilon}{2 C_0 \ccone},
\end{equation}
due to  $d \big(\bo{v}, \bo{\tilde{v}}
\big) \leq \frac{\epsilon n}{C_0}$
and $\vert\vert \bo{v} \vert\vert \geq \ccone n$.

By (\ref{angwv}), (\ref{angwvt}) and (\ref{angvv}),
$$
{\rm ang} \Big( 
    \bo{w} -    \bo{v}  ,
- \bo{v}^{\perp}
       \Big) 
\leq  \pi/2 - 2\qzero  \frac{12 \pi^2}{\ccone C_0} +   \frac{\pi \epsilon }{2 C_0 \ccone},
$$
which gives ${\rm ang} \big( 
    \bo{w} -    \bo{v}  ,
- \bo{v}^{\perp}
       \big) \leq \pi/2 - \qzero$ if we choose $C_0 > \frac{\pi}{2\qzero \ccone}\big( 24 \pi + \epsilon \big)$. 
Hence, $\bo{w} \in \clum{\bo{v}}$, which confirms (\ref{wvinc}). We have obtained (\ref{vciinc}), and thus (\ref{rwrg}). \qed
\noindent{\bf Proof of Lemma \ref{lemlu}.} The set $\tcir$ is a curve in $\R^2$ which may be parametrized in polar coordinates as a continuous function of the angular variable. This readily implies that the intersection of $\tcir + B_{\epsilon n/C_0}\big( \bo{0} \big)$ with any semi-infinite line segment emanating from the origin is a line segment. Let $V$ denote this planar line segment in the case of $\ell_1^*$, i.e., $V = \ell_1^* \cap \big( \tcir  + B_{\epsilon n/C_0}(\bo{0}) \big)$. We will show that $\vert V \vert \leq C \epsilon n/C_0$, for a large constant $C$. We may view $V$ as being formed continuously by moving a ball $B$ of radius $\epsilon n/C_0$ along the curve $\tcir$, from a location where the ball does not intersect $\ell_1^*$, until it has entirely past through $\ell_1^*$. We may then grow $V$ by including the current intersection of $B$ with $\ell_1^*$ to an accumulating line interval, which ends being equal to $V$ when the ball finishes passing through $\ell_1^*$. As a convenience, we orient $\R^2$ so that $\ell_1^*$ is the positive $y$-axis. We parametrize the motion of $B$ by declaring that first contact of $B$ with $\ell_1^*$ occurs at time $t = 0$, and that the ball moves along $\ell_1^*$ counterclockwise, with its $x$-coordinate moving at unit speed. Thus, $B$ finally stops intersecting $\ell_1^*$ at time $2\epsilon n/C_0$, because this is the diameter of $B$.
  
  We write $\bo{\chi(t)}$ for the centre of the ball at time $t$, and $V_t = B(t) \cap \ell_1^*$. We have then that $V = \cup_{t = 0}^{2\epsilon n/C_0} V_t$. Now, $V_t$ is a line segment of length at most $2\epsilon n/C_0$. It is convenient to consider instead $W_t$, which we define to be the line segment of length $2\epsilon n/C_0$ with the same midpoint as $V_t$. 
  
  Setting $W = \cup_{t=0}^{2\epsilon n/C_0} W_t$, we have $V \subseteq W$, since $V_t \subseteq W_t$ for each $t$. The intervals  $W_t$ have lower and upper endpoints $W_t^l$ and $W_t^u$ whose motion is simple to describe. Let $\theta_t$ denote the angle relative to the positive $y$-direction made by the velocity $\bo{\chi'(t)}$ of the centre of $B$ at time $t$. Then $\big\vert \frac{d}{dt} W_t^l \big\vert = \big\vert \frac{d}{dt} W_t^u \big\vert = \cot \theta_t$. Thus,
  $$
  \Big\vert W_{2\epsilon n/C_0}^l - W_0^l \Big\vert \leq \frac{2\epsilon n}{C_0}
   \sup \Big\{ \Big\vert  \frac{d}{dt} W_t^l \Big\vert: 0 \leq t \leq 2\epsilon n/C_0 \Big\} = 
 \frac{2\epsilon n}{C_0}
   \sup \Big\{ \cot \theta_t  : 0 \leq t \leq 2\epsilon n/C_0 \Big\}. 
  $$
  The same conclusion holds for $\big\vert W_{2\epsilon n/C_0}^u - W_0^u \big\vert$. Since $W \subseteq \big[ \inf_{t \in [0,2\epsilon n/C_0]} W_t^l, \sup_{t \in [0,2\epsilon n/C_0]} W_t^u \big]$, we find that the conclusion $\vert V \vert \leq C \epsilon n/C_0$ will follow if we show that $\sup \big\{ \cot \theta_t: 0 \leq t \leq 2\epsilon n/C_0 \big\}$ is bounded away from infinity, uniformly in $n$. To this end, let $\phi_t = \ang \big(\bo{\chi(t)},\bo{\chi'(t)} \big)$.
In the notation of (\ref{supang}), $\phi_t = \pi/2 - \ang \big( w_{\bo{z}}, \perpu{\bo{z}} \big)$, where $\bo{z} = \bo{\chi(t)}$.
Hence, by (\ref{supang}), $\inf_{t \in \big[ 0,2\epsilon n/C_0 \big]} \phi_t \geq \delta$. Note further that $\big\vert \theta_t - \phi_t  \big\vert$ is equal to the angle $\psi_t$ between $\ell_t^*$ (which is $\pi/2$ in the considered coordinate frame) and $\bo{\chi(t)}$. We will now argue that $\psi_t$ is at most $\delta/2$ whenever $0 \leq t \leq 2\epsilon n/c_0$. This will yield 
\begin{equation}\label{inftheta}
\inf_{t \in [0,2\epsilon n/C_0]} \theta_t \geq \delta/2,
\end{equation}  
so that $\vert V \vert \leq 2 \cot(\delta/2) \epsilon n/C_0$. Recall that $\cir \cap B_{\ccone n} = \emptyset$. This forces $\tcir \cap B_{\ccone n/2} = \emptyset$: for otherwise, $\globdis \geq \ccone n/2$, whereas we have that $\globdis \leq \epsilon n/C_0$. (We impose $\epsilon < \ccone C_0/2$.)

Note that $\bo{\chi(t)}$ has $x$-coordinate $-t$ and, since 
 $\bo{\chi(t)} \in \tcir$, $\vert\vert \bo{\chi(t)} \vert\vert \geq \ccone n/2$. Thus, $\sin \psi_t \leq t/(\ccone n/2)$, whence $\psi_t \leq \frac{2\epsilon n/C_0}{\ccone n/2} = \frac{4\epsilon}{C_0\ccone}$ whenever $0 \leq t \leq 2\epsilon n/C_0$. Hence, $\psi_t \leq \delta/2$ for such $t$ if we insist that $\epsilon < \delta C_0 \ccone/8$. This completes the proof that $\vert V \vert \leq C \epsilon n/C_0$. \qed
\noindent{\bf Proof of Lemma \ref{lemedgep}.} The path $P_1 = \partial_{\bo{u_0},\bo{u_1}}^+ \ell_1^*$ starts at $\bo{u_0}$ and ends at $\bo{u_1}$, and is oriented, in the sense that each of its horizontal, and its vertical, displacements is in the same direction. Hence, $\vert E(P_1) \vert = \vert u_1(1) - u_0(1) \vert + \vert u_1(2) - u_0(2) \vert$. We bound the length of the planar line segment $\big[ \bo{u_0}, \bo{u_1 }\big]$. 
Note that, since $\bo{u_0},\bo{u_1} \in \partial^+ \ell_1^*$, we have that
$d \big( \bo{u_0}, \ell_1^* \big), d \big( \bo{u_1}, \ell_1^* \big) \leq 2$. To bound $\big\vert \bo{u_1} - \bo{u_0} \big\vert$, we wish to show that $\bo{u_0}$ and $\bo{u_1}$ are close not merely to $\ell_1^*$, but to the set 
$V = \ell_1^* \cap \big( \tcir +  B_{\epsilon n/C_0}(\bo{0}) \big)$ defined in the proof of Lemma \ref{lemlu}. To see this, note that $\bo{u_0} \in \cir$ implies that $\bo{u_0} \in \tcir + B_{\epsilon n/C_0} \big( \bo{0} \big)$. Hence, using again the moving ball $B$, with the same coordinate frame as in the earlier proof, note the ball contains $\bo{u_0}$ at some time.
If we mark this point of contact on
the ball $B$ at this time, and then evolve time backwards, the marked point will follow a trajectory to a point $\bo{q} \in \ell_1^* \cap \big(  \tcir + B_{\epsilon n/C_0} ( \bo{0} ) \big) = V$, and this trajectory will be completed in time $d \big( \bo{u}, \ell_1^* \big) \leq 2$. In the coordinate frame in question, the vertical displacement of $\bo{q}$ from $\bo{u_0}$ is in absolute value at most
$\int_0^2 \cot \theta_t d t$, which is at most $\cot \big( \delta/2  \big) t \leq 2 \cot \big( \delta/2 \big)$, by (\ref{inftheta}), (since $n$ is high). Thus, $d \big( \bo{u_0}, V \big) \leq d \big( \bo{u_0}, \bo{q} \big) \leq \sqrt{4 + 4 \cot^2 \big( \delta/2 \big)}$. The same is true for $d \big( \bo{u_1}, V \big)$.
We obtain
$$
   \big\vert  \bo{u_1} - \bo{u_0} \big\vert \leq d \big( \bo{u_1}, V \big) + \big\vert V \big\vert + d \big( \bo{u_1}, V \big) \leq C \big( 1 + \epsilon n/C_0 \big).
$$
Recalling that 
$\vert E(P_1) \vert = \big\vert  u_1(1) - u_0(1) \big\vert +  \big\vert  u_1(2) - u_0(2) \big\vert$, we obtain 
$\big\vert E(P_1)  \big\vert \leq 2C \big( 1 + \epsilon n /C_0 \big)$. The inequality holds for $\big\vert E(P_2) \big\vert$, yielding the result. \qed
\end{section}
\begin{section}{The area-gain mechanism}\label{secagm}
In the upcoming proof of Theorem \ref{thmmaxrg}, we will work with an application of the sector-storage replacement operation in which the circuit located in the full-plane output traps an area that is only known to exceed a function that is slightly smaller than $n^2$. We will need to know that such a circuit is not much more probable than one that does trap an area of $n^2$. Hence, in this section, 
we will prove: 
\begin{prop}\label{propag}
Let $P = \P_{\beta,q}$ with $\beta < \hat{\beta}_c$ and $q \geq 1$. 
There exists $c > 0$ and $C > 0$, such that, for
$C  \log n  \leq t \leq cn$, we have that
$$
P \Big( \exc \geq n t  \Big\vert \acon \Big) \geq \exp \big\{- ct \big\},
$$
where, on the event $\acon$, 
we define the area-excess $\exc$ to be the quantity $\big\vert \intg \big\vert - n^2$.
\end{prop}
\noindent{\bf Sketch of the proof.} Consider a random operation that  maps $\zoz$ to $\zoz$. The input of the operation is taken to be subcritical percolation, for the purposes of this sketch.  The plane is cut in two along the $y$-axis, the left-hand portion is kept unchanged, and the right-hand portion is shifted by $t$ units to the right, where $t$ is the parameter appearing in the statement of the proposition. In between the two portions is an infinite vertical strip of width $t$. The output configuration in the plane is completed by filling in the strip with an independent percolation (of the same parameter). Clearly, the operation maps the percolation measure to itself. Adopting a similar notation as in the proof of Proposition \ref{proprgmac}, we make use of the terms satisfactory input and successful action. Define input to be satisfactory if $\acon$ occurs. If $\acon$ does occur, then  the circuit $\cir$ is highly likely to   
have a diameter of order $n$, with the origin being close to its centre. We locate points above and below the origin on the $y$-axis, at which the circuit is cutting through. After the shift, each point has two images, lying on opposites sides of the strip, with a horizontal line connecting the two images. Action is successful if, in the new percolation in the strip, these two horizontal paths are both open. If the input is satisfactory (which has probability $P(\acon)$) and action is successful (probability $\exp \big\{ - ct \big\}$), then the output contains an open  circuit that traps area at least $n^2 + ctn$, because it traps an area equal to the area trapped in the domain circuit in the two portions, and an additional $ctn$ in the strip. Since the operation maps percolation to itself, this appears to complete the sketch. Of course, similarly to the discussion in Section 1.1 of \cite{hammondone}, the argument is flawed. The $y$-axis cut must not run through the circuit several times on either the north-side or the south-side, because, otherwise, the attempt to find a circuit in the output may fail. This fault is corrected by selecting points at which regeneration sites are present. We are relying on the limited information so far gathered regarding the presence of such sites: Proposition \ref{proprgmac} tells us that they may be found in any small angle cone. This is enough: we will find one close to the northerly direction, and another close to the southerly direction, and then use a variant of the shifting operation, which allows for the two points, on the north side and the south, possibly having differing $x$-coordinates. 

The shifting operation is now defined formally.

\begin{definition}\label{defsro}
Let $P$ be a given measure on configurations $\zoz$. 
Let $A,B \subseteq \R^2$ and $\bo{x} \in \Z^2$ be such that $E(A) \cap E(B) = \emptyset$
and $E(A) \cap \big( E(B) + \bo{x} \big) = \emptyset$. The shift-replacement operation
$\xi =  \xi_{A,B,\bo{x}}$ is the following random map
$\xi: \zoz \to \zoz$.
Let $\omega \in \zoz$ denote the input of $\phi$. 
We define $\xi(\omega) \big\vert_{E(A)} = \omega \big\vert_{E(A)}$.
The input in $B$ is displaced by $\bo{x}$ 
and recorded as $\xi(\omega)$: 
for each $y \in E(B)$, we set $\xi(\omega) ( y  ) =
\omega \big( y - \bo{x} \big)$. The configuration $\xi(\omega)$ is 
completed by assigning 
$\xi(\omega) \big\vert_{E(\Z^2) \setminus \big( E(A) \cup ( E(B) + \bo{x}) \big)}$
to be random, its law being the marginal on
$E(\Z^2) \setminus \big( E(A) \cup ( E(B) + \bo{x}) \big)$ of the
  conditional distribution of $P$
given the already assigned values 
$\xi(\omega) \big\vert_{E(A) \cup \big( E(B) + \bo{x} \big)}$.
\end{definition}
We will take the operation to act regularly, in the same sense as Definition \ref{defregact}:
\begin{definition}\label{defregactn}
Let $P$, $A$, $B$ and $\bo{x}$ be specified as in Definition \ref{defsro}.
The shift-replacement operation 
$\xi_{A,B,\bo{x}}$ will be said to act regularly if 
\begin{itemize}
 \item the input configuration has the distribution $P$, 
 \item 
the randomness of the action is chosen such that, 
given  
$\xi(\omega) \big\vert_{E(A) \cup \big( E(B) + \bo{x} \big)}$, or, equivalently,  
$\omega \big\vert_{E(A) \cup \big( E(B) \big)}$, the configuration
$\xi(\omega) \big\vert_{E(\Z^2) \setminus \big( E(A) \cup ( E(B) + \bo{x}) \big)}$
is conditionally independent of the configuration 
$\omega \big\vert_{E(\Z^2) \setminus \big( E(A) \cup E(B) \big)}$.
\end{itemize}
\end{definition}
We now give the proof. \\
\noindent{\bf Proof of Proposition \ref{propag}.}
By Proposition \ref{proprgmac}, 
$$
\P \Big( W_{\bo{e_2},\pi/3} \big( \bo{0} \big) \cap \reg \not= \emptyset, 
     W_{- \bo{e_2},\pi/3} \big( \bo{0} \big) \cap \reg \not= \emptyset
     \Big\vert \acon \Big) \geq 1 - \exp \big\{ - cn \big\},
$$
where $\bo{e_1}$ and $\bo{e_2}$ denote the Cartesian unit vectors.
\begin{figure}\label{figareag}
\begin{center}
\includegraphics[width=0.7\textwidth]{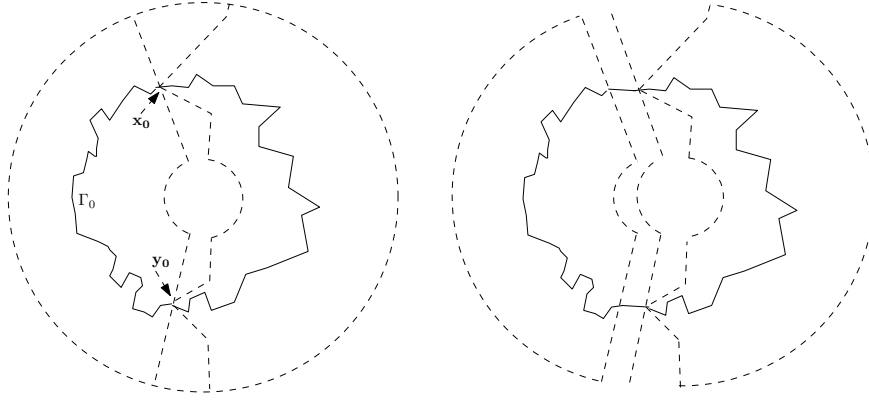} \\
\end{center}
\caption{A pictorial summary of the proof of Proposition \ref{propag}: the circuit $\cir$ in a satisfactory input, and the circuit $\Gamma$ in the output of a successful action of  $\sroapp$.} 
\end{figure}
By this estimate and Lemma \ref{lemmac}, we may find
 $\xo \in \ccup$ and $\yo \in \ccdown$  that satisfy
\begin{equation}\label{xoyointineq}
P \Big(  \big\{ \xo  \big\} \cup  \big\{ \yo  \big\} \subseteq \reg  
     \Big\vert \acon, 
 \cir \cap B_{\ccone n} = \emptyset, 
  \cir \subseteq B_{\cctwo n} \Big) \geq \frac{1}{2\pi^2 \cctwo^4 n^4}.
\end{equation}
Let $H_0 \subseteq \R^2$ consist of those points $\bo{x} \in
\R^2$
such that $\ccone n \leq \vert\vert \bo{x} \vert\vert \leq \cctwo n$ and either
\begin{itemize}
\item $\argu(\yo) + c_0/2 \leq \argu(\bo{x}) \leq \argu(\xo) - c_0/2$, or
\item  $\argu(\xo) - c_0/2 \leq \argu(\bo{x}) \leq \argu(\xo)$ 
   and $\bo{x} \in \clum{\xo}$, or
\item   $\argu(\yo) 
         \leq \argu(\bo{x}) \leq \argu(\yo) + c_0/2$
 and $\bo{x} \in \clu{\yo}$.
\end{itemize}
See Figure 5.
%\ref{figareag}.

We claim that, if $\Gamma$ is a circuit satisfying $\centre(\Gamma) = \bo{0}$, 
 $\Gamma \subseteq B_{\cctwo n} \setminus B_{\ccone n}$ and 
 $\big\{ \xo  \big\} \cup  \big\{ \yo  \big\} \subseteq {\rm RG} \big( \Gamma \big)$, then
\begin{equation}\label{eqnvah}
\Gamma \cap  A_{\yo,\xo} \subseteq H_0.
\end{equation}
Indeed, suppose that $\bo{x} \in \Gamma$.
If  $\argu(\yo)  \leq \argu(\bo{x}) \leq \argu(\yo) + c_0/2$,
then $\bo{x} \in \clu{\yo}$, because $\yo \in {\rm RG}\big( \Gamma \big)$. 
If  $\argu(\xo) - c_0/2  \leq \argu(\bo{x}) \leq \argu(\xo)$,
then $\bo{x} \in \clum{\xo}$, due to $\xo \in \reg$.
From $\ccone n \leq \vert\vert \bo{x}  \vert\vert \leq \cctwo n$ we obtain (\ref{eqnvah}).

Let $Q : = E^* \big( A_{\xo,\yo} \cap B_{\cctwo n} \big)$, where note that the sector 
$A_{\xo,\yo}$ contains the negative $x$-axis. Set $H = E \big( H_0 \big)$.
It is readily seen that there exists $\clemkam \in \N$ and $\clemkac \in \N$ such
that, for all $h \in \Z$, $h \geq 0$ and $n \in \N$,  the sets
$Q$ and $H + (h,0)$ are $(\clemkam,\clemkac)$-well separated, as defined in Lemma \ref{lemkapab}.

Fix $h \in [0,\cctwo n ]$.
We will apply the shift-replacement operation $\sroapp:\zoz \to \zoz$ with a regular action. 
%(Recall from Definition \ref{defsro} that this
%operation constitutes the first component of the storage-shift-replacement
%operation defined there.) 
We will make use of:
\begin{lemma}\label{lemabvar}
There exists  $\ccona = C\big(\clemkam,\clemkac\big)$ such that
the following holds. Let $A,B \subseteq E \big( B_{\cctwo n} \big)$ and $\bo{x} \in \Z^2$
be such that the pairs of edge-sets $(A,B)$ and $(A,B+\bo{x})$ are both $\big(\clemkam,\clemkac\big)$-well
separated. Let $\tilde{P}$ denote the measure on $\zoz$
of the output of an application of $\xi_{A,B,\bo{x}}$ to an input having
the law $P$. Then, for all $\omega \in \zoz$, 
$$
 \ccona^{-1} \leq \frac{d \tilde{P}}{d P} \big( \omega \big)
 \leq \ccona. 
$$ 
\end{lemma}
\noindent{\bf Proof.} 
This follows straightforwardly from Lemma \ref{lemkapab} and the translation invariance of $P$. \qed  
The input of $\sroapp$ is defined to be satisfactory if it realizes the event
$$
{\rm SAT} : =  
  \Big\{    \acon  ,
 \cir \cap B_{\ccone n} = \emptyset,   
\cir \subseteq B_{\cctwo n}, 
\big\{ \xo  \big\} \cup  \big\{ \yo  \big\} \subseteq \reg  
  \Big\}.
$$
Let $P_1$ denote the horizontal path that connects $\xo$ to $\xo + (h,0)$. Let $P_2$
denote the horizontal path that connects $\yo$ to $\yo + (h,0)$.
The action of $\sroapp$ is defined to be successful if 
$E(P_1) \cup E(P_2)$ is open under $\sroapp(\omega)$.
(Note that $E(P_1) \cup E(P_2) \subseteq R$,
where 
$$
R = E(\Z^2) \setminus \bigg( Q \cup \Big( H + \big( h,0 \big) \Big) \bigg),
$$ 
is the region in which $\sroapp(\omega)$ assigns an updated configuration.)

In the case of  successful action on a satisfactory input, we claim that
the output has the property that
\begin{itemize}
\item there exists an open circuit $\Gamma$ such that 
$\Gamma \subseteq B_{2\cctwo n}$
and
$\big\vert {\rm INT} \big( \Gamma \big) \big\vert \geq n^2 + \sqrt{3} hn$.
\end{itemize}
Let $\Gamma$ be such that its edge-set is given by
$$
\Big( E(\cir) \cap E \big( A_{\xo,\yo} \big) \Big) \cup \Big( E(\cir) \cap E \big( A_{\yo,\xo} \big) + (h,0)
\Big) \cup E(P_1) \cup E(P_2)
$$
This edge-set is necessarily open in the output of successful action on satisfactory input.
Note that $\xo,\yo \in \reg$ and $\bo{0} \in \intg$ imply that
$\cir \cap \axyo$ and $\cir \cap A_{\xo,\yo}^c$ are each connected. Thus,
so is $\Gamma$. From $\cir \subseteq B_{\cctwo n}$ and $\vert E(P_1) \vert =
\vert E(P_2) \vert = h \leq \cctwo n$, we obtain  $\Gamma \subseteq B_{2\cctwo n}$.

To bound $\big\vert {\rm INT} \big( \Gamma \big) \big\vert$, note that
\begin{equation}\label{agac}
 \bo{x} \in \intg \cap \axyo \implies \bo{x} \in {\rm INT} \big( \Gamma \big)
 \cap \axyo,
\end{equation}
\begin{equation}\label{agacshift}
 \bo{x} \in \intg \cap \axyo^c \implies \bo{x} + (h,0) \in
 {\rm INT} \big( \Gamma \big) \cap \axyo^c,
\end{equation}
and further that the region 
\begin{equation}\label{agacs}
S = \Big\{  \bo{x} \in \Z^2:  y_0(2) \leq x_2 \leq  x_0(2), \bo{x} \not\in \axyo, \bo{x} - (h,0) \in \axyo \Big\}
\end{equation}
satisfies $S \subseteq \intg$ and is also disjoint from 
$\axyo \cup \big( \axyo^c + (h,0) \big)$. (Here, we use the notation $\bo{v} = (v(1),v(2))$ for $\bo{v} \in \Z^2$.) 
The claims (\ref{agac}), (\ref{agacshift}) and (\ref{agacs}) are straightforward, and similar to (3.24) of \cite{hammondone}.
%similarly to (\ref{intaxyo}). 
We omit the details, and hope that Figure 5
%\ref{figareag} 
is more enlightening.  
By using these three claims, we find that
$$
\big\vert {\rm INT} \big( \Gamma \big) \big\vert 
 \geq  \big\vert \intg \cap \axyo \big\vert +
 \big\vert \intg \cap \axyo^c \big\vert + \big\vert S \big\vert
 = \big\vert \intg \big\vert + \big\vert S \big\vert \geq n^2 + \big\vert S \big\vert. 
$$
Note that
$$
 \big\vert S \big\vert = h \Big(  \big( \xo \big)_2 - \big( \yo
\big)_2  \Big) \geq \sqrt{3} h n. 
$$
Indeed, $\bo{x_0} \in \cir$, $\big\vert \argu(\xo) - \pi/2 \big\vert < \pi/6$ and $\cir \cap B_{\ccone n} = \emptyset$ imply that $\big( \xo \big)_2 \geq \sqrt{3} \ccone n$. Similarly, $\big( \yo \big)_2 \leq - \sqrt{3} \ccone n$.
This establishes the property claimed for the output of $\sroapp$.

In acting on a configuration having the law $P$, the
input of $\sroapp$ is satisfactory with a probability that satisfies
$$
P \big( {\rm SAT} \big) \geq \frac{1}{4 \pi^2 \cctwo^4 n^4} 
 P \Big( \acon \Big)
$$
by (\ref{xoyointineq}) and Lemma \ref{lemmac} in the second inequality. Conditionally on satisfactory input,  
the action of $\sroapp$ is successful with probability at least $c^{\vert
  E(P_1) \vert + \vert
  E(P_2) \vert} = c^{2h}$, by the bounded energy property of $P$. 
Thus,
\begin{eqnarray}
 & & \tilde{P} \Big( \exists \, \textrm{an open circuit $\Gamma: \Gamma \subseteq
   B_{2 \cctwo n}$, $\big\vert {\rm INT} \big( \Gamma \big) \big\vert \geq n^2 + \sqrt{3} h n$} \Big) \nonumber
 \\
 & \geq & \mathbb{P} \Big( \, \textrm{the input is satisfactory and the
   action is successful} \Big) \nonumber \\
 & \geq & P \Big( \acon \Big) \frac{1}{4 \pi^2 \cctwo^4 n^4} 
 c^{2h}, \nonumber
\end{eqnarray}
with $\mathbb{P}$ denoting the law of the regular action of $\sroapp$.
By Lemma \ref{lemabvar}, and Lemma \ref{lemcirprop} with $A = n^2 +  \sqrt{3} h n$,
%\begin{eqnarray}
% & & \frac{1}{8 \pi C^2 n^2}  P  \Big( \exists \, \textrm{an open circuit $\Gamma$: $\Gamma \subseteq
%   B_{2 \cctwo n}$, $\big\vert {\rm INT}\big( \Gamma \big) \big\vert \geq n^2 +  \sqrt{3} h n$} \Big) %\nonumber
% \\
% & \leq &  P \Big( \big\vert \intg \big\vert \geq n^2 +  \sqrt{3} h n \Big), \nonumber
%\end{eqnarray}
we find that
$$
P \Big( \exc \geq  \sqrt{3} h n \Big\vert \acon \Big) \geq \frac{1}{\ccona}
\frac{1}{20 \pi (2\cctwo)^2 n^2} \frac{1}{4 \pi^2 \cctwo^4 n^4}  c^{2h}.
$$
If $h \geq C \log n$ in addition to $h \leq Cn$, we have that
$$
P \Big( \exc \geq  \sqrt{3} h n  \Big\vert \acon \Big) \geq \exp \big\{ - 3\log \big( c^{-1} \big) h \big\},
$$
which yields the statement of the proposition. \qed
\end{section}
\begin{section}{The profusion of regeneration sites in the circuit}\label{secproof}
We are ready to give the proof of Theorem \ref{thmmaxrg}. To explain our approach, consider a configuration $\omega$ satisfying $\acon \cap  \big\{ \mar > u/n \big\}$. 
There exist consecutive elements $\bo{x},\bo{y} \in \reg$, with $\argu(\bo{y}) > \argu(\bo{x}) + u/n$. It is tempting to try to prove the theorem by applying the sector storage-replacement operation $\sigma_{\bo{x},\bo{y}}$ to an input such as $\omega$ (which we would regard as satisfactory input), while defining the successful action of the operation to be the presence of an open path from $\bo{x}$ to $\bo{y}$ with a typical fluctuation in the newly sampled configuruation in $\axy$ in the full-plane output. One would hope that the path $P : = \cir \cap \axy$ from $\bo{x}$ to $\bo{y}$ in the sector output manifests a peculiarity additional to merely connecting these two points, since $\reg \cap \axy = \big\{ \bo{x},\bo{y }\big\}$ in the input. We might aim to show something resembling (\ref{rgcemp}) in the proof of Proposition \ref{proprgmac}, for example, that
 ${\rm CRG}_{c,K,\phi}^{\bo{x},\bo{y}}(P) \cap \axy = \big\{ \bo{x},\bo{y} \big\}$, for some small value of $c$ (and a suitable choice of lattice animal $\phi 
\subseteq B_K$). However, this approach does not work. If $\bo{y}$ happens to lie very close to one of the boundary segments of $\clu{\bo{x}}$, then the path $P$ may be very close to straight (and full of connection regeneration points), but such that the tiny additional fluctuation in the path is enough to ensure $\reg \cap \axy = \big\{ \bo{x}, \bo{y} \big\}$.

A second problem is that the circuit in the full-plane output may capture an area rather less than $n^2$, because the path $P$ sways outwards more than its replacement occurring under successful action. 

A proposal to resolve the first difficulty is to try to locate two regeneration sites $\bo{x'}$ and $\bo{y'}$ on opposite sides of the sector $\axy$ such that $\bo{y'} \in C_{\pi/2 - 2\qzero}^F(\bo{x'})$, namely,  such that $\bo{y'} - \bo{x'}$ is closer to $\perpu{\bo{x'}}$ by a small positive angle (that we happen to choose to be $\qzero$) than is mandated by the fact that $\bo{y'} \in \clu{\bo{x'}}$ (which holds provided that $\ang \big( \bo{x'},\bo{y'} \big) \leq c_0$). It is now realistic to show an analogue of (\ref{rwrg}) 
(see the upcoming (\ref{csetcont})), 
so that the absence of elements of $\reg$ having arguments in the long interval $\big[ \argu(\bo{x}),\argu(\bo{y}) \big]$
corresponds to a long gap in connection regeneration points for $P = \cir \cap A_{\bo{x'},\bo{y'}}$. As such, Lemma \ref{lemmaxreg} may then be invoked to show that this feature is indeed an additional peculiarity for a path in the sector output. 
This means that we may modify the previous argument, now applying 
$\sigma_{\bo{x'},\bo{y'}}$ in place of $\sigma_{\bo{x},\bo{y}}$. Acting on an input satisfying 
 $\big\{ \acon \big\} \cap  \big\{ \mar > u/n \big\}$, the sector output of $\sigma_{\bo{x'},\bo{y'}}$ contains a path 
(that was $\cir \cap A_{\bo{x'},\bo{y'}}$ in the input) from $\bo{x'}$ to $\bo{y'}$ with a peculiarity of conditional probability $\exp \big\{ - cu \big\}$ given the presence of open path from $\bo{x'} \leftrightarrow \bo{y'}$.

In searching for such a pair $(\bo{x'},\bo{y'})$, we must keep in mind the second difficulty, that of area-shortfall in the full-plane output. We want to keep this shortfall as small as possible, by choosing $(\bo{x'},\bo{y'})$ so that $\cir \cap A_{\bo{x'},\bo{y'}}$ does not stray too far outside the triangular region $T_{\bo{0},\bo{x'},\bo{y'}}$. Hence, the following definition:
\begin{definition}\label{defpertp}
Recall $T_{\bo{0},\bo{x},\bo{y}} \subseteq \R^2$ from Definition \ref{defclosedtri}.
Let $\Gamma$
denote a circuit for which $\centre \big( \Gamma \big) = \bo{0}$.
\begin{itemize}
\item
A pair $\big(\bo{u},\bo{v}\big) \in \Z^2 \times \Z^2$,
$\arg \big(\bo{u}\big) < \arg \big(\bo{v}\big)$, $\bo{u}, \bo{v} \in {\rm RG} \big( \Gamma \big)$,
is said to be well-aligned if one of the conditions 
$\bo{u} \in  \clumb{\bo{v}}$
or
$\bo{v} \in  \club{\bo{u}}$
pertains.
\item A pair $(\bo{u},\bo{v}) \in \Z^2 \times \Z^2$,
$\arg\big(\bo{u}\big) < \arg\big(\bo{v}\big)$, $\bo{u}, \bo{v} \in {\rm RG} \big( \Gamma \big)$, is
said to be outward-facing if
$$
{\rm RG} \big( \Gamma \big) \cap \Big( A_{\bo{u},\bo{v}} \setminus T_{\bo{0},\bo{u},\bo{v}} \Big)
= \emptyset.
$$
That is, the regeneration sites of $\Gamma$ having an argument value between that of  $\bo{u}$ and $\bo{v}$
belong to $T_{\bo{0},\bo{u},\bo{v}}$. 
\end{itemize} 
Let $\bo{v},\bo{w} \in \reg$, $\argu(\bo{v}) < \argu(\bo{w})$,
satisfy $A_{\bo{v},\bo{w}} \cap \reg = \big\{ \bo{v},\bo{w} \big\}$,
with $\argu(\bo{w}) - \argu(\bo{v})$ being maximal subject to this
condition. 
(An arbitrary rule should be used to find $\bo{v}$ and $\bo{w}$ if there is
a choice to be made.) Note that $\argu(\bo{w}) - \argu(\bo{v}) = \mar$.
If there exists a pair $\bo{x},\bo{y} \in \Z^2$, 
$\argu(\bo{x}) < \argu(\bo{y})$, of points such
that 
$\big( \bo{x},\bo{y} \big)$
form a well-aligned outward-facing pair, with $A_{\bo{v},\bo{w}} \subseteq \axy$,
let the pair of such points 
for which $\axy$ is minimal be called the {\bf pertinent pair}. 
\end{definition}
\noindent{\bf Remark.} If a pertinent pair $\big( \bo{u}, \bo{v} \big)$, 
 $\argu(\bo{u}) < \argu(\bo{v})$,
satisfies $\ang\big(\bo{u},\bo{v}\big) < \theta$, 
then it is readily seen that $\bo{u} \in C_{\pi/2 - 2\qzero + \theta}^B(\bo{v})$ and $\bo{v} \in C_{\pi/2 - 2\qzero + \theta}^F(\bo{u})$.  

Our new strategy, then, consists of taking $\big(\bo{x'},\bo{y'} \big)$ to be the pertinent pair. We need to know that the angular separation of the elements in the pertinent pair is typically less than a small $n$-independent constant. 
In this regard, we have the following result, whose proof appears after that of 
Theorem \ref{thmmaxrg}.
\begin{lemma}\label{lemexwaofp}
Let $\epsilon > 0$. A circuit $\Gamma$ for which $\centre(\Gamma) = \bo{0}$ is said to be $\epsilon$-full, 
if, for all $\bo{u} \in S^1$, there exists a
well-aligned outward-facing pair $\big( \bo{v}, \bo{w} \big)$ 
for which $\bo{v} \in W^-_{\bo{u},\epsilon/2}$ and
 $\bo{w} \in W^+_{\bo{u},\epsilon/2}$. 
For $\epsilon > 0$ sufficiently small,
\begin{displaymath}
  \Big\{ \cir \subseteq B_{\cctwo n} \setminus B_{c_1 n}  \Big\}
 \cap \Big\{ {\rm GD}\big( \cir \big) \leq
\epsilon n \Big\} \cap \Big\{ \mar \leq \epsilon \Big\} \subseteq  \Big\{ \cir \, \, \textrm{is $4\epsilon$-full} 
  \Big\}.
\end{displaymath}
Specifically, the event on the left-hand-side ensures that the pertinent pair $\big( \bo{x'},\bo{y'} \big)$ exists and satisfies
 $\ang \big( \bo{x'},\bo{y'} \big) \leq 4 \epsilon$.
\end{lemma}
We begin the formal argument.\\
%\begin{prop}\label{propmaxrg}
%Let $P$ satisfy the ratio weak mixing property.
%There exist $c > 0$ and $C > 0$ such that
%$$
%P \Big(  \mar > u/n   \Big\vert \acon \Big)
% \leq  \exp \Big\{ - c u \Big\} 
%$$
%for $C \log n   \leq u \leq c n$.
%\end{prop}
\noindent{\bf Proof of Theorem \ref{thmmaxrg}.} \\
\noindent{\bf Definition of satisfactory input, successful action and of operation parameters:}
Let
$\omega \in \zoz$  denote a configuration. If $\cir = \emptyset$ under $\omega$,
or $\omega$ is such that the pertinent pair does not exist, then $\omega$ is not a satisfactory input. 
If the pertinent pair exists, we record it by $(\bo{x'},\bo{y'})$,
$\arg(\bo{x'}) < \arg(\bo{y'})$.

Let  ${\rm SAT}_1$ be given by
\begin{eqnarray}
  & & \Big\{ \big\vert  \intg \big\vert  \geq n^2 \Big\}
 \cap \Big\{ 2u n^{-1} \geq \mar \geq u n^{-1}   \Big\} \nonumber \\
  & & 
 \cap \Big\{ \textrm{the pertinent pair $(\bo{x'},\bo{y'})$ exists, with $\ang (\bo{x'},\bo{y'}) \leq \epsilon$}  
  \Big\} 
\cap  \Big\{ \cir \subseteq B_{\cctwo n} \setminus B_{\ccone n} \Big\}. \nonumber
\end{eqnarray}
Note that ${\rm SAT}_1  \subseteq \big\{ \cir \subseteq B_{\cctwo n} \big\}$ permits us to find two points $\bo{x_0},\bo{y_0} \in B_{\cctwo n}$ for which
$$
 P \Big(   \bo{x'} = \bo{x_0}, \bo{y'} = \bo{y_0}  \Big\vert {\rm SAT}_1 \Big) \geq \frac{1}{2 \pi^2 \cctwo^4 n^4}.
$$
These points will be the parameters of the sector storage-replacement operation
$\sigma_{\bo{x_0},\bo{y_0}}$ that we will use. Set
 ${\rm SAT}_2 =  \big\{ \bo{x'} = \bo{x_0}, \bo{y'} = \bo{y_0} \big\}$.
The input 
$\omega$ 
will be called satisfactory if it realizes the event 
${\rm SAT} :=  {\rm SAT}_1 \cap {\rm SAT}_2$.

The operation is defined to act successfully if the updated configuration $\omega_1 \big\vert_{\axyoe}$ used by $\sigma_{\bo{x_0},\bo{y_0}}$ realizes the event ${\rm GAC} \big( \xo,\yo \big)$. \\
\noindent{\bf Output properties.}
Fix $\phi \subseteq E^*(B_K)$
that contains a path from $\partial B_K \cap W_{-(\bo{y_0} - \bo{x_0}),\qzero/2}$ to  $\partial B_K \cap W_{\bo{y_0} - \bo{x_0},\qzero/2}$ and satisfies 
$\phi  \cap \big( W_{\bo{x_0},2\qzero} \cup  W_{- \bo{x_0}, 2\qzero}   \big) = \emptyset$. 
That such a $\phi$ exists is ensured by 
\begin{equation}\label{angxoyo}
  \ang \big( \bo{y_0} - \bo{x_0} , \perpu{\bo{x_0}} \big) < \pi/2 - 3\qzero/2.
\end{equation}
This condition is satisfied due to ${\rm SAT} \subseteq \big\{ \ang(\xo,\yo) \leq C_0 \epsilon \big\}$, $C_0 \epsilon < \qzero/2$ and the remark following Definition \ref{defpertp}.
%We defer verifying (\ref{angvouo}) 
%for now, since we will anyway obtain it in due course.

We now claim that, if the input configuration $\omega$ is satisfactory
and the operation $\sigma_{\xo,\yo}$ acts successfully, then the two parts
of the output
$\sigma_{\bo{x_0},\bo{y_0}}(\omega) = \big( \omega_1,\omega_2 \big)$ have the properties that:
\begin{enumerate}
\item {\bf Full-plane output property:} the full-plane configuration $\omega_1 \in \zoz$ has an open circuit $\Gamma$ for which 
$\Gamma \subseteq B_{5 \cctwo n}$,  
%cgac
and
 $$
  \Big\vert  {\rm INT}  \big( \Gamma \big) \Big\vert \geq n^2 - 
 \cfp \big( \qzero/2 \big) \epsilon n u, \qquad \textrm{with $\cfp =  2 \cctwo \big( \cctwo C_0 + 4 \cctwo C_0 \big) \csc^2 \big( \qzero/2 \big)$,}
 $$
\item {\bf Sector open-set property:} and the sector configuration $\omega_2 \in \big\{ 0,1
  \big\}^{\axyoe}$ 
realizes the event, to be denoted by ${\rm SOSP}'$,
that there exists a connected open set in $\axyo$,
to be called  $C_{\bo{x_0},\bo{y_0}}$, such that 
$\{ \xo \} \cup \{ \yo \} \subseteq C_{\bo{x_0},\bo{y_0}}$,
$C_{\bo{x_0},\bo{y_0}}  \cap W_{\xo,Cn^{-1}\log n}^+
\subseteq \clu{\xo}$, 
$C_{\bo{x_0},\bo{y_0}}  \cap W_{\yo,Cn^{-1}\log n}^-
\subseteq \clu{\yo}$,
and $\maxozrg^{\xo,\yo}_{\qzero/2,K,\phi}  \big( C_{\bo{x_0},\bo{y_0}} \big) \geq \frac{c_1}{2\pi} u$.
\end{enumerate} 
\noindent{\bf Remark}. We comment further on the main ideas of the proof before continuing. The quantity  $\cfp \epsilon n u$ in the full-plane circuit property  is the area-shortfall we mentioned before the proof. This area-shortfall is given in essence by the area of the region trapped inside $\big( \cir \cap A_{\bo{x_0},\bo{y_0}} \big) \setminus T_{\bo{0},\bo{x_0},\bo{y_0}}$, a region of ``length'' $\vert\vert \xo - \yo \vert\vert \leq C_0 \epsilon n$ and a ``width'' of the order of the greatest distance between regeneration points in $A_{\xo,\yo}$, namely, $\Theta(u)$. (See (\ref{eqnatcu}) for the precise statement.) This $\Theta(\epsilon n u)$ shortfall has to be purchased back by a use of Proposition \ref{propag}, at a probabilistic cost of $\exp \big\{ - \Theta(\epsilon u) \big\}$, which cost must be paid for by the budget provided by the additional peculiarity  
 $\maxozrg^{\xo,\yo}_{\qzero/2,K,\phi}  \big( C_{\bo{x_0},\bo{y_0}} \big) \geq \frac{c_1}{2\pi} u$ of the path $C_{\bo{x_0},\bo{y_0}}$ in the sector open-set property. This budget is of order $\exp \big\{ - \Theta(u) \big\}$ by Lemma \ref{lemmaxreg}.
That is, due to the pertinent pair having small angular displacement of at most $C_0 \epsilon$, the peculiarity of the sector open-set property is much more unusual than the possibility of removing the area-shortfall appearing in full-plane circuit property. Another way of expressing this point is as follows. There are two competing mechanisms for gaining area by surgery: one is the the area-gain mechanism introduced in Proposition \ref{propag}, and the second is the outward fluctuation of the circuit from the line segment interpolating the pertinent pair. The first of these permits the purchase of an area-gain of
$\Theta(un)$ at the expense of the occurrence of an event of probability $\exp \big\{ - \Theta(u) \big\}$ (the event being the opening of the two 
horizontal paths of length $\Theta(u)$ depicted in Figure 5).
%\ref{figareag}). 
The second permits the purchase of an area-gain of only $\Theta(\epsilon un)$ at the expense of the occurrence of  an event of probability $\exp \big\{ - \Theta(u) \big\}$ (the event being the presence of a  $\cir$-regeneration block of size $\Theta(u)$ in the circuit between $\bo{x'}$ and $\bo{y'}$, since, as we will argue, such a block forces $\maxozrg^{\xo,\yo}_{\qzero/2,K,\phi}  \big( C_{\bo{x_0},\bo{y_0}} \big) \geq \Theta(u)$). That is, the first mechanism is more efficient, so the event entailed by the second does not typically occur. \\
\noindent{\bf Proof of the sector open-set property.} Note that $\omega_2$ is equal to the input configuration $\omega\big\vert_{\axyoe}$ in $\axyoe$. 
We set $C_{\xo,\yo} = \cir \cap \axyo$. Note that
$\big\{ \bo{x'} = \bo{x_0} \big\} \cap \big\{  \bo{y'} = \bo{y_0}\big\}$ implies that $\bo{x_0},\bo{y_0} \in \reg$, which, alongside $\bo{0} \in \intg$, ensures that $C_{\xo,\yo}$ is a connected set, as well as the three inclusions in the definition of ${\rm SOSP}'$.
We restate the final requirement of sector open-set property  as a lemma whose proof appears at the end of the present one:
\begin{lemma}\label{lemxoyobd} 
Satisfactory input and successful action of $\sigma_{\xo,\yo}$ yield a sector output $\omega_2 \in \{ 0,1 \}^{E(\axyo)}$ for which
$$
\maxozrg^{\xo,\yo}_{\qzero/2,K,\phi} \big( C_{\bo{x_0},\bo{y_0}} \big) \geq \frac{c_1}{2\pi} u.
$$
\end{lemma}
\noindent{\bf Proof of the full-plane circuit property.}
%As in the proof of Proposition \ref{propmscbexc}, 
Let $\gamma_{\bo{x_0},\bo{y_0}} \subseteq \axyo$, $\{ \xo \} \cup \{ \yo \} \subseteq \gamma_{\bo{x_0},\bo{y_0}}$
denote the $\omega_2$-open path from $\xo$ to $\yo$ whose existence is ensured by 
$\omega_2 \big\vert_{\axyoe}$
satisfying ${\rm GAC}(\xo,\yo)$. We set 
$\tilde\Gamma = \big( \cir \cap \axyo^c \big) \cup \gamma_{\bo{x_0},\bo{y_0}}$. It is easy (and identical to the corresponding argument in the proof of Proposition 2 of \cite{hammondone}) to obtain $\tilde{\Gamma} \subseteq B_{5 \cctwo n}$. 
%That is,
%$\Gamma$ is the path formed from $\cir$ by replacing $\cir \cap \axyo$
%(which touches $\partial \axyo$ just at $\xo$ and $\yo$) with $\gamma_{\xo,\yo}$
%(which has the same intersection with $\partial \axyo$). 
%
%Let $\cir = \cir(\omega)$ be the circuit defined by virtue of the input $\omega$ satisfying $\acon$.
%Set $\tilde\Gamma = \big( \cir \cap \axyo^c \big) \cup \gamma$, where 
% $\gamma$ denotes the common open cluster of $\bo{x_0}$ and $\bo{y_0}$
%in $\csec$ of the updated sector configuration $\omega_2\big\vert_{\csece}$ in the event of successful %action.
%Let $\Gamma$ denote the outermost open circuit of $\tilde\Gamma$. Note that
%$\big\vert {\rm INT} \big( \Gamma \big) \big\vert = \big\vert {\rm INT} \big(  \tilde\Gamma \big) \big\vert$. 

Note the following properties of $\cir$ and $\Gamma$:
$$
 \Big\vert {\rm INT} (\Gamma) \Big\vert = 
 \Big\vert {\rm INT} (\Gamma) \cap \aarg{x_0}{y_0} \Big\vert
 +   \Big\vert {\rm INT} (\Gamma) \cap \Big( \R^2
 \setminus \aarg{x_0}{y_0} \Big) \Big\vert,
$$
$$
 {\rm INT} (\Gamma) \cap \Big( \R^2
 \setminus \aarg{x_0}{y_0} \Big) = 
 {\rm INT} \big(\cir\big) \cap \Big( \R^2
 \setminus \aarg{x_0}{y_0} \Big),
$$
\begin{equation}\label{intloghalf}
 \Big\vert {\rm INT} \big( \Gamma \big) \cap \aarg{x_0}{y_0} \Big\vert
 \geq   \big\vert T_{\bo{0},\bo{x_0},\bo{y_0}} \big\vert  +
 \frac{\vert\vert \bo{x_0} - \bo{y_0} \vert\vert^{3/2}}{10} \big( \log
 \vert\vert  \xo - \yo \vert\vert \big)^{1/2},
\end{equation}
\begin{equation}\label{eqnatcu}
 \Big\vert {\rm INT} (\cir) \cap \aarg{x_0}{y_0} \Big\vert
 \leq   \big\vert T_{\bo{0},\bo{x_0},\bo{y_0}} \big\vert  + 
 \Big( \vert\vert \bo{x_0} - \bo{y_0} \vert\vert + 
 4 \cctwo \csc \big( \qzero/2  \big) u \Big)  2 \cctwo \csc \big( \qzero/2 \big) u.
\end{equation}
Only (\ref{eqnatcu}) requires justification beyond that provided in the proof of Proposition 2 of \cite{hammondone}.
To this end, $h > 0$, let $\ell_h$ denote the planar line
parallel to $\ell_{\xo,\yo}$ at distance $h$ from $\ell_{\xo,\yo}$ and in the
opposite connected component of $\R^2 \setminus \ell_{\xo,\yo}$ to that to
which $\bo{0}$ belongs. 
Let $\bo{x_0^*}$ and $\bo{y_0^*}$ denote the intersections of $\ell_{2 \cctwo \csc \big( \qzero/2 \big) u}$ with
$\ell_{\bo{0},\xo}$ and $\ell_{\bo{0},\yo}$. We claim that 
$\cir \cap \axyo \subseteq T_{\bo{0},\bo{x_0^*},\bo{y_0^*}}$, in other words, that,
$\cir \cap \axyo$ lies in the bounded component of $\axyo \setminus
\ell_{2 C \csc \big( \qzero/2 \big) u}$. Indeed, let $\bo{v} \in \cir \cap \axyo$ satisfy that $\bo{v}$
lies in the unbounded component of $\axyo \setminus \ell_{\xo,\yo}$. We must
show that $d \big( \bo{v}, \ell_{\xo,\yo} \big) \leq 2 \cctwo \csc(\qzero/2) u$.
Consider the sector $W_{\bo{v},2un^{-1}}^+$. Since $\mar \leq 2 u
n^{-1}$, $\reg \cap W_{\bo{v},2un^{-1}}^+ \not= \emptyset$. Let $\bo{x_1}$
denote the element of this last set of minimal argument. Then
$\argu(\bo{v}) \leq \argu(\bo{x_1}) \leq\argu(\bo{y_0})$, because $\yo \in \reg$. Note that
$\bo{x_1} \in T_{\bo{0},\xo,\yo}$, because $\bo{x_1} \in \reg \cap
A_{\xo,\yo}$, and the pair $(\xo,\yo)$ is outward-facing. Therefore, 
$d \big( \bo{v}, \ell_{\xo,\yo} \big) \leq d \big( \bo{v}, \bo{x_1}
\big)$. 
Take $\bo{x} =\bo{x_1}$ and $\bo{y} = \bo{v}$ in Lemma \ref{lemdistang}. The hypotheses are satisfied, because 
$\ang \big( \bo{v},\bo{x_1 }\big) \leq 2 u n^{-1} \leq c_0$, so that $\bo{x_1} \in \reg$ gives  
$\ang \big( \bo{v} - \bo{x_1} , - \perpu{\bo{x_1}} \big) \leq \pi/2 - \qzero$.
From the lemma, we find that  $d\big(\bo{x_1},\bo{v}\big) \leq d\big(\bo{x_1},\bo{q}\big) \csc \big( \qzero/2 \big)$. We have that
$d\big(\bo{x_1},\bo{q}\big) = \vert\vert \bo{x_1} \vert\vert \sin \ang \big( \bo{x_1} , \bo{v} \big) \leq 2 \cctwo  u$, since $\bo{x_1} \in \cir \subseteq B_{\cctwo n}$. This shows that $d \big( \bo{x_1}, \bo{v} \big) \leq 2 \cctwo \csc(\qzero/2) u$, as required for  
$\cir \cap \axyo \subseteq T_{\bo{0},\bo{x_0^*},\bo{y_0^*}}$. 
 
From this inclusion, it clearly follows that 
\begin{equation}\label{intato}
\intg \cap \axyo \subseteq  T_{\bo{0},\bo{x_0^*},\bo{y_0^*}}.
\end{equation} 
Note that 
\begin{equation}\label{ttxyineq}
 \big\vert T_{\bo{0},\bo{x_0^*},\bo{y_0^*}} \big\vert \leq 
 \big\vert T_{\bo{0},\xo,\yo} \big\vert + \Big( \vert\vert \xo - \yo
 \vert\vert + \vert\vert \bo{x_0^*} - \xo  \vert\vert + \vert\vert \bo{y_0^*} - \yo  \vert\vert \Big) 2 \cctwo \csc \big( \qzero/2 \big) u,
\end{equation} 
because $T_{\bo{0},\bo{x_0^*},\bo{y_0^*}} \setminus T_{\bo{0},\bo{x_0},\bo{y_0}}$ is contained in a rectangle of width $2 \cctwo \csc \big( \qzero/2 \big) u$ and length bounded above by the quantity in brackets in (\ref{ttxyineq}).
Note that 
$\vert\vert \bo{y_0^*} - \yo  \vert\vert 
= 2 \cctwo \csc \big(\qzero/2\big) u \csc \big( \psi
\big)$, where $\psi = \ang \big( \yo^* - \yo, \yo^* - \xo^* \big) = \ang
\big( \yo, \yo - \xo \big)$. Note that one of $\yo \in \club{\xo}$ and 
$\xo \in \clumb{\yo}$ applies, because the pair $(\xo,\yo)$ is
well-aligned. From $\ang \big( \xo,\yo \big) \leq C_0 \epsilon \leq \qzero$, we learn that $\yo \in \clu{\xo}$. Hence, 
 \begin{eqnarray}
 & & \psi = \ang \big( \yo,\yo - \xo \big) \geq \ang \big( \yo, \xo^{\perp} \big) - \ang \big( \xo^{\perp}, \yo - \xo \big) \nonumber \\
  & = & \pi/2 -  \ang \big( \yo, \xo \big) - \ang \big( \xo^{\perp}, \yo - \xo \big) \geq \qzero/2, \nonumber
 \end{eqnarray}
 where, in the last inequality, we used $\ang \big( \yo,\xo \big) \leq C_0 \epsilon \leq \qzero/2$ and
   $\ang \big( \xo^{\perp}, \yo - \xo \big) \leq \pi/2 - \qzero$ (which is equivalent to $\yo \in \clu{\xo}$). 
Hence, $\vert\vert \bo{y_0^*} - \yo  \vert\vert \leq 2 \cctwo \csc^2 \big( \qzero/2 \big)u $.  
Returning to (\ref{ttxyineq}) with this bound and the same one on  $\vert\vert \bo{x_0^*} - \xo  \vert\vert$, we obtain
$$
 \big\vert T_{\bo{0},\xo^*,\yo^*} \big\vert \leq 
 \big\vert T_{\bo{0},\xo,\yo} \big\vert 
 \Big( \vert\vert \bo{x_0} - \bo{y_0} \vert\vert + 
 4 \cctwo \csc \big( \qzero/2  \big) u \Big)  2 \cctwo \csc \big( \qzero/2 \big) u.
$$
By (\ref{intato}), we obtain (\ref{eqnatcu}).

Assembling, we obtain
\begin{equation}\label{areagbd}
 \big\vert {\rm INT} (\Gamma)  \big\vert
 \geq \big\vert {\rm INT} (\cir)  \big\vert  
- \Big( \vert\vert \bo{x_0} - \bo{y_0} \vert\vert + 
 4 \cctwo \csc \big( \qzero/2  \big) u \Big)  2 \cctwo \csc \big( \qzero/2 \big) u.
\end{equation}
Note that, here, we chose to omit the positive term of order 
$\vert\vert \bo{x_0} - \bo{y_0} \vert\vert^{3/2 + o(1)}$ that arises via (\ref{intloghalf}) from successful action entailing the good area capture event ${\rm GAC} \big( \xo,\yo \big)$. It is true that, if the maximum distance between consecutive elements of $\reg$ (which is $\Theta(u)$) is much smaller than the square-root $\vert\vert \xo - \yo \vert\vert^{1/2}$ of the distance between the elements of the pertinent pair, then this positive term is sufficient to ensure $\big\vert {\rm INT} (\Gamma)  \big\vert \geq n^2$. That is, the problem of area-shortfall does not exist in this case. However, the other case must be handled, and, here, this extra term is of no value.  

Note that (\ref{areagbd}) implies that
$\big\vert {\rm INT} (\Gamma)  \big\vert \geq n^2  - 
 \cfp \epsilon n u$,
due to 
$$
u \leq n \mar \leq n \big( \argu(\yo) - \argu(\xo) \big) \leq C_0 \epsilon n
$$
and to
$\big\vert \intg \big\vert \geq n^2$ and
$\vert\vert \yo - \xo \vert\vert \leq 
  \vert\vert \xo \vert\vert 
 \csc \big(\qzero/2 \big) \ang \big( \xo,\yo \big) \leq \cctwo C_0 \csc \big( \qzero/2 \big) \epsilon n$. (The first inequality is due to Lemma \ref{lemdistang} whose hypotheses are satisfied for $(\xo,\yo)$ by virtue of  $\ang(\xo,\yo) \leq C_0 \epsilon \leq c_0$ and $\xo \in \reg$;  the second inequality is due to $\xo \in \cir \subseteq B_{\cctwo n}$.)
We see that, indeed, the full-plane circuit property holds. \\
\noindent{\bf The upper bound on the probability of the two output properties.} 
Note that, since the input configuration has law $P$, the full-plane configuration $\omega_1 \in \zoz$ in the output also has this law. 
%By the same argument yielding (\ref{fpcpub}) and (\ref{proponebd})[], 
By Lemma \ref{lemcirprop},  the full-plane circuit property is thus satisfied by the output with probability at most
 \begin{eqnarray}
 & & P \Big( \exists \, \textrm{an open circuit} \, \, \Gamma: \Gamma \subseteq B_{\cctwo n}, 
  \big\vert {\rm INT} \big( \Gamma \big) \big\vert \geq n^2 - 
 \cfp \epsilon n u
\Big) \nonumber \\
   & \leq &  \cpi \cctwo^2  n^2    P \Big( 
     \big\vert \intg \big\vert \geq n^2  -  
 \cfp \epsilon n u
  \Big). \nonumber
 \end{eqnarray}
We now use Proposition \ref{propag} to gauge the probabilistic cost of 
recovering the area-shortfall  $\cfp \epsilon n u$ in the full-plane circuit property.
The proposition implies that
$$
P \Big( \acon \Big) \geq 
 \exp \Big\{ - c \cfp \epsilon u  \Big\} P \Big(   \big\vert {\rm INT} \big( \Gamma_0 \big) \big\vert \geq 
 n^2 -  \cfp \epsilon u n \Big),
$$   
provided that $C \log n \leq u \leq cn$, which is
our assumption. 
%(Note that the values of $C$ and $c$ may differ from those
%that appear in Proposition \ref{propag}.) 
Thus, the full-plane circuit property is satisfied by the output with
probability at most 
$$
\cpi \cctwo^2 n^2  \exp \Big\{ c \cfp \epsilon u  \Big\} 
P \Big( \acon \Big). 
$$
The action of $\sigma_{\xo,\yo}$ being regular,
the conditional probability that sector open-set property property  is satisfied, given that the full-plane circuit property occurs, is, as in (3.45)
%ptomf 
of \cite{hammondone}, at most
$$
 \sup  \Big\{ P_{\tilde\omega} \big( H \big):
     \tilde\omega \in \{ 0,1 \}^{E(\Z^2) \setminus \axyoe} \Big\}
$$  
which, by Lemma \ref{lemkapab}, is at most
\begin{eqnarray}
& & \cfp  
 P \Big( \bo{x_0} \leftrightarrow  \bo{y_0}, \maxozrg^{\qzero/2,K}_\phi \big( C_{\xo,\yo} \big)
 \geq \frac{c_1}{4\pi}u \Big) \nonumber \\
& \leq &   \cfp  
 P \Big( \bo{x_0} \leftrightarrow  \bo{y_0} \Big) \exp \Big\{ - c 
 \frac{c_1}{2\pi} u \Big\},
\end{eqnarray}
the inequality by Lemma \ref{lemmaxreg} and $u \geq C \log \vert\vert \xo
\leftrightarrow \yo \vert\vert$.

In summary, in acting on an input with law $P$, $\sigma_{\bo{x_0},\bo{y_0}}$ will return an output having the full-plane circuit and sector open-set properties with probability at most 
\begin{equation}\label{prsigubdt}
 \cpi \cctwo^2 n^2   \exp \Big\{ c  \cfp \epsilon u  \Big\} 
P \Big( \acon \Big)
  \cfp P \Big( \bo{x_0} \leftrightarrow  \bo{y_0} \Big)  \exp \Big\{ - 
 \frac{c c_1}{2\pi } u \Big\}.
 \end{equation}
\noindent{\bf The lower bound on the probability of satisfactory input and successful action.}
This derivation coincides with the corresponding one for Proposition 2 of \cite{hammondone}. 
We must derive that $\xo$ and $\yo$ satisfy the hypotheses of Lemma \ref{lemgac}. Note that, under ${\rm SAT}_1$, 
$(\xo,\yo)$ satisfies 
$\vert\vert \xo - \yo \vert\vert \leq 2 \cctwo n$. Moreover, $(\xo,\yo)$ is the pertinent pair and 
and $\ang(\xo,\yo) \leq \epsilon \leq \qzero$. 
Hence, the remark following Definition \ref{defpertp}
implies that $\yo \in \clu{\xo}$ and $\xo \in \clum{\yo}$. From $\ang(\xo,\yo) \geq \mar \geq u n^{-1} \geq C n^{-1} \log n$ and $\xo,\yo \not\in B_{\ccone n}$, we obtain $\vert\vert \xo - \yo \vert\vert \geq 2 \pi^{-1} \ccone C n^{-1} \log n$. Hence, the hypotheses of Lemma \ref{lemgac} are indeed satisfied, by fixing the constants used in the present proposition high enough relative to $\clemgac$.
We learn that the input is satisfactory
and the operation acts successfully with probability at least
\begin{eqnarray}
& & P \Big( {\rm SAT}_1 \Big) \times \frac{1}{2 \pi^4 \cctwo^4 n^4} \times \inf_{\tilde\omega \in \{ 0,1 \}^{E(\Z^2) \setminus \axyoe}}
 P_{\tilde\omega} \Big(  {\rm GAC} \big( \xo, \yo \big) \Big) \nonumber \\
& \geq & 
 P \Big( {\rm SAT}_1 \Big) \times \frac{1}{2 \pi^4 \cctwo^4 n^4} n^{-C}
  P \Big( \bo{x_0} \build\leftrightarrow_{}^{\axyo}  \bo{y_0} \Big). \label{satsucbd}
 \end{eqnarray}
\noindent{\bf Conclusion by comparison of the obtained bounds.}
Since successful action on satisfactory input forces the output is have
the two properties,
the quantity in (\ref{satsucbd}) is at most that in (\ref{prsigubdt}). That is,
 $P \big( {\rm SAT}_1 \big)$ is at most
\begin{eqnarray}
 & & 
%4 \pi^6 C^6 n^6 \cfp n^C 
 n^C \exp \Big\{ c  \cfp \epsilon u  \Big\} 
  \frac{P \big( \xo \leftrightarrow \yo \big)}{P \big( \xo \build\leftrightarrow_{}^{\axyo} \yo \big)} 
 P \Big(  \big\vert \intg \big\vert \geq   n^2 \Big)  \exp \Big\{ - c \cdot
 \frac{c_1}{2\pi} u \Big\} \nonumber \\
 & \leq & 
%4 \pi^6 C^6 n^6 \cfp n^C  
 n^C \exp \Big\{ c  \cfp \epsilon u  \Big\} 
 P \Big(  \big\vert \intg \big\vert \geq   n^2  \Big)  \exp \Big\{ - c \cdot
 \frac{c_1}{2\pi} u \Big\}, \label{satozbd}
 \end{eqnarray}
the inequality due to the bound 
\begin{equation}\label{xyarbd}
 \frac{P\big( \xo \build\leftrightarrow_{}^{\axyo} \yo \big)}{P \big( \xo
   \leftrightarrow \yo \big)} \geq c,
\end{equation}
which follows straighforwardly from Lemma \ref{lemoznor} with $\delta$ any positive value
less than $\qzero/2$, (since, as we have seen, 
$\yo \in \clu{\xo}$ and $\xo \in \clum{\yo}$).
By choosing $\epsilon > 0$
fixed and small, and noting that $u \geq C \log n$, we obtain
\begin{equation}\label{psatbdrg}
P \Big( {\rm SAT}_1 \Big) \leq \exp \Big\{ - \frac{c c_1 u}{4 \pi} \Big\} P \Big( \acon \Big).
\end{equation}
That is, the second exponential term in (\ref{satozbd}) is more significant than the other term, if $\epsilon > 0$ is small. This reflects the efficacy of the area-gain mechanism of Proposition \ref{propag} over that of local bulging at the circuit boundary due to large regeneration clusters which was discussed in remark following the statement of the two output properties.
By Lemma \ref{lemexwaofp}, Proposition \ref{propglobdis}, Lemma \ref{lemmac} and (\ref{psatbdrg}), we obtain
\begin{eqnarray}
& & P \Big(  2u n^{-1} \geq \mar \geq u n^{-1} \Big\vert \acon \Big) \nonumber \\
& \leq &  \exp \Big\{ - \frac{c c_1 u}{4 \pi} \Big\}  +  \exp \big\{
 - c(\epsilon)  n \big\} \leq  \exp \Big\{ - \frac{c c_1 u}{8 \pi} \Big\} ,   \nonumber
\end{eqnarray}
the latter inequality due to $u \leq cn$.
By replacing $u$ by $2^i u$ for those $i \in \N$ such that $2^i u \leq cn$,
and summing, we obtain
$$
P \Big(  c \geq \mar \geq u n^{-1} \Big\vert \acon \Big)
 \leq   \exp \Big\{ - \frac{c c_1 u}{16 \pi} \Big\}.  
$$
Alongside Proposition \ref{proprgmac}, this completes the proof, subject to the following: \\
\noindent{\bf Proof of Lemma \ref {lemxoyobd}.}
We begin by showing that
\begin{equation}\label{csetcont}
{\rm CRG}_{\qzero/2,K,\phi}^{\xo,\yo} \big( C_{\bo{x_0},\bo{y_0}} \big) \setminus \Big( B_K(\bo{x_0}) \cup  B_K(\bo{y_0})  \Big)
 \subseteq \reg.
\end{equation}
Let $\bo{v}$ be an element on the left-hand-side. We must verify that
\begin{equation}\label{vgcont}
  \cir  \cap \Big( W_{\bo{v},c_0}\big( \bo{0} \big) \setminus B_K \big( \bo{v} \big)  \Big)
\subseteq \clus{\bo{v}}{\bo{v}},
\end{equation}
and that
$$
\cir  \cap  B_K \big( \bo{v} \big)   \subseteq \clus{\bo{v}}{\bo{v}}.
$$
Regarding (\ref{vgcont}), note that the set
$\cir  \cap \Big( C_{\bo{v},c_0}\big( \bo{0} \big) \setminus B_K \big( \bo{v} \big)  \Big)$
is contained in the disjoint union
$$
 \Big( \cir  \cap A_{\argu(\bo{v}) - c_0,\argu(\bo{x_0}})
 \Big) \, \cup \,
  \Big(   C_{\bo{x_0},\bo{y_0}}   \setminus  B_K(\bo{v})  \Big)  \, \cup \, 
   \Big(  \cir \cap A_{\argu(\bo{y_0}),\argu(\bo{v}) + c_0} \Big)  
$$
because  $C_{\bo{x_0},\bo{y_0}}   =   \cir  \cap A_{\bo{x_0},\bo{y_0}}$.
We note that
\begin{equation}\label{cxyinc}
  C_{\bo{x_0},\bo{y_0}}   \setminus  B_K \big( \bo{v} \big) \subseteq \clus{\bo{v}}{\bo{v}}
\end{equation}
follows from the condition (implied by $\bo{v} \in 
{\rm CRG}_{\qzero/2,K,\phi} \big( C_{\bo{x_0},\bo{y_0}} \big)$)
\begin{equation}\label{cxyinctnew}
   C_{\bo{x_0},\bo{y_0}}  \setminus  B_K \big( \bo{v} \big) \subseteq 
 W_{-(\bo{y_0} - \bo{x_0}),\qzero/2} \big( \bo{v} \big) \cup  W_{\bo{y_0} - \bo{x_0},\qzero/2} \big( \bo{v} \big),
\end{equation}
along with 
$W_{-(\bo{y_0} - \bo{x_0}),\qzero/2} \big( \bo{v} \big) \subseteq \clum{\bo{v}}$
and
$W_{\bo{y_0} - \bo{x_0},\qzero/2} \big( \bo{v} \big) \subseteq \clu{\bo{v}}$.

We now show the last two inclusions. It suffices to show that
\begin{equation}\label{stscon}
\ang \big(  \bo{y_0} - \bo{x_0} , \bo{v}^{\perp} \big) \leq \pi/2  - \frac{3\qzero}{2}.
\end{equation}
The occurrence of ${\rm SAT}$ implies that $\big( \xo,\yo \big)$ is the pertinent pair, so that either $\yo \in  \club{\xo}$ or $\xo \in \clumb{\yo}$. In the first case, we have that
$\ang \big(  \bo{y_0} - \bo{x_0} , \bo{x_0}^{\perp} \big) \leq \pi/2 - 2\qzero$. Note then that
$$
\ang \big(  \bo{v}^{\perp}, \bo{x_0}^{\perp} \big) = 
\ang \big(  \bo{v},  \bo{x_0} \big) \leq 
\ang \big(  \bo{y_0},  \bo{x_0} \big) \leq C_0 \epsilon \leq \qzero/2,  
$$
the second inequality by means of the occurrence of ${\rm SAT}$. 
The case that $\xo \in \clumb{\yo}$ is similar. This 
yields (\ref{stscon}) and, thus, (\ref{cxyinc}).

For (\ref{vgcont}), it remains to verify that
\begin{equation}\label{vincone}
  \cir  \cap A_{\argu(\bo{v}) - c_0,\argu(\bo{x_0})}
  \subseteq  \clum{\bo{v}} 
\end{equation}
and that
\begin{equation}\label{vinctwo}
  \cir  \cap A_{\argu(\bo{y_0}),\argu(\bo{v}) + c_0}
  \subseteq  \clu{\bo{v}}. 
\end{equation}
Note that, by virtue of $\xo,\yo \in \reg$, 
$\cir  \cap A_{\argu(\xo) - c_0,\argu(\bo{x_0})}
  \subseteq \clum{\xo}$
and
$\cir \cap A_{\argu(\yo),\argu(\bo{y_0}) + c_0}
  \subseteq \clu{\yo}$.
Due to 
$A_{\argu(\bo{v}) - c_0,\argu(\xo)} \subseteq 
 A_{\argu(\xo) - c_0,\argu(\xo)}$
and
$A_{\argu(\yo),\argu(\bo{v}) + c_0} \subseteq A_{\argu(\yo),\argu(\yo) + c_0}$, 
it suffices for (\ref{vincone}) and (\ref{vinctwo}) to show that
$$
 \clum{\xo} \cap 
 A_{\argu(\bo{v}) - c_0,\argu(\xo)} \subseteq \clum{\bo{v}}, 
$$
and that
$$
 \clu{\yo}  \cap 
 A_{\argu(\yo),\argu(\bo{v}) + c_0} \subseteq \clu{\bo{v}}.
$$
We now show these two statements. 

Let $\bo{z_1},\bo{z_2}$ denote the points of intersection of 
$\partial 
 \clum{\xo}$
and $\big\{ \bo{w} \in \R^2: \argu (\bo{w}) = \argu (\bo{v}) - c_0 \big\}$.

It suffices to show that 
\begin{equation}\label{zvve}
{\rm ang} \Big( \bo{z_1} - \bo{v}, - \bo{v}^{\perp} \Big) \leq \pi/2 - \qzero
\end{equation}
and
$$
{\rm ang} \Big( \bo{z_2} - \bo{v}, - \bo{v}^{\perp} \Big) \leq \pi/2 - \qzero.
$$
To this end, 
note that each of $\partial \clum{\xo}$ and $\partial \clum{\bo{v}}$
consists of a union of two semi-infinite line segments. In the case of
$\partial \clum{\xo}$, one of the two line segments attains the closest
approach to $\bo{0}$ among points in $\partial \clum{\xo}$. Let $\ell_1$
denote the planar line that contains this segment. Let $\ell_2$ denote the
corresponding planar line in the case of $\partial \clum{\bo{v}}$.
Note that ${\rm ang}\big(\ell_1,\ell_2\big) = {\rm ang}(\xo,\bo{v})$. 

We parametrize the line $\ell_1$ by a real variable $t$, in such a way that the point $t = 0$
is the unique point $\bo{p} \in \ell_1$ for which the line $\ell_1^\perp$
orthogonal to $\ell_1$ through $\bo{p}$ intersects $\bo{v}$. 
We write $\ell_1(t)$ for the point on $\ell_1$ parametrized by $t$.
We choose the orientation of the line so that $t(\xo) > 0$, and use the arc-length parametrization. 
The line $\ell_2$ may then be described by means of an affine 
function $h: \R \to \R$,
so that the intersection of $\ell_2$ with the displacement of $l_1^\perp$
through $\ell_1(t)$ has a displacement of $h(t)$ from $\ell_1(t)$. (We take $h(0)
> 0$.) 
There is a unique value $t^* > 0$ for which $h(t^*) = 0$. The lines $\ell_1$ and $\ell_2$ intersect at
$\ell_1 (t^*)$.
\begin{figure}\label{figmanyrg}
\begin{center}
\includegraphics[width=0.4\textwidth]{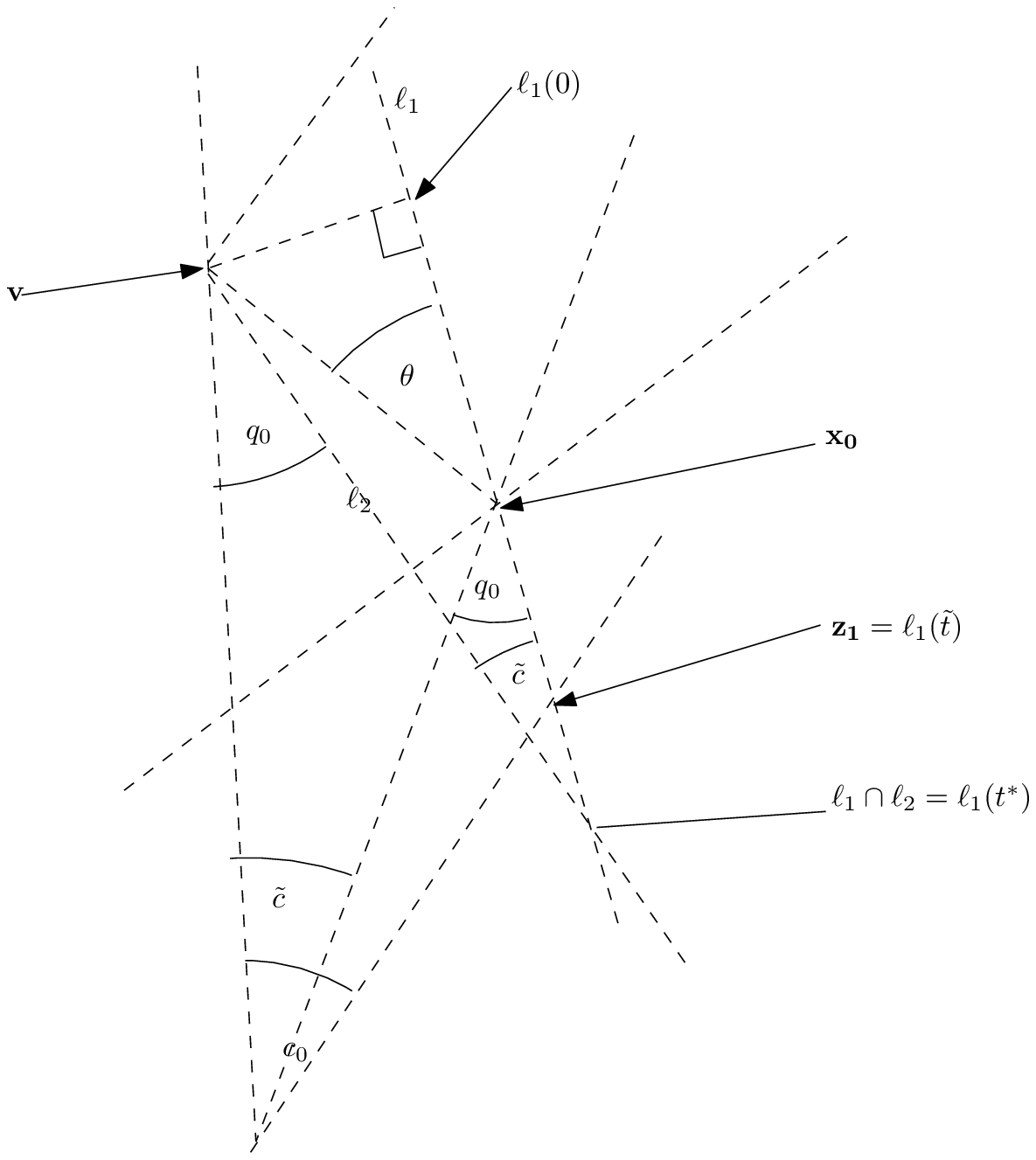} \\
\end{center}
\caption{Deriving (\ref{zvve})}
\end{figure}
The condition 
${\rm ang} \big( \bo{z_1} - \bo{v}, - \bo{v}^{\perp} \big) \leq \pi/2 - \qzero$ is equivalent to $h(\tilde{t}) \geq 0$ for that $\tilde{t}$ for which $\ell_1(\tilde{t}) = z_1$, which is simply the condition that $\tilde{t} \leq t^*$. To obtain (\ref{zvve}), then, we must verify that $\tilde{t} \leq t^*$.

See Figure 6. 
Let $\theta = \ang \big( \bo{v} - \xo, \ell_1(\bo{0}) - \xo \big)$.
We will use the bound $\theta \geq \qzero/4$ to argue that $t^*$ is bounded below by a constant multiple of $n$, whatever the value of $\bo{v} - \xo$, or, specifically,
\begin{equation}\label{eqtstar}
 t^* \geq \frac{c_1}{4} \sin \big( \qzero/4 \big) n.
\end{equation}
(In essence, $\bo{v} \in {\rm CRG}_{\qzero/2,K,\phi}^{\xo,\yo}(C_{\xo,\yo})$ implies that $\bo{v} - \xo$ points in a direction close to $\yo-\xo$. The pair $(\xo,\yo)$ being well-aligned, this forces the angle $\theta$ depicted to satisfy a lower bound $\theta \geq q_0/4$.) 

To obtain (\ref{eqtstar}), note that $h'(t) = \tan(\tilde{c})$, where $\tilde{c} = \ang \big(\ell_1,\ell_2 \big)$ is also given by $\tilde{c} = \ang \big( \bo{v},\xo \big)$. Now, $h(0) =  d(\bo{v},\ell_1(0)) = \vert\vert \bo{v} - \xo \vert\vert \sin (\theta) \geq  \vert\vert \bo{v} - \xo \vert\vert \sin \big( \qzero/4 \big)$, by $\theta \geq \qzero/4$. The function $h$ being affine,
$$
h(t) = h(0) - th'(0) \geq   \vert\vert \bo{v} - \xo \vert\vert \sin \big( \qzero/4 \big) - t \tan (\tilde{c}) 
 \geq  \vert\vert \bo{v} - \xo \vert\vert \frac{\qzero}{2\pi} - 4\pi^{-1} \tilde{c} t,
$$
the second inequality by $\tilde{c} \leq c_0 < \pi/4$. However, $\vert\vert \xo - \bo{v} \vert\vert \geq d\big( \bo{x_0},\ell_{\bo{0},\bo{v}} \big) = \vert\vert \xo \vert\vert \sin(\tilde{c}) \geq 2 \pi^{-1} c_1 n \tilde{c}$, the equality by considering the right-angled triangle with hypotenuse $\big[ \bo{0},\bo{x_0}\big]$ and a third vertex on $\ell_{\bo{0},\bo{v}}$ attaining the closest approach to $\xo$ on this line. The second inequality follows from $\bo{x_1} \in \cir \subseteq B_{\ccone n}^c$. Hence, $h(t) \geq \pi^{-2} \ccone \qzero \tilde{c} n - 4 \pi^{-1} \tilde{c} t$. This means that $h(t)\geq 0$ if $t \in [0,\frac{c_1 \qzero n}{4 \pi}]$, so that we obtain (\ref{eqtstar}). The key point is that our condition for $h(t) \geq 0$ does not depend on the angle $\tilde{c} = \ang \big( \xo, \bo{v} \big)$: if $\bo{v}$ is very close to $\bo{x_0}$, for example, the distance $h(0)$ that $\ell_2(t)$ must catch up with $\ell_1(t)$ as $t$ increases from $0$ is small, but the rate $h'(t) \approx \ang \big( \xo, \bo{v} \big)$ of this catching up scales to be small as well.

For (\ref{eqtstar}), it remains to prove $\theta \geq \qzero/4$. 
Note that $\xo,\yo \in C_{\xo,\yo}  \setminus B_K (\bo{v})$ and (\ref{cxyinctnew}) imply that
$$
{\rm ang} \Big( \bo{x_0} - \bo{v}, - \big( \bo{y_0} - \bo{x_0} \big) \Big) \leq \qzero/2,
$$
(and also that ${\rm ang} \big( \bo{y_0} - \bo{v}, \bo{y_0} - \bo{x_0} \big) \leq \qzero/2$).
We know that
${\rm ang} \big( \bo{y_0} - \bo{x_0},  \bo{x_0}^{\perp} \big) \leq \pi/2 - 2\qzero + \qzero/4$, due to  the pair $\big(\xo,\yo  \big)$ being well-aligned, $\ang(\xo,\yo) \leq C_0 \epsilon$, the remark following Definition \ref{defpertp}, and $C_0 \epsilon < \qzero/4$. We deduce that
\begin{equation}\label{vminxo}
{\rm ang} \Big( \bo{v} - \bo{x_0},  \bo{x_0}^{\perp} \Big) \leq \pi/2 - 2\qzero + \qzero/2 + \qzero/4 = \pi/2 \ 5\qzero/4.
\end{equation} 
Standing at $\xo$, facing in the direction $\perpu{\xo}$ and then turning clockwise by a right-angle, we face towards $\bo{v}$ after rotating by at most $\pi/2 - 5\qzero/4$ (by (\ref{vminxo})), then towards $\ell_1(\bo{0})$ after further rotating by $\theta$, and then continue to rotate by $\qzero$ (since $\ang\big( \ell_1,\xo \big) = \qzero$). Hence, $\theta \geq \qzero/4$, as we sought to show.

We now show that
\begin{equation}\label{eqttilde}
  \tilde{t} \leq 4 \pi^{-1} \sin \big( \qzero/8 \big) c_0 n. 
\end{equation}
Note that (\ref{eqtstar}), (\ref{eqttilde}) and the assumption that $c_0 < \frac{c_1 \qzero^2}{128 \pi C}$ indeed ensures that $t^* \leq \tilde{t}$ and thus, by the chain of reductions, (\ref{csetcont}).

 Note that $\tilde{t} = \vert\vert \bo{z_1} - \bo{l_1}(0) \vert\vert$.
 Note that $\bo{z_1},\bo{l_1}(0)$ and $\bo{v}$ are the vertices of a
 right-angled triangle whose hypotenuse is $\big[ \bo{v} , \bo{z_1} \big]$. Thus,
 $\vert\vert \bo{z_1} - \bo{l_1}(0)  \vert\vert \leq \vert\vert  \bo{v} - \bo{z_1} \vert\vert$. Now,
 \begin{equation}\label{vzonex}
  \vert\vert  \bo{v} - \bo{z_1} \vert\vert \leq  \vert\vert  \bo{v} - \xo \vert\vert +  \vert\vert  \bo{z_1} - \xo \vert\vert. 
 \end{equation}
Consider the right-angled triangle whose vertices include $\bo{0}$ and
$\xo$, with hypotenuse containing part of the line segment $\big[ \bo{0},\bo{v}
\big]$. Let $\bo{q}$ denote the third vertex of the triangle. We have that 
$\vert\vert \bo{q} - \xo \vert\vert = \vert\vert \xo \vert\vert  \tan
\big( \tilde{c} \big)$. Let $\bo{\tilde{q}}$ denote the point $\ell_{\bo{0},\bo{v}}^+$ that is closest to $\xo$.

From (\ref{vminxo})
and  $\ang\big( \bo{\tilde{q}} - \xo, \xo^{\perp} \big) = \tilde{c}$,
we find that  $\ang\big( \bo{v} - \xo, \bo{\tilde{q}} - \xo \big) \leq \pi/2 - 9\qzero/8$, 
by $\tilde{c} \leq  \ang \big( \xo , \yo \big)  \leq C_0 \epsilon < \qzero/8$. 

We have that
\begin{eqnarray}
 & & \vert\vert \bo{v} - \bo{x_0}  \vert\vert \leq
 \vert\vert \bo{\tilde{q}} - \bo{x_0}  \vert\vert \sec \big( \pi/2 - 9\qzero/8 \big) \label{ineqvone} \\
 &  \leq & \vert\vert \bo{q} - \bo{x_0}  \vert\vert \sec \big(  \pi/2 - 9\qzero/8 \big)
 = \vert\vert \xo \vert\vert \tan(\tilde{c})  \sec \big(  \pi/2 - 9\qzero/8 \big), \nonumber
\end{eqnarray}
where, in the second inequality, we used $\bo{q} \in \ell_{\bo{0},\bo{v}}$ and the
definition of $\bo{\tilde{q}}$.

Now consider the right-angled triangle whose vertices include $\bo{0}$ and
$\xo$, and whose hypotenuse contains the line segment $\big[ \bo{0},\bo{z_1}
\big]$. Let $\bo{q'}$ denote the third vertex of the triangle. We have that 
\begin{equation}\label{adhus}
\vert\vert \bo{q'} - \xo \vert\vert = \vert\vert \xo \vert\vert 
\tan \big(
 c_0 - \tilde{c} \big). 
\end{equation}
Let $\bo{q''}$ denote the point on the semi-infinite line
segment from $\bo{0}$ through $\bo{z_1}$ that is closest to $\xo$.

From $\ang\big( \bo{z_1} - \xo, \bo{q'}  - \xo \big) = \pi/2 - \qzero$
and  $\ang\big( \bo{q''} - \xo, \bo{q'} - \xo \big) = c_0 - \tilde{c}$,
we find that  $\ang\big( \bo{z_1} - \xo, \bo{q''} - \xo \big) \leq \pi/2 - \qzero + c_0
- \tilde{c} \leq \pi/2 - \qzero + c_0 \leq \pi/2 - \qzero/2$, 
since $c_0 \leq \qzero/2$. 

We have that
\begin{eqnarray}
 & & \vert\vert \bo{z_1} - \bo{x_0}  \vert\vert \leq
 \vert\vert \bo{q''} - \bo{x_0}  \vert\vert \sec \big( \pi/2 - \qzero/2 \big)
  \nonumber \\
 & \leq &  \vert\vert \bo{q'} - \bo{x_0}  \vert\vert \sec \big( \pi/2 - \qzero/2 \big)
 = \vert\vert \xo \vert\vert \tan\big(c_0 - \tilde{c}\big)  \sec \big( \pi/2 - \qzero/2 \big), \label{ineqvtwo}
\end{eqnarray}
where, in the second inequality, we used that 
$\bo{q'} \in \ell_{\bo{0},\bo{z_1}}$ and the
definition of $\bo{q''}$. In the third inequality, we used (\ref{adhus}). 

From (\ref{vzonex}), (\ref{ineqvone}) and (\ref{ineqvtwo}), we obtain  
$$
  \vert\vert  \bo{v} - \bo{z_1} \vert\vert \leq \vert\vert \xo \vert\vert \tan (\tilde{c}) \sec (\pi/2 - 9\qzero/8) +
     \vert\vert \xo \vert\vert \tan \big( c_0 - \tilde{c} \big) \sec (\pi/2 - \qzero/2).
$$
Note that $\tan(c) \leq 4 \pi^{-1}$ for $c \in (0,\pi/4)$. By 
$\tilde{c} = \ang \big( \xo, \bo{v} \big) \leq  \ang \big( \xo, \yo \big) \leq c_0 \leq \pi/4$, we have that $\tan \big( c_0 - \tilde{c} \big) \leq 4 \pi^{-1} \big( c_0 - \tilde{c} \big)$ and $\tan (\tilde{c}) \leq 4 \pi^{-1} \tilde{c}$.
%By $\tilde{c} = \ang\big( \xo, \bo{v} \big) \leq \ang \big( \xo,\yo \big)
%\leq \qzero/8$, the last inequality of which is ensured by the occurrence of
%${\rm SAT} \cap \big\{ \bo{x'} = \xo, \bo{y'} = \yo \big\}$, and by $\epsilon
%\leq \qzero/8$, 
%we have that $\tan(\tilde{c}) \leq 2 \tilde{c}$ and $\tan \big( c_0 -
%\tilde{c} \big) \leq 2 \big( c_0 - \tilde{c} \big)$. 
Thus,
$$
 \vert\vert  \bo{v} - \bo{z_1} \vert\vert \leq 4 \pi^{-1} \vert\vert \xo \vert\vert
 \sec \big(\pi/2 - \qzero/2 \big) c_0.
$$
By $\tilde{t} \leq \vert\vert \bo{v} - \bo{z_1} \vert\vert$ and $\vert\vert
\xo \vert\vert \leq \cctwo n$ (a consequence of $\cir \subseteq B_{\cctwo n}$), 
we find that 
$$
\tilde{t} \leq 4 \pi^{-1} \cctwo \sec \big(\pi/2 - \qzero/2  \big) c_0 n,
$$
which is (\ref{eqttilde}). This completes the derivation of (\ref{zvve}) and thus  (\ref{csetcont}).
%
%It suffices that for $\tilde{t} \leq t^*$ that
%$$
%2C \sec\big(\pi/2 - \qzero + \qzero/2  + \qzero/8 + \qzero/8 \big) c_0 \leq \fr%ac{c_1}{4} \sin (1/8),$$
%or $c_0 \leq \frac{c_1}{8C} \frac{\sin(1/8)}{\sec\big(\pi/2 - \qzero + \qzero/2% + \qzero/8 + \qzero/8 \big)}$. This is ensured by Definition [RG$\cir$]. 
%This completes the proof of (\ref{csetcont}).
In order to obtain the statement of the lemma, 
let $P = \big\{ \xo = \bo{p_0},\bo{p_1},\ldots,\bo{p_r} = \yo \big\}$
satisfy $\bo{p_i} \in {\rm CRG}_{\qzero/2,K,\phi} \big( C_{\xo,\yo}\big)$ for $i \in \big\{ 1,\ldots,r-1 \big\}$. We must find 
$i \in \{ 1,\ldots, r \}$ for which $\vert\vert \bo{p_i} - \bo{p_{i-1}} \vert\vert \geq \frac{c_1}{4\pi}u$.

By the occurrence of $\mar \geq u n^{-1}$, the definition of the pertinent
pair
 $(\bo{x'},\bo{y'})$, and 
$\big\{ \bo{x'} = \xo \big\} \cap \big\{ \bo{y'} = \yo \big\}$, there exists $\bo{w} \in S^1$ for which
\begin{equation}\label{cainc}
W_{\bo{w},\frac{u}{2n}}\big( \bo{0} \big) \subseteq A_{\xo,\yo}
\end{equation}
and
\begin{equation}\label{crgem}
W_{\bo{w},\frac{u}{2n}}\big( \bo{0} \big) \cap \reg = \emptyset.
\end{equation}
By means of $\vert\vert \xo \vert\vert, \vert\vert \yo \vert\vert \geq c_1 n$ (which follows from
$\cir \cap B_{c_1 n} = \emptyset$), 
each of the balls $B_K(\xo)$ and $B_K(\yo)$ has an angular width of at most $\frac{\pi K}{2c_1 n}$ as viewed from $\bo{0}$, so that, by $u \geq C \log n$,
we infer from (\ref{cainc}) that
\begin{equation}\label{combb}
W_{\bo{w},\frac{u}{4n}}\big( \bo{0} \big) \cap  \Big( B_K \big(\xo \big) \cup
B_K \big(\yo \big)  \Big) = \emptyset.
\end{equation}
By  (\ref{csetcont}), (\ref{crgem}) and (\ref{combb}), we find that
$$
W_{\bo{w},\frac{u}{4n}}\big( \bo{0} \big) \cap {\rm CRG}_{\qzero/2,K,\phi}^{\xo,\yo} \big( C_{\xo,\yo} \big) = \emptyset.
$$
Note that $\axyo \setminus 
W_{\bo{w},\frac{u}{4n}}\big( \bo{0} \big)$
consists of two connected components, $C_1$ and $C_2$, each having the form of an angular sector rooted at $\bo{0}$, and with $\xo \in C_1$ and $\yo \in C_2$, say. 

Let 
$$
i = \inf \Big\{ j \in \{ 1,\ldots,r \}: \bo{p_j} \in C_2 \Big\}.
$$
Then $\bo{p_{i-1}} \in C_1$ and $\bo{p_i} \in C_2$. Thus,
$$
\argu \big(\bo{p_i}\big) - 
\argu \big(\bo{p_{i-1}}\big) \geq \frac{u}{4n}.
$$
However, the occurrence of {\rm SAT} entails that
$\cir \cap B_{c_1 n} = \emptyset$, so that 
 $\vert\vert \bo{p_{i-1}} \vert\vert,\vert\vert \bo{p_i} \vert\vert \geq
 c_1 n$. Hence,
$$
\vert\vert \bo{p_i} - \bo{p_{i-1}} \vert\vert \geq \frac{2}{\pi}
 \min \Big\{  
\vert\vert \bo{p_i} \vert\vert  ,
\vert\vert \bo{p_{i-1}} \vert\vert  \Big\}
 \ang \big(  \bo{p_i} ,  \bo{p_{i-1}} \big) \geq \frac{2}{\pi} c_1 n \frac{u}{4n} =
 \frac{c_1 u}{2 \pi},
$$
as required. \qed
\noindent{\bf Proof of Lemma \ref{lemexwaofp}.}
We introduce a procedure, to be denoted by  {\rm SEARCH}, that will aim to 
detect an appropriate
well-aligned outward-facing pair. The input of ${\rm SEARCH}$ consists of a
configuration $\omega \in \zoz$ and a direction $\bo{u} \in S^1$.
If ${\rm SEARCH}$ terminates successfully, then its output will consist of
a pair 
$(\bo{v},\bo{w}) \in \Z^2 \times \Z^2$, $\argu(\bo{v}) \leq \argu(\bo{u})
\leq \argu(\bo{w})$, which is a well-aligned outward-facing pair for the
outermost $\omega$-open circuit $\cir$ enclosing $\bo{0}$.
(We continue to write $\cir$ for this circuit in the following.)

To specify {\rm SEARCH}, we define two subprocedures. These are the  {\rm counterclockwise} and {\rm
  clockwise sweep}, $\sweepm$ and $\sweepp$. Each sweep has as input a
configuration $\omega \in \zoz$ and a specified element $\xo \in \reg$. 
For each sweep, the sweep will be well-defined on a certain subset of 
$(\omega,\bo{x}) \in \zoz \ltimes \reg$. (A sweep will be said to succeed
if it returns an output, and to fail otherwise.) In acting on  
$(\omega,\bo{x})$, the counterclockwise sweep $\sweepp$ will return an
element $\bo{y} \in \reg$ satisfying $\argu(\bo{y}) > \argu(\bo{x})$ if it
succeeds, whereas $\sweepm$  will return an
element $\bo{y} \in \reg$ for which  $\argu(\bo{y}) < \argu(\bo{x})$. 

We now specify the action of $\sweepp$ on $(\omega,\bo{x})$. Set 
$C = \club{\bo{x}}$.
The
boundary $\partial C$ is comprised of two semi-infinite line segments,
one of which attains the closest approach to $\bo{0}$ among points in
$\partial C$. Let $\ell = \ell(C)$ denote the planar line of which this line
segment forms a part. We write $\ell^*$ for the planar line in which the other
of the two semi-infinite line segments is contained. The plane with the
line $\ell$ removed  
is composed of two disjoint half-planes, one of which contains $\bo{0}$.
Let $H$ denote this half-plane. Let $\bo{x'}$ denote 
 the first element of $H^c \cap \reg$
encountered in the counterclockwise sense from $\bo{x}$, provided that such a
point exists. If $\bo{x'}$ does exist, then the action of $\sweepp$ on
$(\omega,\bo{x})$
is defined to succeed, with its output  $\sweepp \big(\omega,\bo{x} \big)$
being given by $\bo{x'}$. If, in addition, $\bo{x'} \in C$, the output is
declared to be ${\it good}$. If no such $\bo{x'}$ exists, then $\sweepp$
 is declared to fail on input $(\omega,\bo{x})$.

The action of the clockwise sweep $\sweepm$ on $(\omega,\bo{x})$
 is defined verbatim, with $C$ now being given by  $\clumb{\bo{x}}$, and
 with $\bo{x'}$ being sought in a clockwise search.

We now define the action of $\search$ on 
input $\big(\omega,\bo{u} \big)
\in \zoz \times S^1$. Similarly to the subprocedures, $\search$ will return
an output (and will be said to act successfully) only for a subset of
inputs $\big(\omega,\bo{u} \big)$. We will specify a finite sequence 
$\big\{ \bo{x_0}, \bo{x_1},\ldots \big\}$ of elements of $\reg$. As a
preliminary step, we take $\xo$ to be the first element of $\reg$
encountered in a clockwise sense that begins in the direction $\bo{u}$ (and
has $\bo{0}$ as its centre). If no such element exists (due to $\reg =
\emptyset$), then $\search$ is declared to fail. 
We then iteratively construct $\bo{x_1},\bo{x_2},\ldots$ by alternately
applying $\sweepp$  and $\sweepm$, with the input being $(\omega,\bo{x})$,
where $\bo{x}$ is the last element in the presently constructed sequence
$\big\{ \bo{x_0}, \bo{x_1}, \ldots \big\}$. 
This procedure stops at the end of the first sweep which is performed
successively and with a good output. In this event, $\search$ succeeds,
and returns as output the final two
elements in the presently constructed sequence $\big\{
\bo{x_0},\bo{x_1},\ldots \big\}$ (with these two elements ordered according
to increasing argument). If ever a sweep is performed that fails, then
$\search$ fails on the input  $\big(\omega,\bo{u} \big)$ in question.

We now note some relevant properties of these procedures. 

\noindent{\bf Claim 1:} 
if a sweep, of either type, in acting on an input 
 $\big(\omega,\bo{x} \big) \in \zoz \ltimes \reg$, does succeed, 
then the 
output (that we call $\bo{x'}$), is such that the pair 
$\big(\bo{x},\bo{x'} \big) \in \reg^2$ 
is outward-facing. \\
\noindent{\bf Claim 2:} if either sweep returns a good output $\bo{x'}$ in
acting on $\big(\omega,\bo{x} \big)$, then  
the pair $\big(\bo{x},\bo{x'} \big)$ is well-aligned. \\
\noindent{\bf Claim 3:} consider the sequence $\big( \bo{x_0}, \ldots, 
\bo{x_m} \big)$ eventually constructed in the action of $\search$ on input 
$\big(\omega,\bo{u} \big)$, as well as   the sequence 
of angular intervals 
$\big\{ \big( \argu(\bo{x_i}), \argu(\bo{x_{i+1}})  \big) : 0 \leq i \leq m
- 1 \big\}$. (Naturally, we actually order
$\argu(\bo{x_i}) <  \argu(\bo{x_{i+1}})$  in the
above list.) This sequence is strictly nested, and the increase from
one term to the next occurs alternately between the two endpoints. \\
\noindent{\bf Claim 4:} 
if $\omega \in \zoz$ satisfies the event
\begin{equation}\label{gdev}
\big\{ \globdis \leq \epsilon n/C_0 \big\} \cap \big\{ 
 \cir \cap B_{c_1 n} = \emptyset \big\} \cap \big\{ \cir \subseteq B_{\cctwo n}
 \big\}
 \cap \big\{ \mar \leq \epsilon \big\},
\end{equation}
then any sweep, acting on $\big(\omega,\bo{x} \big)$, will return an output,
for any choice of $\bo{x} \in \reg$. Moreover, the angular displacement
between input and output will be at most $2\epsilon$. 

Admitting these claims for now, we complete the argument. The claims imply that, if the input satisfies (\ref{gdev}), 
then procedure ${\rm SEARCH}$, in acting on an input $(\omega,\bo{u})$ with
$\bo{u} \in S^1$, will necessarily succeed, and will
return a well-aligned outward-facing pair of the form
$\big(\bo{x_m},\bo{x_{m+1}}\big)$, with
$\argu(\bo{x_m}) < \argu(\bo{x_{m+1}})$ (by reordering 
$\bo{x_m}$ and $\bo{x_{m+1}}$ if necessary), such that 
$\argu(\bo{x_{m+1}}) - \argu(\bo{x_m}) < 2 \epsilon$ and
$\bo{u} \in A_{\bo{x_m},\bo{x_{m+1}}}$. 
Indeed, by Claim 4, none of the sweeps made during ${\rm SEARCH}$ will fail. Nor can the sequence 
$\big\{ \bo{x_i}: i \geq 0 \big\}$ continue indefinitely: by Claim 3, each 
term in the sequence  
$\big\{ \bo{x_i}: i \geq 0 \big\}$ is distinct from all the earlier
ones, so that the sequence must be finite, since the set $\reg$ is finite. The
sequence finishes when the output of one of the sweeps is declared to be
good. At such a time, the final pair  $\big(\bo{x_m},\bo{x_{m+1}}\big)$ in the sequence is outward-facing and well-aligned, by Claims 1 and 2.
The property $\bo{u} \in A_{\bo{x_m},\bo{x_{m+1}}}$ follows by inductively
showing that $\bo{u} \in A_{\bo{x_i},\bo{x_{i+1}}}$ for each $0 \leq i \leq
m$.

It remains to verify the four claims. \\
\noindent{\bf Proof of Claim 1.} 
To see the claim for $\sweepp$, we retain the notation $C,\ell,H$ and $\bo{x'}$
from the definition of the sweep. (Note that this use of $\bo{x'}$ is
consistent with its use in the statement of the claim.) We have that 
$$
\reg \cap \Big( A_{\bo{x},\bo{x'}} \setminus \big( \{ \bo{x} \} \cup \{ \bo{x'} \} \big)  \Big) \subseteq H.
$$
However, $H \cap
A_{\bo{x},\bo{x'}} \subseteq T_{\bo{0},\bo{x},\bo{x'}}$, so that
$$
\reg \cap \Big( A_{\bo{x},\bo{x'}} \setminus \big( \{ \bo{x} \} \cup \{
\bo{x'} \} \big)  \Big) \subseteq   T_{\bo{0},\bo{x},\bo{x'}},
$$ 
and $\big( \bo{x}, \bo{x'} \big)$ is outward-facing. \\
\noindent{\bf Proof of Claim 2.} This is by construction. 
In the case of $\sweepp$, for example, we have that $\bo{x'} \in
\club{\bo{x}}$. \\
\noindent{\bf Proof of Claim 3.} We validate the claim in a generic instance. 
\begin{figure}\label{figwaotp}
\begin{center}
\includegraphics[width=0.3\textwidth]{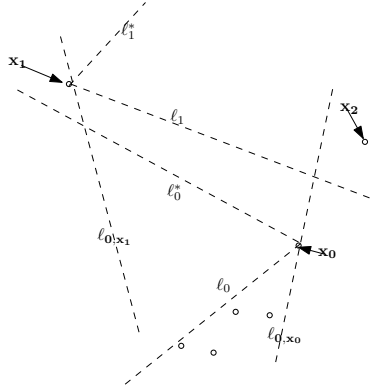} \\
\end{center}
\caption{Illustrating the proof of Claim 3. The circular marks represent elements of $\reg$. The circuit $\cir$ itself is not depicted.}
\end{figure}
Suppose that, in the action of
$\search$ on 
$\big(\omega,\bo{u} \big)$, we have that $m \geq 2$, so that at least the
element $\bo{x_2}$ is constructed. In the course of the action of $\search$,
$\sweepp$ acts on $(\omega,\xo)$ to return $\bo{x_1}$, and 
$\sweepm$ acts on $(\omega,\bo{x_1})$ to return $\bo{x_2}$. 
We use the notation $C,\ell,\ell^*$ and $H$ as in the definition of a sweep, the
subscripts $0$ and $1$ in what follows respectively referring to the action of
$\sweepp$ on $(\omega,\xo)$ and that of  
$\sweepm$ on $(\omega,\bo{x_1})$.
We will argue that $\argu(\bo{x_2}) < \argu(\bo{x_1})$, in this way,
showing that the clockwise sweep $\sweepm$ causes the angular interval $\big[
\argu(\bo{x_0}), \argu(\bo{x_1}) \big]$ to be extended on the left to
become  $\big[
\argu(\bo{x_2}), \argu(\bo{x_1}) \big]$. 

By the definition of the action of $\sweepm$ on $(\omega,\bo{x_1})$, we must
verify that
\begin{equation}\label{rghax}
 \reg \cap H_1^c \cap A_{\bo{x_0},\bo{x_1}} = \big\{ \bo{x_1} \big\}.
\end{equation} 
Indeed, in this case, the clockwise search that $\sweepm$ performs will begin
at $\bo{x_1}$ and will continue past $\bo{x_0}$ before locating $\bo{x_2}$.
To see (\ref{rghax}), note that 
\begin{equation}\label{rghaxx}
 \reg \cap A_{\bo{x_0},\bo{x_1}} \subseteq H_0 \cup  \big\{ \bo{x_0} \big\}
 \cup  \big\{ \bo{x_1} \big\},
\end{equation} 
by the fact that $\sweepp(\omega,\bo{x_0}) = \bo{x_1}$.
The construction of the sequence $\big\{ \bo{x_0},\bo{x_1},\ldots \big\}$
did not stop at $\bo{x_1}$. Therefore, the action of  $\sweepp$ on
$(\omega,\bo{x_0})$ was not declared to be good. We learn that $\bo{x_1}
\not\in C_0$. Since $\bo{x_1} \not\in H_0$ (because this vertex is the output
 $\sweepp(\omega,\bo{x_0})$), we see that $\bo{x_1}$ lies in the opposite
 half-plane as does $\bo{0}$ in the division provided by $\R^2 \setminus
 \ell_0^*$. Note that the region $C_1$ never meets the line $\ell_0^*$:
the line $\ell_1$, to which the line segment in the boundary of $C_1$ that is
closer to the origin belongs, meets the line $\ell_0^*$ that contains the line
segment in the boundary of $C_0$ that is further from the origin, at a
location counterclockwise to $\bo{x_1}$, because $\argu(\bo{x_1}) >
\argu(\xo)$. 

This property implies that $T \subseteq H_1$, where $T$ denotes 
the triangular region bounded by the lines $\ell_{\bo{0},\bo{x_0}}$,
$\ell_{\bo{0},\bo{x_1}}$ and $\ell_0^*$. (Recall that $\ell_{\bo{x},\bo{y}}$ denotes
the planar line containing $\bo{x},\bo{y} \in \R^2$.)

Noting that $\big( H_0 \cap A_{\bo{x_0},\bo{x_1}} \big) \cup \big\{
\bo{x_0} \big\} \subseteq T$, we find from (\ref{rghaxx}) that
$\reg \cap A_{\bo{x_0},\bo{x_1}} \subseteq T \cup \big\{ \bo{x_1 }
\big\}$. From $T \subseteq H_1$, we arrive at (\ref{rghax}). This completes
the argument for Claim 3. \\ 
\noindent{\bf Proof of Claim 4.} 
By (\ref{czercond}), for all 
$\bo{x} \in \partial \wulff$,
\begin{equation}\label{eqnparbox}
\partial \wulff \cap W_{\bo{x},c_0} \subseteq C_{\pi/2 - 3\qzero}^F \big( \bo{x} \big)
\cup  C_{\pi/2 - 3\qzero}^B \big( \bo{x} \big). 
\end{equation}
Recall that $\tcir = n \partial \wulff$ 
denotes the dilation of $\partial \wulff$ attaining 
$\globdis$. Set $0 < \cwinf : = \inf \big\{ \vert\vert \bo{x} \vert\vert: \bo{x} \in \partial \wulff \big\}$,
and $\cwsup : = \sup \big\{ \vert\vert \bo{x} \vert\vert: \bo{x} \in \partial \wulff \big\} < \infty$. 
%The hypotheses on $\cir$ imply that there exist $0 < c < C < \infty$ such that $cn \leq t_0 \leq Cn$. 
By $\cir \subseteq \tcir + B_{\globdis}$, it follows from Lemma \ref{lemlu} 
%and $0 < \cwinf : = \inf \big\{ \vert\vert \bo{x} \vert\vert: \bo{x} \in \partial \wulff \big\}$,
%$\cwsup : = \sup \big\{ \vert\vert \bo{x} \vert\vert: \bo{x} \in \partial \wulff \big\} < \infty$ 
that
\begin{equation}\label{gamtoc}
\cir \subseteq  
\big( n + c' \epsilon n \big) \wulff \setminus  \big( n - c' \epsilon n \big) \wulff,
\end{equation}
where $c' < \frac{2 \cwinf \qzero}{\pi^2 \cwsup}$ may 
be ensured by choosing $C_0$ high.
We now argue that, for $\bo{x} \in \tcir$,
\begin{equation}\label{eqntzeroc}
 \Big( 
 \big( n + c' \epsilon n \big) \wulff \setminus 
 \big( n - c' \epsilon n \big) \wulff \Big)
 \cap A_{\argu(\bo{x})+\epsilon,\argu(\bo{x}) + 2\epsilon}
 \subseteq C^F_{\pi/2 - 5\qzero/2}\big( \bo{x} \big).
\end{equation}
To this end, let $\bo{y}$ belong to the left-hand-side, and let 
$\bo{v} \in \ell_{\bo{0},\bo{y}}^+ \cap n \partial \wulff$. By (\ref{eqnparbox}), 
$\ang \big( \bo{v} - \bo{x} , \perpu{\bo{x}} \big) \geq 3 \qzero$.
To establish (\ref{eqntzeroc}), we must show that 
$\ang \big( \bo{y} - \bo{x} , \bo{v} - \bo{x} \big) \leq \qzero/2$.
Note that $\bo{y} \in t \partial \wulff$ with $\big\vert t - n \big\vert \leq c' \epsilon n$, so that $\vert\vert \bo{y} - \bo{v} \vert\vert \leq \cwsup c' \epsilon n$.
We have that $\vert\vert \bo{v} - \bo{x} \vert\vert \geq \vert\vert \bo{x} \vert\vert \sin \ang \big( \bo{x}, \bo{v} \big) \geq 2 \pi^{-1} \cwinf  n \epsilon$, from $\vert \vert \bo{x} \vert\vert \geq \cwinf  n$ and $\ang \big( \bo{x} , \bo{v} \big) \geq \epsilon$. We find then that 
$\sin \ang \big( \bo{y} - \bo{x} , \bo{v} - \bo{x} \big) \leq \frac{\pi \cwsup c'}{2 \cwinf}$, so that  
$\ang \big( \bo{y} - \bo{x} , \bo{v} - \bo{x} \big) \leq 
 \frac{\pi^2 \cwsup c'}{4\cwinf } \leq \qzero/2$, confirming (\ref{eqntzeroc}). 

Let $\bo{y} \in \big( n + c' \epsilon n \big) \wulff \setminus \big( n - c'
\epsilon n \big) \wulff$. 
Set $\bo{x}$ to be the point of intersection 
$\tilde\Gamma_0 \cap \ell_{\bo{0},\bo{y}}^+$. 
From 
$\vert\vert \bo{x} \vert\vert \geq  \cwinf n$ 
%$\vert\vert \bo{x} \vert\vert \geq  c \cwinf n$ 
and $\vert\vert \bo{y} - \bo{x} \vert\vert \leq c' \cwsup \epsilon n$, by a short argument that we  omit,
\begin{equation}\label{doubinc}
C_{\pi/2 - 5\qzero/2}^F\big(\bo{x}\big) 
\cap \big( W_{\bo{x},\epsilon}^+ \big)^c \subseteq
C_{\pi/2 - 2\qzero}^F\big(\bo{y}\big) \qquad \textrm{and} \qquad
C_{\pi/2 - 5\qzero/2}^B\big(\bo{x}\big) 
\cap \big( W_{\bo{x},\epsilon}^- \big)^c \subseteq
C_{\pi/2 - 2\qzero}^B\big(\bo{y}\big),    
\end{equation}
provided that $c'$ is chosen small enough (so that $C_0$ is chosen to be high): in fact, 
$$
c' < \frac{\pi \cwinf}{8\cwsup} \frac{\tan(5\qzero/2) - \tan(2\qzero)}{\tan(2\qzero)\tan(5\qzero/2)}
$$ 
suffices.

By (\ref{gamtoc}), (\ref{eqntzeroc}) and (\ref{doubinc}) we find that, for any $\bo{y} \in \Gamma_0$,
\begin{equation}\label{vacincl}
\cir \cap A_{\argu(\bo{y}) + \epsilon,\argu(\bo{y}) + 2\epsilon}
 \subseteq  C_{\pi/2 - 2\qzero}^F \big( \bo{y} \big).
\end{equation}
Without loss of generality,
consider the clockwise sweep acting on  $\big(\omega,\bo{y} \big)$. Note
that $A_{\argu(\bo{y}) + \epsilon,\argu(\bo{y}) + 2\epsilon} \cap \reg
\not= \emptyset$, because $\mar \leq \epsilon$.
Note that the sweep will necessarily finish before or on
reaching any element of the set $\reg \cap \club{\bo{y}}$. We have then that 
(\ref{vacincl}) implies that the sweep will return on output 
$\bo{y'} \in \reg$ satisfying $\argu(\bo{y}) \leq \argu(\bo{y'})
\leq \argu(\bo{y}) + 2\epsilon$, as we sought to show. \qed
\end{section}

\bibliographystyle{plain}

\bibliography{mlrbib}

\end{document}